\documentclass[10pt]{article}

\usepackage{latexsym}
\usepackage{amssymb}
\usepackage{amscd}
\usepackage{amsmath}
\usepackage{amsthm}
\usepackage[usenames,dvipsnames]{color}
\usepackage{graphicx}
\usepackage{wrapfig}
\usepackage{tikz}
\usetikzlibrary{shapes.geometric}
\usepackage{url}
\usepackage{fancyvrb}

\usepackage[mathscr]{euscript}

\usepackage{hyperref}
\hypersetup{
 colorlinks,
 citecolor=Violet,
 linkcolor=Blue,
 urlcolor=Blue}

\def\pf{\noindent{\bf Proof:}~ }
\def\eop{\hfill\rule{2.0mm}{2.0mm}}

\def\bC{{{\mathbb C}}}
\def\bR{{{\mathbb R}}}

\def\geq{\geqslant}
\def\leq{\leqslant}

\newcommand{\goto}{\rightarrow}

\newcommand{\argmin}{\operatornamewithlimits{\arg\hskip -0.02in \min}}

\newcommand{\beq}{\begin{equation*}}
\newcommand{\eeq}{\end{equation*}}
\newcommand{\beqn}{\begin{equation}}
\newcommand{\eeqn}{\end{equation}}
\newcommand{\bea}{\begin{eqnarray}}
\newcommand{\eea}{\end{eqnarray}}

\newcommand{\gap}{\vspace{0.1in}}
\newcommand{\demo}[1]{\pf}
\newcommand{\bx}{{\bf x}}

\newcommand{\cV}{\mathcal{V}}
\newcommand{\cF}{\mathcal{F}}
\newcommand{\cE}{\mathcal{E}}
\newcommand{\euS}{\mathscr{S}}

\newtheorem{proposition}{Proposition}[section]

\newtheorem{theorem}[proposition]{Theorem}
\newtheorem{corollary}[proposition]{Corollary}
\newtheorem{conjecture}[proposition]{Conjecture}
\theoremstyle{definition}
\newtheorem{definition}[proposition]{Definition}
\theoremstyle{remark}
\newtheorem{remark}[proposition]{Remark}

\begin{document}
\title{{\sc Numerical Methods for Biomembranes}: \\conforming subdivision methods
versus non-conforming PL methods}

\date{December 28, 2018 \\ Revised: March 11, 2020}

\author{
Jingmin Chen\thanks{
Citigroup Global Markets Inc., 390 Greenwich Street, New York, NY 10013, U.S.A.. Email: \href{mailto:jingmchen@gmail.com}{jingmchen@gmail.com}.
}
\and
Thomas Yu
\thanks{
Department of Mathematics, Drexel University. Email: \href{mailto:yut@drexel.edu}{yut@drexel.edu}.
He was supported in part by the National Science Foundation grants DMS 0512673 and DMS 0915068.
}
\and
Patrick Brogan\thanks{Department of Mathematics, Drexel University. Email: \href{mailto:pbrogan12@gmail.com}{pbrogan12@gmail.com}. He was supported in part by the Office of the
Provost and the Steinbright Career Development Center of Drexel University.}
\and
Robert Kusner\thanks{Department of Mathematics, University of Massachusetts at Amherst. Email:
\href{mailto:kusner@math.umass.edu}{kusner@math.umass.edu}. He was supported in part by the National Science Foundation grants PHY 1607611, DMS 1439786 and DMS 1440140.}
\and
Yilin Yang\thanks{Center for Computational Engineering, M.I.T., Email: \href{mailto:yiliny@mit.edu}{yiliny@mit.edu}.}
\and
Andrew Zigerelli
\thanks{
Department of Electrical and Computer Engineering, University of Pittsburgh. Email: \href{mailto:anz37@pitt.edu}{anz37@pitt.edu}.
He was supported in part by a 2013 Goldwater scholarship during his study at Drexel University.
}
%
}

\makeatletter \@addtoreset{equation}{section} \makeatother
\maketitle

\centerline{\bf Abstract:}
The Canham-Helfrich-Evans models of biomembranes consist of a family of
geometric constrained variational problems.
In this article, we compare two classes of numerical methods
  for these variational problems based on piecewise linear (PL) and subdivision surfaces (SS).
Since SS methods are based on spline approximation and
 can be viewed as higher order versions of PL methods, one may expect that the only difference between
 the two methods is in the accuracy order. In this paper,
 we prove that
a numerical method based on minimizing any one of the `PL Willmore energies' proposed in the literature
 would fail to converge to a solution of the continuous problem, whereas a method based on
  minimization
of the bona fide Willmore energy, well-defined for SS but not PL surfaces, succeeds.
Motivated by this analysis,
we propose also a regularization method for the PL method  based on techniques from conformal geometry.
We address a number of implementation issues crucial for the efficiency of our solver. A software package called {\sc Wmincon}
accompanies this article, provides parallel implementations of all the relevant geometric functionals. When combined with a
standard constrained optimization solver, the
geometric variational problems can then be solved numerically. To this end, we realize that some of the available
optimization algorithms/solvers are capable of preserving symmetry, while others manage to break symmetry; we
explore the consequences of this observation.

\vspace{.2in} \noindent {\bf Acknowledgments.} TY is indebted to Tom Duchamp and Aaron Yip for extensive
discussions and many of their insightful remarks. We also thank Tim Mitchell, Michael Overton, Justin Smith, and Shawn Walker for help.
TY was partially supported
by NSF grants DMS 0915068 and DMS 1115915. RK was supported in part by the Aspen Center For Physics (funded by NSF-PHY 1607611),
ICERM (funded by NSF-DMS 1439786), and MSRI (funded by NSF-DMS 1440140.)

\vspace{.2in}

\section{Introduction}
Lipid bilayers are arguably the most elementary
and indispensable structural components of biological membranes which form the
boundary of all cells.
It is
known since the seminal work of Canham \cite{Canham:Elastic}, Helfrich \cite{Helfrich:Elastic} and
Evans \cite{Evans1974923} in the
70’s that bending elasticity, induced by curvature, plays the key role in driving the
geometric configurations of such membranes.

\begin{figure}
\begin{center}
\includegraphics[height=3.2cm]{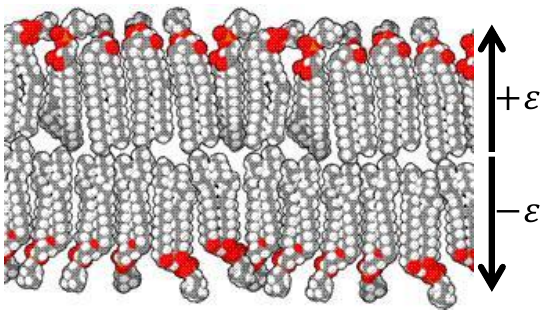}
\end{center}
\vspace{-20pt}
\label{fig:offset}
\caption{The `offset surfaces' of the pivotal surface of a lipid bilayer}
\end{figure}

The so-called spontaneous curvature model of Helfrich
suggests that a biomembrane surface $S$ configures itself to minimize
$\int_S H^2 dA$ subject to the area, volume and area difference (related to the bilayer characteristics) constraints,
i.e. $S$ solves the variational \textbf{Helfrich problem}
\begin{equation}
\begin{aligned}
 \min_S W(S) :=\int_S H^2 \, dA
\text{ s.t. }
\left\{
  \begin{array}{ll}
    \mbox{(i)} & A(S) := \int_S 1 \; dA = A_0, \\
    \mbox{(ii)} & V(S) :=
\frac{1}{3} \int_S [x {\hat{\mathbf{i}}} + y \hat{\mathbf{j}} + z \hat{\mathbf{k}}] \cdot \hat{\mathbf{n}} \: dA =
V_0, \\
    \mbox{(iii)} & M(S) :=-\int_S H \; dA = M_0.
  \end{array}
\right.
\end{aligned}
\label{eq:Helfrich}
\end{equation}
Here $H=(\kappa_1+\kappa_2)/2$ is the mean curvature.
In (ii), $V(S)$ is the enclosed volume, expressed here as a surface integral of $S$ via
the divergence theorem.
The connection of (iii)
to bilayer area difference comes from the
relation $-\int_S H dA = \lim_{\varepsilon \goto 0} \frac{1}{4\varepsilon}({\rm area} (S_{+\varepsilon})- {\rm area} (S_{-\varepsilon}))$,
where $S_{+\varepsilon}$ and $S_{-\varepsilon}$ are the `$\varepsilon$-offset surfaces',\footnote{\label{footnote:sign}We assume that the normal of any closed orientable surface points outward. In particular, it means $H<0$ for a sphere.}
 and
that the thickness of the lipid bilayer, $2 \varepsilon$, is negligible compared to the size of the vesicle; see Figure~\ref{fig:offset}.
The constraint values $A_0$, $V_0$ and $M_0$ are determined by physical conditions (e.g. temperature, concentration).
$W(S)$ is called the \textbf{Willmore energy} of the surface $S$.
When the area-difference constraint (iii) is omitted, the variational problem
 is referred to as the \textbf{Canham problem}.
When even the volume constraint (ii) is omitted, there is essential no constraint as $W$ is scale-invariant; in this
 case the area constraint (i)
only fixes the scale, and we refer to the variational problem as the \textbf{Willmore  problem}.

Due to the scale-invariance of the Willmore energy,
the solution, up to homothety, of any
of the Willmore, Canham or Helfrich problems depends only on %
the
\emph{reduced volume} and \emph{reduced total mean curvature} defined by:
\bea
v_0:= V_0/[(4\pi/3) (A_0/4\pi)^{3/2}], \quad m_0 := M_0/[4\pi(A_0/4\pi)^{1/2}].
\eea
This terminology is used by a group of biophysicists who have done
many computational and physical experiments exploring the shapes of
 phospholipid vesicles, and we shall follow it. Note that $v_0$ is essentially what
a geometer would call the \emph{isoperimetric ratio}. By the isoperimetric inequality, we have
$v_0 \in (0,1]$ and $v_0=1$ is uniquely realized by a round sphere.

It is observed experimentally that no topological change occurs in any accessible time-scale, so we
aim to solve  any of the Helfrich, Canham or Willmore problems when $S$ is assumed to be an orientable closed
surface with a fixed genus $g$. Spherical ($g=0$) vesicles are the most common among naturally
occurring biomembranes, although higher genus ones have been synthesized in the laboratory
\cite{MichaletBensimon:Genus2,LIPOWSKY:Genus2,Seifert:Config}.
The Canham, Helfrich and related models explain the large
variety of shapes observed in even a closed vesicle with a spherical topology
\cite{LIPOWSKY:Nature,Seifert:Config,Lim:Stomatocyte–discocyte–echinocyte,Seifert:Starfish}.

Several numerical treatments of these models have been proposed in the
literature: \cite{Kusner:Brakke,Brakke:evolver}, \cite{Feng2006394,ChenGrundelYu:Sphere},
\cite{Bonito:Biomembranes,Dziuk:Review}, \cite{Qian-2004,Qian-2005,Qian-2006}, \cite{Schumacher:WillmoreFlow}.
Among these,
the methods in \cite{Kusner:Brakke,Brakke:evolver}
were used extensively by biophysicists 
to study real lipid bilayer membranes. While the key ingredients of these algorithms are
implemented in Brakke's well-known
Surface Evolver software \cite{Brakke:evolver}, the overall algorithms were not completely analyzed
 by
the geometers who invented them \cite{Kusner:Brakke,Brakke:evolver,Francis1997} and even less so by
the biophysicists who used them
\cite{LIPOWSKY:Genus2,MichaletBensimon:Genus2,Seifert:Config,Seifert:Starfish}.
As such, there are little understanding of these methods, and
the computational results claimed in the extensive biophysics literature are difficult to reproduce.
Moreover, there is no systematic
comparison of this method with the later ones, at both a theoretical or computational level.

These numerical methods continue to be used extensively in the study of phenomena in biomembranes, see, e.g.
\cite{ZHAO2017164,Baumgart2003,KSM:Confinement1,KSM:Confinement2,BGB:DynamicsBiomembranes} and references therein.
Similar geometric variational problems show up in other scientific areas. A notable example is found in the
quasi-local mass problem of general relativity, in which
maximizers of the Hawking mass -- defined similarly as the Willmore energy -- are sought.

The goal of this paper is to clarify and refine some of these numerical methods, and establish some theoretical understandings of them.
Before we proceed, we mention further related work on
discrete minimal surfaces and discrete elasticae (the 1-D counterpart to Willmore surface),
see, e.g., \cite{Dziuk:DiscretePlateauI,Dziuk:DiscretePlateauII,Schumacher-2019,Bruckstein-2001,scholtes:2019} and the references therein.

\subsection{PL and SS}
 \label{sec:PL_SS}
A standard approach to represent surfaces of arbitrary topology is to use the piecewise linear (PL) approach.
A PL surface can be specified by a mesh $\mathcal{M}=(\mathcal{V},\mathcal{F})$ where
$\mathcal{V} \in \bR^{\#V \times 3}$ records the 3-D coordinates of the vertices of the control mesh, $\#V$ denotes the total
number of vertices, and $\mathcal{F} \in \mathcal{I}^{\#F \times 3}$ is a list of triplets of indices from
$\mathcal{I} := \{1,\ldots,\#V\}$ which
records the vertices of each of the $\#F$ triangle faces in the mesh $\mathcal{M}$. We assume that the PL surfaces realized
by the mesh are closed and orientable.
The orientation can be conveniently encoded in a consistent ordering of the vertices in the face list $\mathcal{F}$.

In a numerical method, $\cF$ is usually fixed and $\cV$ varies. This fits the framework of
our variational problems well, as
fixing $\cF$ also fixes the genus of the surface, and varying $\cV$ means we find the embedding of $\cF$ --
viewed as an abstract simplicial complex
-- that optimizes the Willmore energy under the corresponding constraint(s).

A closed, oriented PL surface has a well-defined area $A$ and enclosed volume $V$, but no
 classically defined normals or mean curvatures, hence it also does not have
 a classically defined total mean curvature $M$ or Willmore energy $W$.
As such, any numerical method for the Willmore, Canham, or Helfrich problems
based on approximating the solution surfaces by PL surfaces may be classified as a \textbf{nonconforming method} in FEM parlance.

\begin{figure}[h]
\centerline{
\begin{tabular}{cc}
\includegraphics[width=2.3in]{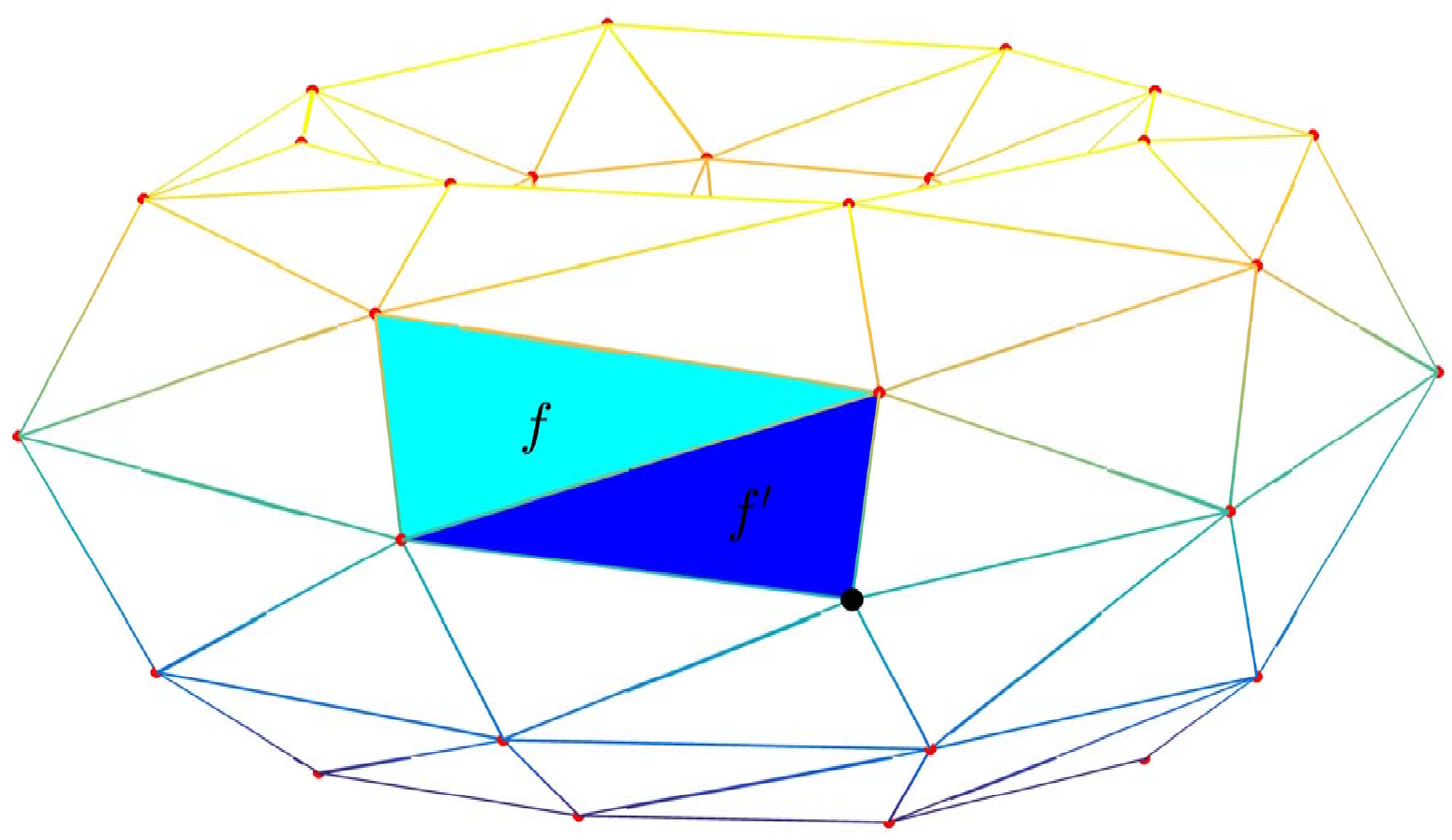} & \includegraphics[width=2.3in]{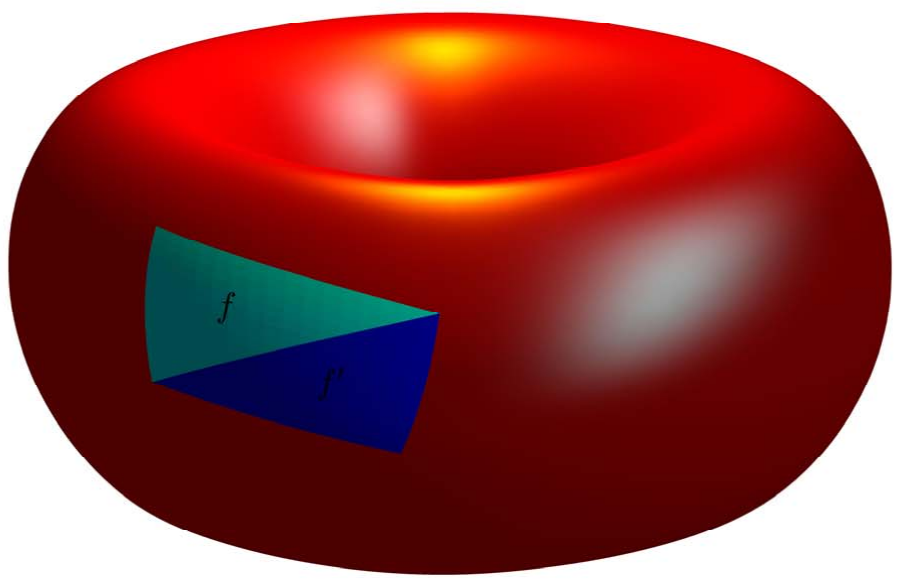}
 \end{tabular}}
\caption{(a) A control mesh $\mathcal{M}$; $f$ is a regular face, $f'$ is an irregular face (b) The Loop subdivision surface
corresponding to $M$; in light blue: the regular patch associated to $f$, in deep blue: the irregular patch associated to $f'$}
\label{fig:CMesh}
\end{figure}

A subdivision surface (SS) is specified by the same data $\mathcal{M}=(\mathcal{V},\mathcal{F})$, except that the
associated surface has enough regularity for a well-defined total mean curvature $M$ or Willmore energy $W$.
We shall primarily be using the Loop and C2g0 subdivision surfaces introduced in \cite{LoopScheme}
and \cite{ChenGrundelYu:Sphere}, respectively.
See Figure~\ref{fig:CMesh}(a) for a genus 0 control mesh, and Figure~\ref{fig:CMesh}(b)
for the corresponding Loop subdivision surface.
For each face $f$ in $\mathcal{M}$, there is a corresponding surface patch; see Figure~\ref{fig:CMesh}(b).

The C2g0 scheme handles only genus 0 and 1 surfaces using control meshes with only valence 3 and 6 vertices; the resulting SS
are $C^2$ everywhere, so there is no question about the well-definedness
of $M$ and $W$. The Loop subdivision scheme handles surface of arbitrary genus and control meshes with arbitrary valences, but the resulting SS
are not $C^2$ everywhere.
They are however regular enough to have well-defined $M$ and $W$. 
(The subdivision functions are in $W^{2,2} \cap C^1$ when expressed in characteristic coordinates; see Section~\ref{sec:positive}.)
Therefore, any numerical method for the Willmore, Canham, or Helfrich problems
based on approximating the solution surfaces by SS may be classified as a \textbf{conforming method}.

\subsection{Contributions of this paper}
This paper contributes to the numerical study of biomembranes in the following ways:
\begin{itemize}
\item[(I)] In Section~\ref{sec:Functionals_PLSS}, we give an exposition of numerical methods based on both PL and SS.
In the PL case, we connect and compare several different ideas developed in the applied geometry literature
\cite{Kusner:Brakke,Brakke:evolver,Bobenko05discretewillmore,Meyer2003,Desbrun:1999:IFI:311535.311576}. 
In the SS case, we explain how all the functionals and their gradients with respect to control vertex coordinates can be efficiently
computed based on a \emph{precomputation of basis functions}.
In both cases, we develop parallel implementations and explore new examples. A Matlab based toolbox named
{\sc Wmincon}, with CUDA and C++ implementations of all key functionals, is
available online, at

\centerline{
\url{http://www.math.drexel.edu/~tyu/Software/Wmincon},}

for reproducing the computational results. These algorithmic developments allow us to use the SS method to attack
many instances of the Willmore, Canham and Helfrich problems that would otherwise be too slow, if not impossible, to solve on existing computers.
\item[(II)] In Section~\ref{sec:analysis}, we present an argument explaining why a conforming method would work. This result relies on the
existence theory of Willmore minimizers pioneered by Simon \cite{Simon:Willmore}. 
In contrast,
we show that a naive minimization of several PL Willmore energies would fail to solve the Willmore problem.
The analysis elucidates the difference between conforming and non-conforming methods.
\item[(III)] In Section~\ref{sec:Regularization}, based on a well-founded principle (Section~\ref{sec:Continuous}) exploiting the
uniformization theorem and the theory of harmonic maps,
we propose a regularization of the PL methods based on penalization by harmonic energy.
Unlike the unregularized PL methods which are doomed to fail, the regularized method
appears to yield  solutions converging to those of the continuous problems.
\item[(IV)] In Section~\ref{sec:Further}, we make the observation that certain optimization algorithms are capable of preserving
symmetry, while others are capable of breaking symmetries. We carry out a number of experiments comparing different
optimization algorithms in conjunction with our discretization methods. The experiments reveal subtle analytic
properties of the optimization problems arising from the SS methods. We formulate a number of conjectures.
\end{itemize}

The kind of conforming and non-conforming methods we study in this article are those in the spirit of `minimizing a discretization', i.e.
the methods under study first discretize the variational problem, in either a conforming or non-conforming way, followed by solving the
resulting finite-dimensional optimization problem.
There are methods, such as those in \cite{Bonito:Biomembranes,Schumacher:WillmoreFlow}, that are
in the spirit of `discretizing a minimization'. These methods first
consider a minimization process in the continuous setting, akin to a gradient flow, followed by strategies to discretize the flow. This last step can also be done in a conforming or non-conforming way.
These methods are all based on explicit representations of surfaces;
there are also methods based on implicit representations, such as the phase field methods of Du et al \cite{Qian-2004,Qian-2005,Qian-2006}.

\section{Numerical Methods Based on PL and SS Functionals} \label{sec:Functionals_PLSS}
Recall that either a PL or SS is specified by a control mesh $\mathcal{M}=(\mathcal{V},\mathcal{F})$.
In our numerical method, we assume that $\mathcal{F}$ is fixed and $\mathcal{V}$
varies.
For most $\mathcal{V}$, an immersed surface, denoted by $S[\cV]$, is defined.
The numerical methods considered here approximate the Helfrich problem
\eqref{eq:Helfrich} by a finite-dimensional analog:
\begin{equation} \label{eq:Helfrich_Subdivision}
\begin{aligned}
 \min_{\mathcal{V}}  W(\cV)
\text{ s.t. }
\left\{
  \begin{array}{ll}
    \mbox{(i)} & A(\cV)  = A_0 \\
    \mbox{(ii)} & V(\cV)  = V_0 \\
    \mbox{(iii)} & M(\cV)  = M_0
  \end{array}
\right..
\end{aligned}
\end{equation}
The numerical methods for the Canham and Willmore problems are similar: simply drop the corresponding constraints.
Already mentioned in Section~\ref{sec:PL_SS}, the PL and SS methods have the following
features and relative pros and cons:
\begin{itemize}
\item For SS (based on any regular enough scheme, such as Loop and C2g0), all four functionals are the
{\it exact}, well-defined, values of the $W$, $A$, $V$ and $M$ of the corresponding subdivision surface.
Their computations, however,
have to be performed based on numerical integration. 
\item For PL, $W$ and $M$ are not well-defined for the corresponding PL surface. We will therefore
 replace $W(\cV)$ and $M(\cV)$ in \eqref{eq:Helfrich_Subdivision} by a certain \emph{consistent discretization},
 to be reviewed below.
 These PL Willmore and total mean curvature energies are relatively
 simple to implement and no numerical integration is required.
\end{itemize}

The materials in this section are mostly not new, some of them are actually quite old. The intention is to unify them at one place in
 order to prepare us for the later sections.

At first glance, one may expect that the PL method is simply less accurate than the SS methods, i.e. a PL method would converge but
at a lower rate compared to a SS method.
Our analysis in Section~\ref{sec:analysis} falsifies this speculation. 
A bulk of this section discusses the definition, properties and computation of these functionals and their gradients.
Efficient computation of these functionals and their gradients are
 necessary for the numerical solution of the Helfrich, Canham and Willmore problems
 using a standard nonlinear optimization solver; see Section~\ref{sec:details}.

\subsection{$W$, $A$, $V$, $M$ for PL surfaces} \label{sec:PL}
The area $A$ and enclosed volume $V$ are of course part of the biomembrane problems.
Their gradients are not only needed for our optimization
solver but also are connected
to the way $M$ and $W$ are defined and computed.
For these reasons, we derive them for the convenience of the readers.

We aim to clarify some not so well-documented details in the literature, such as the
 sign issue of discrete mean curvature, which is irrelevant for
$W$ but crucial for $M$, and the choice of local areas, which is irrelevant for $M$ but impact the behavior for $W$.
 Another goal is to elucidate the connections of a number of
 different discrete mean curvature operators and Willmore
and total mean curvature energies.

\subsubsection{$A$ and $V$}
The area $A$ and enclosing volume $V$ of a closed oriented PL surface can be computed as:
\begin{align} \label{eq:AV}
\begin{split}
A = \frac{1}{2} \sum_{f \in \mathcal{F}} \| (\mathcal{V}_{f_2}-\mathcal{V}_{f_1}) \times (\mathcal{V}_{f_3}-\mathcal{V}_{f_1}) \|,
\quad
V = \frac{1}{6} \sum_{f \in \mathcal{F}} \det([\mathcal{V}_{f_1}, \mathcal{V}_{f_2}, \mathcal{V}_{f_3}]).
\end{split}
\end{align}


For any smooth functional
$F: {\rm domain}(F) \overset{{\rm open}}{\subset} \bR^{\#V \times 3} \cong \prod_{v \in \mathcal{I}} \bR^3 \goto \bR$,
we denote by $\nabla_v F$ ($\in \bR^3$) its gradient
with respect to the coordinates of the vertex indexed by $v\in \mathcal{I}$.
The volume gradient can be expressed as
\begin{align} \label{eq:VolumeGradient}
\nabla_v V & = \frac{1}{6} \sum_{i} \nabla_v \det([\mathcal{V}_{v}, \mathcal{V}_{w_i}, \mathcal{V}_{w_{i+1}}])
= \frac{1}{6} \sum_{i} \mathcal{V}_{w_i} \times \mathcal{V}_{w_{i+1}} \in \bR^3,
\end{align}
where $w_1, \ldots, w_{i}, w_{i+1}, \ldots \in \mathcal{I}$ is
a counterclockwise enumeration (viewed from the outside) of the vertices connected to $v$ (a.k.a. the `1-ring' of $v$).

Note that the gradient formula can be used to show that the formula for $V$ is invariant under rigid motions when the PL surface
is \emph{closed} and \emph{consistently oriented}. (Observe that $\sum_v \langle \nabla_v V, \mathbf{a} \rangle = 0$ for any
constant vector $\mathbf{a}$; the proof relies on both assumptions.)
The volume gradient itself is invariant under translation and equivariant with respect to rotation.
The latter is obvious from the formula; the former is obvious also as $V$ is translation-invariant,
but it helps to see it directly from the formula:
$\sum_i (\mathcal{V}_{w_i} + \mathbf{a}) \times (\mathcal{V}_{w_{i+1}} + \mathbf{a})
= \sum_{i} \mathcal{V}_{w_i} \times \mathcal{V}_{w_{i+1}} + \sum_i \mathcal{V}_{w_i} \times \mathbf{a}
+ \mathbf{a} \times \sum_i  \mathcal{V}_{w_{i+1}}  + \sum_i \mathbf{a} \times \mathbf{a}
= \sum_{i} \mathcal{V}_{w_i} \times \mathcal{V}_{w_{i+1}}.$
Note that the two sums in the middle cancel only because $w_1, \ldots,w_{\rm val(v)}$ form
a closed-loop.

\begin{remark} \label{remark:EffArea}
The vector $\nabla_v V$, in turn, has another geometric interpretation: if the ``base of the pyramid
around $v$'' is coplanar, i.e. the vertices indexed by $w_1, \ldots, w_{{\rm val}(v)}$ lie on the same plane and form a polygon, then
$\nabla_v V$ is a vector orthogonal to the plane and its length is one-third the area of the polygon. (By translation-invariant,
we can assume that the polygon is centered at the origin.)
Clearly, it is independent of the coordinates of vertex $v$.
In general, $3\|\nabla_v V\|$ can be used to define a notion
of the ``area of a non-planar polygon."
This so-called `effective area' is used to define one of the discrete
mean curvatures.
\end{remark}

Next, we have the following derivation for the area gradient:
\begin{align} \label{eq:cotangent}
\begin{split}
\nabla_v A & = \frac{1}{2} \sum_{i} \nabla_v \| (\mathcal{V}_{w_i} - \mathcal{V}_v) \times (\mathcal{V}_{w_{i+1}} - \mathcal{V}_v) \|
=
\frac{1}{2} \sum_{i} \frac{\overbrace{(\mathcal{V}_{w_i} - \mathcal{V}_v)}^{:=p_i} \times
\overbrace{(\mathcal{V}_{w_{i+1}} - \mathcal{V}_v)}^{:=p_{i+1}}}
{\| (\mathcal{V}_{w_i} - \mathcal{V}_v) \times (\mathcal{V}_{w_{i+1}} - \mathcal{V}_v)\|}
\times \overbrace{(\mathcal{V}_{w_{i+1}} - \mathcal{V}_{w_{i}})}^{=p_{i+1}-p_i} \\
&= \frac{1}{2} \sum_{i} \frac{ -\big((p_{i+1}-p_i) \cdot p_{i+1} \big)
p_i + \big((p_{i+1}-p_i) \cdot p_i \big) p_{i+1} }{\| p_i \times p_{i+1} \|} \\
&= -\frac{1}{2} \sum_{i}
\frac{(p_{i+1}-p_i) \cdot p_{i+1}}{\| (p_{i+1}-p_i) \times p_{i+1} \|}  p_i
+
\frac{(p_{i}-p_{i+1}) \cdot p_{i}}{\| (p_{i}-p_{i+1}) \times p_{i} \|}  p_{i+1}
\\
&=
-\frac{1}{2} \sum_{i}
(\cot \angle v w_{i+1} w_i) p_i + (\cot \angle v w_i w_{i+1}) p_{i+1} =
\frac{1}{2} \sum_{i}
(\cot \alpha_i + \cot \beta_i)
(\mathcal{V}_v - \mathcal{V}_{w_i}),
\end{split}
\end{align}
where $\alpha_i$ and $\beta_i$ are the angles opposite the edge $v w_i$ in the two incident triangles.
In above, the second equality can be seen from the chain rule, in which an intermediate map is of the form
$C(\mathbf{x}) = (\mathbf{a} - \mathbf{x}) \times (\mathbf{b}-\mathbf{x})$, which can be simplified to $\mathbf{a} \times \mathbf{b} + (\mathbf{b}-\mathbf{a}) \times \mathbf{x}$
and hence has a constant derivative expressible by a cross product.
The third equality follows from the vector triple product formula
$(\mathbf{a} \times \mathbf{b})\times \mathbf{c} = -(\mathbf{c}\cdot \mathbf{b})\mathbf{a}  + (\mathbf{c} \cdot \mathbf{a})\mathbf{b}$.
The fifth equality follows from
$\mathbf{a} \cdot \mathbf{b}/\| \mathbf{a} \times \mathbf{b} \| =\cot(\mbox{angle between $\mathbf{a}$ and $\mathbf{b}$})$.

Equation~\ref{eq:cotangent} is connected to the well-known cotangent formula for the Laplace-Beltrami operator;
see Remark~\ref{remark:cotangent}.

The following comment will be found useful when computing discrete mean curvature.
\begin{remark} \label{remark:orientation}
\noindent
t is clear that $A$, and hence also $\nabla A$, has nothing do with the global orientation of the PL surface;
in particular, they are well-defined even for a non-orientable PL surface. The direction of $\nabla_v A$ tells `which way
the PL surface is poking' at the vertex $v$.
However, the enclosing volume $V$ requires
the surface to be both closed and orientable, and in this case the formula for $V$ in \eqref{eq:AV} would only give the
enclosing volume if all the faces are oriented in a counter-clockwise fashion when viewed from the outside. In
particular, reversing the orientation of all the faces would flip the sign of $V$ and reverse the direction of each $\nabla_v V$.
\end{remark}

\subsubsection{$W$ and $M$}

We review 5 discrete Willmore and two discrete total mean curvature energy functionals for PL surfaces which we learn from
\cite{Kusner:Brakke,Sullivan2008,Meyer2003,Bobenko05aconformal,Brakke:evolverManual}. We label them as
\begin{align*}
W_{\rm Centroid}, \;W_{\rm Voronoi}, \;W_{\rm EffArea}, \;W_{\rm NormalCur}, \;W_{\rm Bobenko}, \;
M_{\rm Cotan}, \; M_{\rm Steiner}
\end{align*}
in this paper and in the {\sc Wmincon} package.

Recall that for any smooth orientable surface $S$ with continuous unit
normals denoted by $\mathbf{n}(x)$, $x \in S$, we have
\begin{align} \label{eq:AreaVariation}
\frac{d}{dt}\big|_{t=0} {\rm Area} (S_t) = -2 \int_S h(x) H(x) dA, \quad \forall \; h: S \goto \bR,
\end{align}
where
$S_t := \{x+t h(x) \mathbf{n}(x): x \in S \}$, and $H$ is the mean curvature defined relative to the choice of the normals $\mathbf{n}$.
The above functionals, except $W_{\rm Bobenko}$/$M_{\rm Edge}$, can be derived
based on defining
$$
\mbox{`discrete normals' $\mathbf{n}: \mathcal{I}\goto S^2$}, \;\;\; \mbox{`discrete mean curvatures' $H:\mathcal{I} \goto \bR$,}
\;\;\; \mbox{and} \;\;\;
\mbox{`local areas' $a: \mathcal{I} \goto \bR^+$}
$$
at the vertices of a PL surface, indexed by $\mathcal{I} = \{1,\ldots,\#V \}$, that satisfy a discrete analog of \eqref{eq:AreaVariation},
namely,
\begin{align} \label{eq:H_Area}
\lim_{t\goto 0} \frac{1}{t} \left[ A\left(\mathcal{V} + t [h(v) \mathbf{n}(v)]_{v \in \mathcal{I}} \right) - A(\mathcal{V}) \right]
=
-2 \sum_{v \in \mathcal{I}} h(v) H(v) a(v), \quad \forall \;\; h: \mathcal{I} \goto \bR.
\end{align}
The left-hand side is the directional derivative of $A$ at $\mathcal{V}$ in the direction
$[h(v) \mathbf{n}(v)]_{v \in \mathcal{I}}$, which equals
$$
\sum_v \big\langle \nabla_v A(\mathcal{V}), h(v) \mathbf{n}(v) \big\rangle_{\bR^3}.
$$
In order for it
to equal the right-hand side of \eqref{eq:H_Area} for all scalar field $h$, it is necessary and sufficient,
by setting $h(v') = \delta_{v,v'}$, for
$\mathbf{n}(v)$, $H(v)$ and $a(v)$ to satisfy
\begin{align*}
\nabla_v A \cdot \mathbf{n}(v) = -2 H(v) \, a(v), \;\;\; \forall v.
\end{align*}
Once $a(v)$ is assigned, then
$\mathbf{n}(v)$ and $H(v)$
can be chosen so that the
\emph{mean curvature vector} is
\begin{align} \label{eq:MeanCurVec}
\mathbf{H}(v) := H(v) \mathbf{n}(v) = -\frac{\nabla_v A}{2a(v)}.
\end{align}
This only defines $H(v)$ and $\mathbf{n}(v)$ up to a sign;
the appropriate sign must be determined from the \emph{global orientation} of the PL surface;
a natural way is to choose $\mathbf{n}$ so that
$$
\langle \mathbf{n}(v), \nabla_v V \rangle > 0;
$$ recall Remark~\ref{remark:orientation} and Footnote~\ref{footnote:sign}. The sign of $H(v)$ can
then be determined accordingly.

The local areas used in the various schemes are summarized in the following table.

\gap
\centerline{
\begin{tabular}{|l|l|}
  \hline
  Scheme & local area $a(v)$ \\
  \hline
  {\rm Centroid} \cite{Kusner:Brakke} &  ${\rm Area} ({\rm star}(v))/3 =: a_{\rm centroid}(v)$ \\
  {\rm Voronoi} \cite{Meyer2003} &  ${\rm Area}( \mbox{Voronoi cell around } v)$ \\
  {\rm EffAreaCur} \cite{Sullivan2008},\cite[Page 223]{Brakke:evolverManual} & $\| \nabla_v V \|$ \\
  {\rm NormalCur} \cite[Page 223]{Brakke:evolverManual} &  $|\langle \nabla_v V, \nabla_v A \rangle|/ \| \nabla_v A\|$ \\
  \hline
\end{tabular}}

\gap
Recall Remark~\ref{remark:EffArea} for the effective area. The local area used in `normal curvature' is
the length of the projection of $\nabla_v V$ onto the direction of $\nabla_v A$. The rationale for the
use of these local areas are
discussed in \cite{BYY:Enumath17,Brakke:evolverManual}; see also Section~\ref{sec:W_PL}.

The discrete Willmore energy is then defined as
\begin{align}
\label{eq:WM_natural}
 W_{{\rm Centroid}/{\rm Voronoi}/{\rm EffArea}/{\rm NormalCur}} = \sum_v H(v)^2 a(v),
\end{align}
where the four choices of $a(v)$ in the table above correspond to the four discrete $W$-energies.
Note that the sign of $H(v)$ is irrelevant to the definition of $W$.

Similarly, a discrete total mean curvature functionals can be defined\footnote{The label `Cotan'
may not be ideal, but it reflects the fact that it is based on the cotangent formula \eqref{eq:cotangent} for $\nabla A$.}
as
\begin{align}
\label{eq:WM_natural}
 M_{\rm Cotan} = \sum_v H(v) a(v) = \frac{1}{2} \sum_v {\rm sign}(H(v))  \| \nabla_v A \|.
\end{align}
Unlike $W$, the sign of $H(v)$, dependent on orientation, is crucial, but the choice of local area is irrelevant.

An alternative discrete total mean
curvature, based on Steiner's polynomial, is defined by
\begin{align}
M_{\rm Steiner} := \sum_{e} {\rm length}(e) \,\theta(e),
\end{align}
where $\theta(e) \in (-\pi, \pi)$ is the signed angle between the normals to the adjacent faces at $e$;
see \cite[Section 4.4]{Sullivan2008}, \cite[Figure 6]{Bobenko05aconformal}, \cite[Page 227]{Brakke:evolverManual}.

\begin{remark} \label{remark:cotangent}
Besides the area-variation characterization \eqref{eq:AreaVariation}, we
also have the characterization of mean curvature based on the Laplace-Beltrami operator:
\begin{align} \label{eq:DiscreteLaplace}
\Delta_S \mathbf{X} (x) = 2 \mathbf{H}(x), \quad x \in S,
\end{align}
where $\mathbf{X}: S \goto \bR^3$ is the position function of the surface $S$.
It is just a matter of taste to derive a discrete mean curvature based on a discrete Laplace operator or a discrete area variation.
For our purpose here, we choose the latter simply because we need the area variation anyway for our solver.
In fact, by combining \eqref{eq:cotangent}, \eqref{eq:MeanCurVec} and \eqref{eq:DiscreteLaplace}
one can retrieve the cotangent formula for the discrete Laplace-Beltrami operator.
For yet another connection of the cotangent formula with area and Dirichlet energy, see Section~\ref{sec:PLConformalEnergy}.
Also, see \cite{Wardetzky2008} for the convergence properties of the cotangent formula.
\end{remark}

Bobenko's Willmore energy is based on a rather different philosophy: it is designed to satisfy an exact M\"obius invariant
property and measures a `degree of sphericity' \cite[Proposition 2]{Bobenko05aconformal}.
It is defined as
$$W_{\rm Bobenko} := \frac{1}{2} \sum_v W(v) + 4\pi(1-g),$$ where $W(v) = \sum_{e \ni v} \beta(e) - 2\pi$ and $\beta(e)$ is
an angle formed by the circumscribed
circles of the two triangles sharing the edge $e$ \cite[Definition 1]{Bobenko05aconformal}.\footnote{Since $W_{\rm Bobenko}$,
like all other discrete Willmore energies here, intends to approximate $\iint H^2 \,dA$,
whereas Bobenko's definition of $W$ in \cite{Bobenko05aconformal}
intends to be a discrete analog of $\iint H^2-K \,dA$, the two differ by $4\pi(1-g)$
according to the Gauss-Bonnet theorem. Of course, there is no difference in the genus $g=1$ case.}

A discrete Willmore energy $W_{\rm PL}$ or total mean curvature $M_{\rm PL}$ should have
a consistency property in the sense
that $W_{\rm PL}(\mathcal{M}^n) \goto W(S)$ and $M_{\rm PL}(\mathcal{M}^n) \goto M(S)$ for
 any sequence of PL surfaces $\mathcal{M}^n$ converging to a smooth surface $S$ in an appropriate
sense.
 $W_{\rm Bobenko}$  is known to be consistent with the continuous Willmore energy only
in a very restrictive sense \cite{Bobenko05discretewillmore}.
From our preliminary analysis, the other PL Willmore energies are better behaved in terms of consistency.
We shall report on these in a separate report.
Although such a result is not directly needed in this article,
we believe that it will be necessary for the analysis of
the PL method proposed in Section~\ref{sec:Regularization}.


\subsection{$W$, $A$, $V$, $M$ for subdivision surfaces} \label{sec:SS}
The two specific subdivision schemes used in our solver are the Loop and C2g0  schemes.
Here, we present the details of Loop's scheme \cite{LoopScheme}; the paper
\cite{ChenGrundelYu:Sphere} contains similar details for the C2g0 scheme.
Our presentation will be brief,
but contains the necessary
 implementation details when read in conjunction with the paper \cite{Stam:Loop} by Stam.

\subsubsection{Subdivision surfaces} \label{sec:SSIntro}
Following the subdivision surface literature, a
vertex is called {\it ordinary} if it has valence 6, otherwise it is called an \textit{extraordinary} vertex.
We assume that extraordinary vertices in $\mathcal{M}$ are isolated, i.e. no two extraordinary vertices can be neighbor
of each other. If $\mathcal{M}$ lacks this property, one can simply
apply a mid-point subdivision to $\mathcal{M}$ to resurrect that.

For each triangle face $f$ in $\mathcal{M}$,
we call $f$ a {\it regular} face if all its three bounding vertices are ordinary,
otherwise, under our
assumption, exactly one of the three vertices is extraordinary and we call $f$ an
{\it irregular} face. The corresponding surface patches will be simply referred to
as \textit{regular} and \textit{irregular patches}. See again Figure~\ref{fig:CMesh}.

\begin{figure}[h]
\centerline{
\begin{tabular}{c}
\includegraphics[width=1.55in]{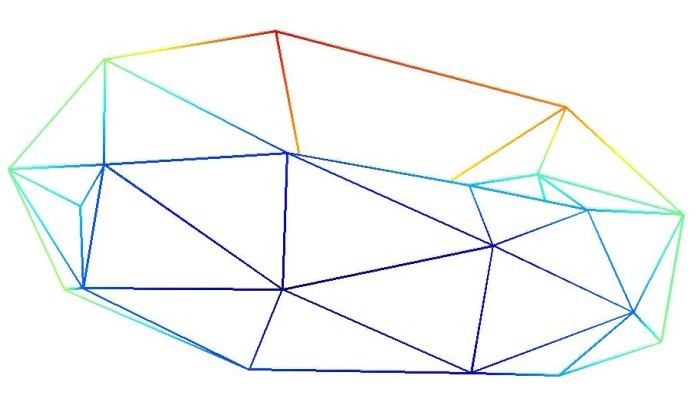}  \includegraphics[width=1.5in]{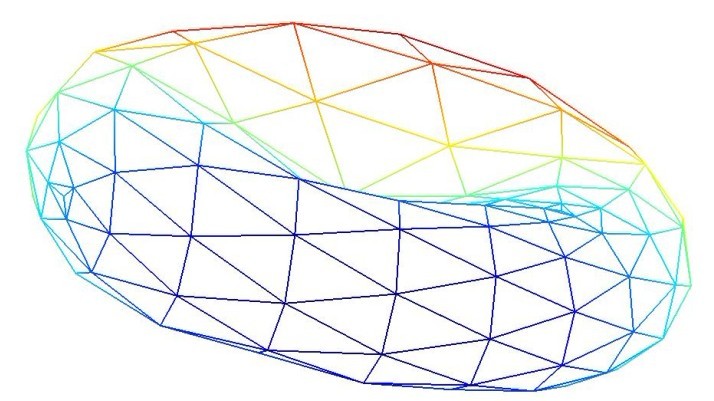}
\includegraphics[width=1.6in]{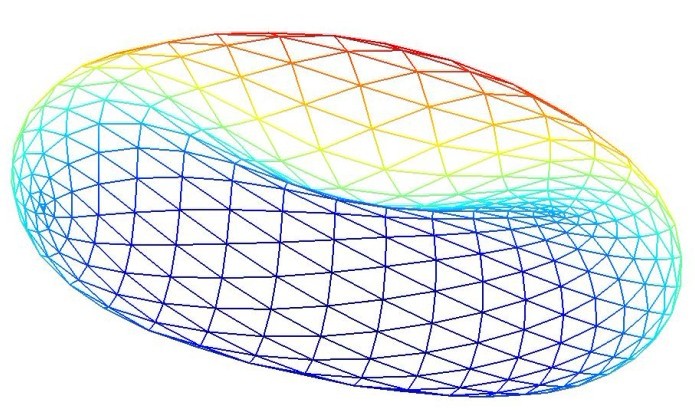}  \includegraphics[width=1.5in]{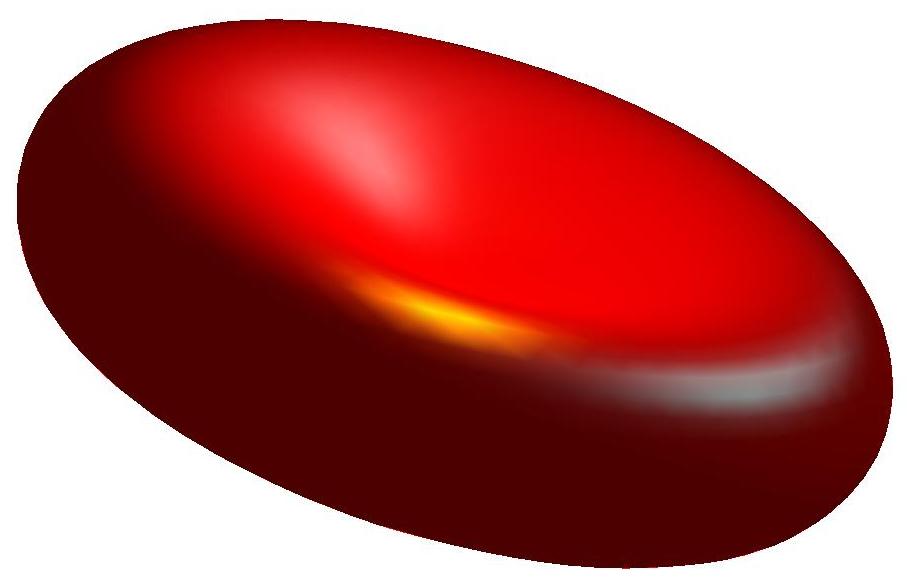}
 \end{tabular}}
\caption{Iteration of subdivision steps}
\label{fig:RedBloodCell}
\end{figure}

\begin{figure}[h]
\centerline{
\begin{tabular}{c}
\includegraphics[height=1.5in]{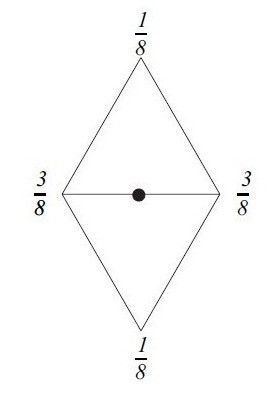} (a) \includegraphics[height=1.5in]{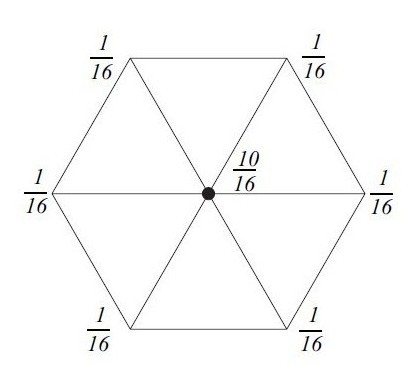} (b)
\includegraphics[height=1.5in]{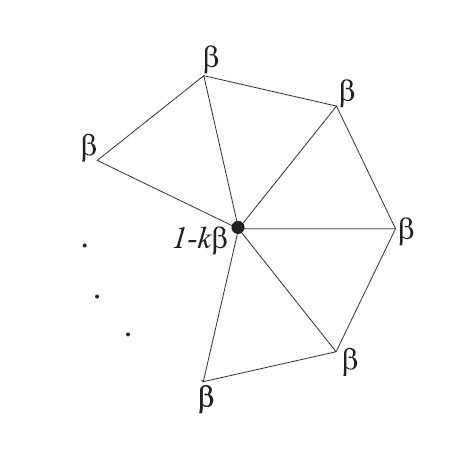} (c) \\
\end{tabular}}
\caption{(a) edge rule; (b) vertex rule for an ordinary vertex; (c) vertex rule for an extraordinary vertex, where $k$ is the valence and $\beta = \frac{1}{k}(\frac{5}{8} - (\frac{3}{8}+ \frac{1}{4}\cos(\frac{2\pi}{k}))^2)$. (In the text, we use `$N$' to denote
the valence of an extraordinary vertex, so as to be consistent with Stam's paper.)}
\label{fig:SubdivisionRules}
\end{figure}

Although the parametric description is more important for us, it is useful 
to recall the popular algorithmic
description of a Loop
surface as the limit of an iteration of subdivision steps, see Figure~\ref{fig:RedBloodCell};
the subdivision step can be succinctly described by the diagrams in Figure~\ref{fig:SubdivisionRules}.

\begin{figure}[h]
\centerline{
\begin{tabular}{ccccc}
\includegraphics[width=1.3in]{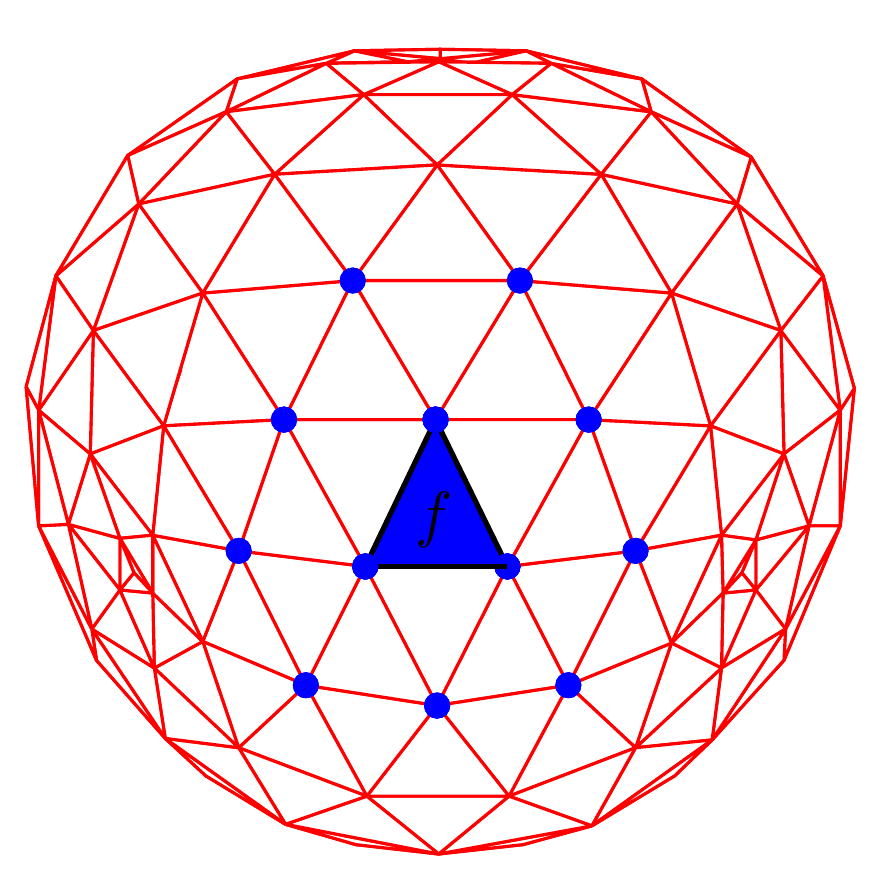} & \hspace{-.2in}
\includegraphics[width=1.3in]{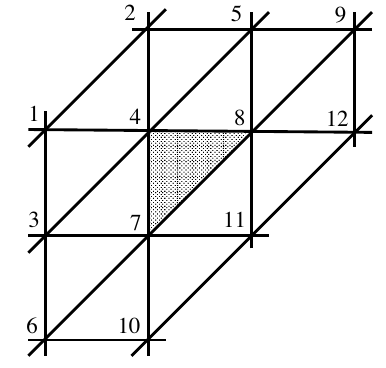} & \hspace{-.2in}
\includegraphics[width=1.3in]{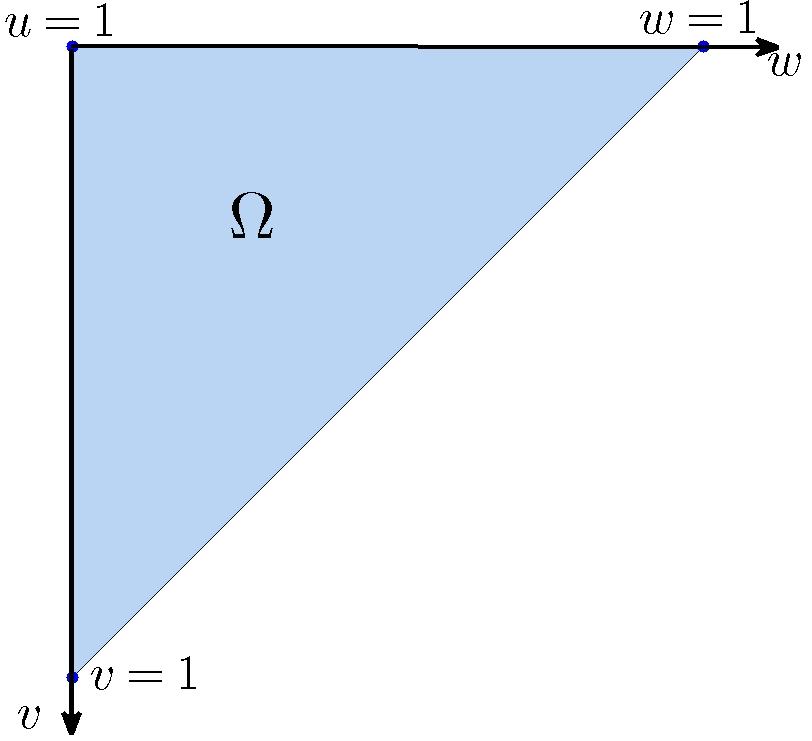} & \hspace{-.2in}
\includegraphics[width=1.3in]{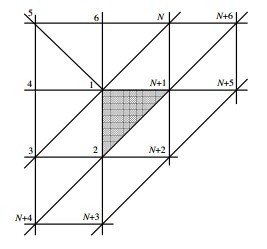} & \hspace{-.2in}
\includegraphics[width=1.3in]{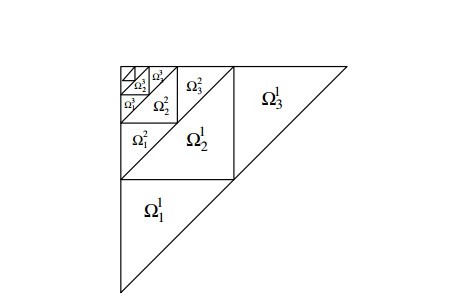}
\\
(a) &\hspace{-.2in} (b) &\hspace{-.2in} (c) &\hspace{-.2in} (d) &\hspace{-.2in} (e)
 \end{tabular}}
\caption{(a) A regular face and its neighboring vertices (b) ordering of the twelve control vertices
(c) the parameter domain $\Omega$. Note: for any point $(v,w) \in \Omega$, its barycentric
coordinates w.r.t. to the bounding vertices of $\Omega$, listed in the order $(0,0)$,
$(1,0)$, $(0,1)$, are simply $(1-v-w,v,w)$.
(d) Ordering of the $N+6$ control vertices around an irregular face (e) Partition of the parameter domain $\Omega$
}
\label{fig:OrdinaryPatch}
\end{figure}

\noindent
{\bf Parametrization of a regular patch.} 
The surface patch associated with a regular face $f$ can be parameterized by
a linear combination of $12$ polynomials with coefficients $\{ \mathbf{c}_{f,i}\}_{i=1}^{12}$
being the coordinates of the vertices in $f$ and their immediate neighbors ordered
as in Figure~\ref{fig:OrdinaryPatch}(b).\footnote{When $f$ is
close to an extraordinary vertex, it is possible that some of these 12 control vertices coalesce.
In such a degenerate situation, one simply repeats the
coalescing vertices when using the formula \eqref{eq:OrdinaryPatch}.}(see Figure~\ref{fig:OrdinaryPatch})
\bea \label{eq:OrdinaryPatch}
\bR^3 \leftarrow \Omega :\mathbf{s}_f(v,w) = \sum_{i=1}^{12} \mathbf{c}_{f,i}^T \mathbf{b}_i(v,w)
\eea
where $\Omega:= \{(v,w): v \in [0,1] \mbox{ and } w \in [0, 1-v]\}$,
and $\mathbf{b}_1(v,w), \ldots, \mathbf{b}_{12}(v,w)$ are the twelve degree 4 polynomials
as shown in \cite[Page 10-11]{Stam:Loop}. (We do not
copy these polynomials from Stam's paper, but mention that
they come from the so-called $M_{222}$ box-spline.)
For notational convenience, we organize
each $\mathbf{c}_{f,i}$ as a row vector of length $3$, write
$\mathbf{c}_f := [\mathbf{c}_{f,1}; \ldots; \mathbf{c}_{f,12}] \in \bR^{12 \times 3}$,
and define  $\mathbf{b} := \mathbf{b}^{6} := [\mathbf{b}_1;\ldots; \mathbf{b}_{12}]$ as a column vector of functions
of length 12. Then
\eqref{eq:OrdinaryPatch} simplifies to
\bea \label{eq:OrdinaryPatchC}
s_f = \mathbf{c}_f^T \mathbf{b}^{6}.
\eea

\noindent
{\bf Parametrization of an irregular patch.} 
The surface patch associated with an irregular face $f$ admits a parametrization
$\mathbf{s}_f: \Omega \goto \bR^3$ which is controlled by $N+6$ control vertices around the face $f$,
where
$N$ is the valence of the extraordinary vertex of $f$. Following Stam's convention,
these $N+6$ vertices are ordered as in Figure~\ref{fig:OrdinaryPatch}(d).
Like the regular case, $\mathbf{s}_f$ is \emph{linearly} related to the control vertices, so
\bea \label{eq:IrrParDifficulty}
\mathbf{s}_f(v,w) = \sum_{i=1}^{N+6} \mathbf{c}_{f,i}^T \mathbf{b}^{N}_i(v,w)
\eea
for some basis functions $\mathbf{b}^{N}_i$, $i=1,\ldots,N+6$,
implicitly defined by the subdivision process. Unlike the regular ($N=6$)
case, none of these basis functions is
a single polynomial anymore. Instead, it is an infinite piecewise polynomial,
with pieces being the (recursively defined) sub-triangles $\Omega^j_k$, $j=1,2,\ldots$, $k=1,2,3$,
as shown
 in Figure~\ref{fig:OrdinaryPatch}(e). Note that in this figure
 the origin $(0,0)$ corresponds to the extraordinary vertex of $f$.

\begin{figure}[h]
\centerline{
\begin{tabular}{ccccc}
\includegraphics[width=1.2in]{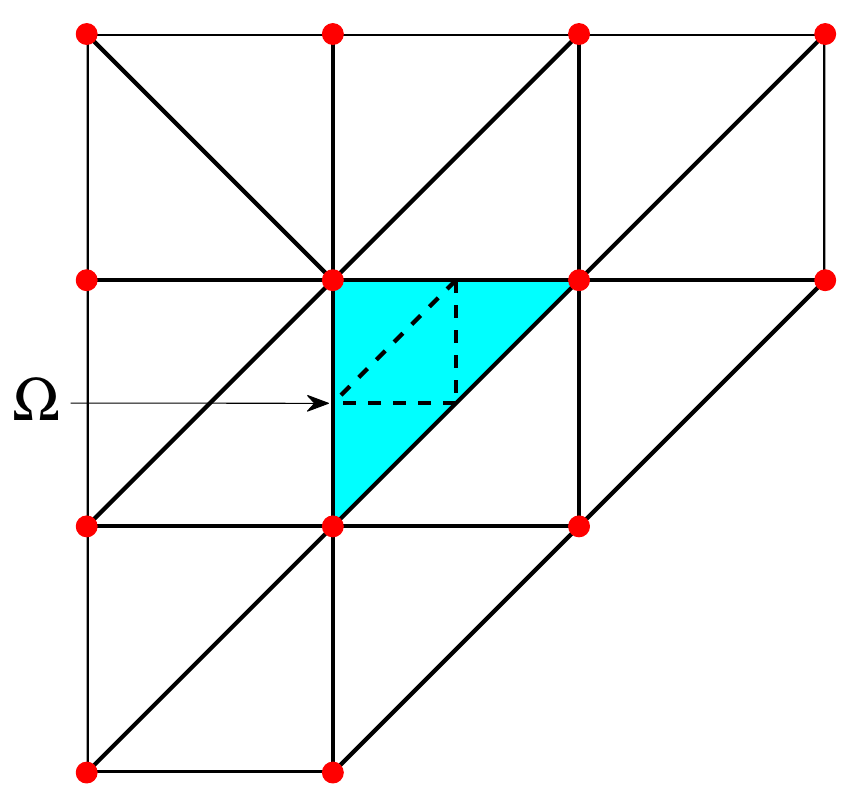} &
\includegraphics[width=1.2in]{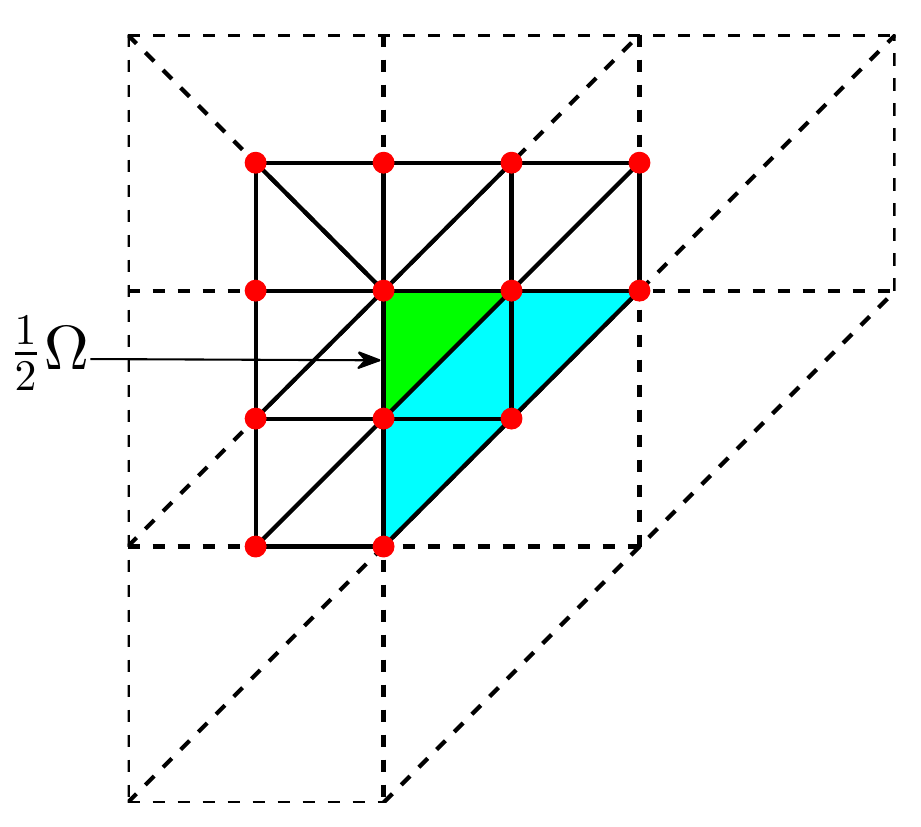} &
\includegraphics[width=1.2in]{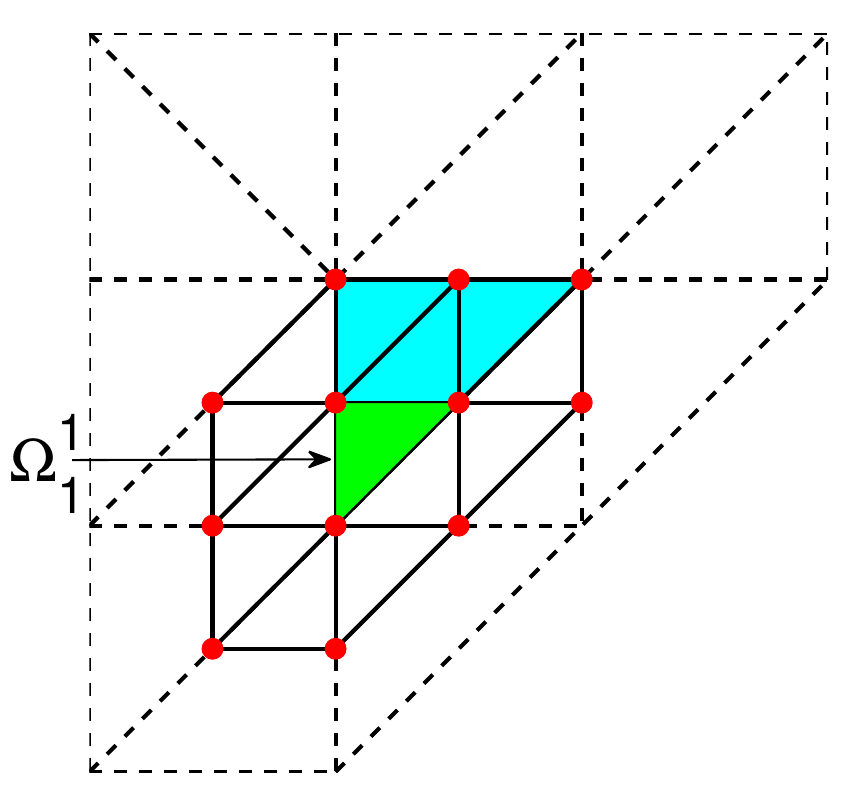} &
\includegraphics[width=1.2in]{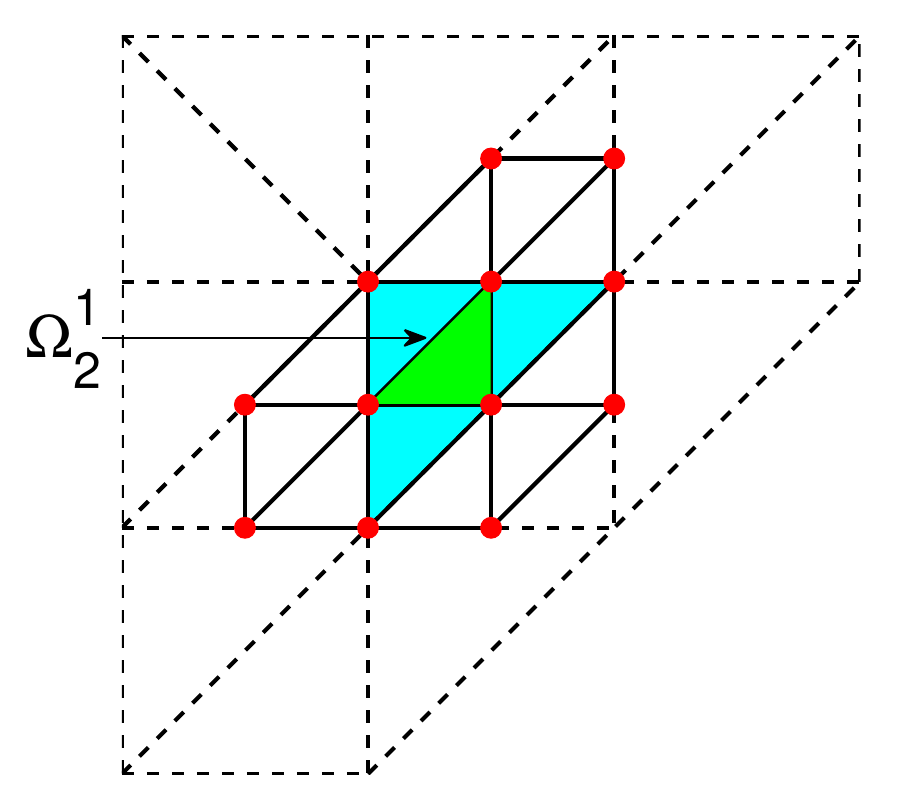} &
\includegraphics[width=1.2in]{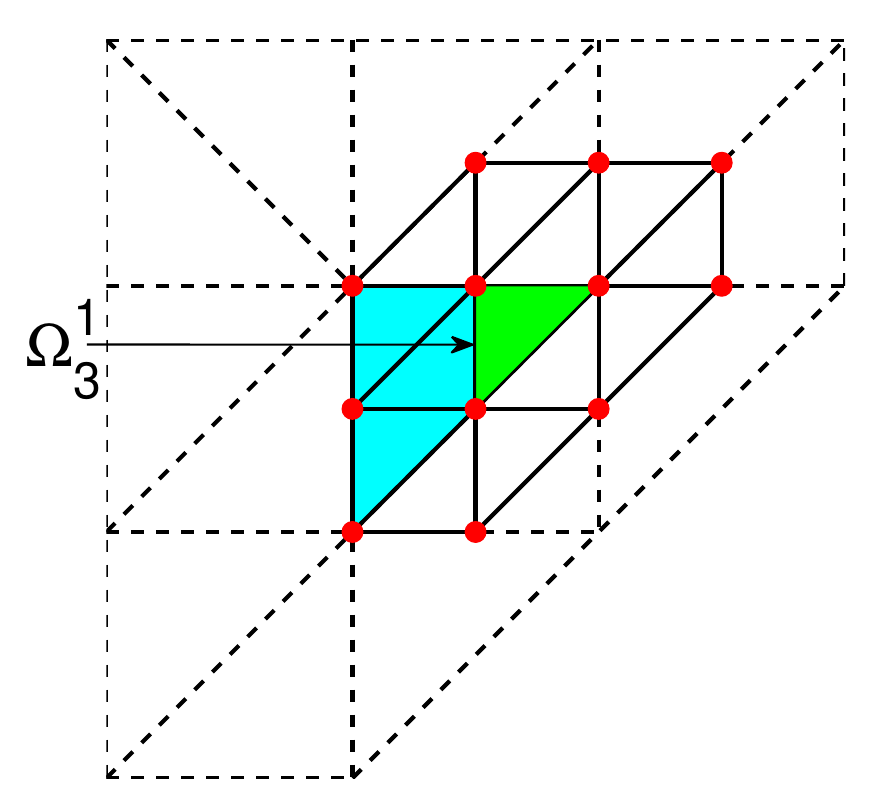} \\
(a) & (b) & (c) & (d) & (e)
 \end{tabular}}
\caption{(a) the $N+6$ control vertices around the extraordinary face $f$ that
determine $\phi_i: \Omega \goto \bR$
(b) the $N+6$ control vertices at the next subdivision level that determine
$\phi_i|_{\frac{1}{2} \Omega}$
(c)-(e) the $12$ control vertices at the next subdivision level
that determine $\phi_i|_{\Omega^1_k}$, $k=1,2,3$}
\label{fig:InfLayer}
\end{figure}

The parametrization \eqref{eq:IrrParDifficulty} is tricky to compute; and this is where
Stam's idea \cite{Stam:Loop,Stam:1998:EEC:280814.280945} comes in.
As in the regular case, write $\mathbf{c}_{f} := [\mathbf{c}_{f,1};\ldots;\mathbf{c}_{f,N+6}] \in
\bR^{(N+6)\times 3}$, and $\mathbf{b}^{N}:=[\mathbf{b}^{N}_1;\ldots;\mathbf{b}^{N}_{N+6}]$.
In a nutshell, Stam's method transforms the control data $\{\mathbf{c}_{f,i}\}$ into
`eigen- control data' $\widehat{\mathbf{c}}_{f} = V^{-1}\mathbf{c}_{f}$, so
\bea \label{eq:EigenXform}
\mathbf{s}_f = \mathbf{c}_f^T \mathbf{b}^{N} = \widehat{\mathbf{c}}_{f}^T
\underbrace{V^T \mathbf{b}^{N}}_{=: \Phi}
=\widehat{\mathbf{c}}_{f}^T \Phi.
\eea
Here $V \in \bR^{(N+6) \times (N+6)}$ is the matrix of (generalized) eigenvectors of the
matrix $A$ (same notation as in Stam's paper) that maps the $N+6$ control
vertices around $f$ to $N+6$ control points in the next subdivision level as shown in
Figure~\ref{fig:InfLayer}(b), so $AV = V \Lambda$ where $\Lambda$ is in a Jordan canonical form.
For the Loop scheme, $\Lambda$ is diagonal when the valence $N$ is greater than 3, but has a
Jordan block of size 2 when $N=3$.
Since the subdivision process is linear and {\it stationary}, i.e. the same linear subdivision rules
are used across different scales, recall Figure~\ref{fig:RedBloodCell}, we have
$$
\mathbf{c}_f^T \mathbf{b}^{N}(v,w) = (A \mathbf{c}_f)^T \mathbf{b}^{N}(2v, 2w),
\quad (v,w) \in \frac{1}{2} \Omega.
$$
Putting these together, we have
\bea \label{eq:ScalingLaw}
\Phi(v, w) = \Lambda^T \Phi(2v,2w),
\quad (v,w) \in \frac{1}{2} \Omega.
\eea
The key point is that these \textit{eigenbasis functions} $\Phi=[\phi_1; \ldots; \phi_{N+6}]$
are easier to evaluate compared to the original basis functions $\mathbf{b}^{N}$: When $N>3$,
$\Lambda={\rm diag}(\lambda_1,\ldots,\lambda_{N+6})$
is diagonal, and
we have $\phi_i(v,w) = \lambda_i \phi(2v,2w)$, apply this recursively we have
\bea \label{eq:ScalingLaw2}
\phi_i(v,w) = \lambda_i^{n-1} \phi(2^{n-1} v,2^{n-1} w), \quad \mbox{when }
(v,w) \in \Omega_k^n.
\eea
(Recall Figure~\ref{fig:OrdinaryPatch}(e).) As a result, each $\phi_i$ is
specified by three -- not infinitely many -- polynomials. Also, in virtue of
\eqref{eq:ScalingLaw2}, $\phi_i$ and its derivatives
can be easily evaluated at arbitrary parameter values after the three polynomials are specified.

These three polynomials can be evaluated based on the same polynomial basis $\mathbf{b}$ from the regular case \eqref{eq:OrdinaryPatch}.
We write $V_i$ as the (generalized) eigen-vector associated to the eigen-basis function $\phi_i$.
For $k=1,2,3$, there are suitable linear maps $M_k$ (expressed as $P_k \overline{A}$ in Stam's paper) so that $M_k V_i$ contains
the data at the $12$ control vertices that determine the polynomial $\phi_i |_{\Omega^1_k}$; see Figure~\ref{fig:InfLayer}(c)-(e).
With an appropriate affine
reparametrization of $\Omega^1_k$ by $\Omega$, denoted as $t_{1,k} : \Omega^1_k \goto \Omega$ by Stam, this polynomial can be expressed as
\bea \label{eq:phi_pieces}
\phi_i|_{\Omega^1_k}(v,w) = (M_k V_i)^T \mathbf{b}(t_{1,k}(v,w)).
\eea
See \cite{Stam:Loop} for details, e.g. on how to exploit the circulant structure in the matrix $A$ in order to facilitate the computation of the related matrices
$V$, $V^{-1}$, $M_k$.

\subsubsection{Computation of $W(\mathcal{V})$, ${A}(\mathcal{V})$, ${V}(\mathcal{V})$, $M(\mathcal{V})$ and their gradients} \label{sec:GradientComp}
\newcommand{\cf}{{\mathbf{c}_f}}
\newcommand{\chf}{{\widehat{\mathbf{c}}_f}}
\newcommand{\cn}{\mathbf{n}}
\newcommand{\Br}{\mathcal{B}_{6}}
\newcommand{\ssf}{{\mathbf{s}_{f}}}
\newcommand{\sfv}{{\mathbf{s}_{f,v}}}
\newcommand{\sfw}{{\mathbf{s}_{f,w}}}
\newcommand{\sfvv}{{\mathbf{s}_{f,vv}}}
\newcommand{\sfww}{{\mathbf{s}_{f,ww}}}
\newcommand{\sfvw}{{\mathbf{s}_{f,vw}}}
\newcommand{\bn}{{\mathbf{n}}}
\newcommand\s[2]{ \mathbf{s}_{#1,#2} }

To summarize the previous section,
a regular patch of a Loop subdivision surface is parameterized by a single degree 4 polynomial on
the reference triangle $\Omega$ for each of the three spatial components,
whereas an irregular patch admits a more complicated parametrization over $\Omega$.
In either case, an
efficient algorithm exists
for evaluating the parametrization and its derivatives at arbitrary parameter values $(v,w) \in \Omega$.
Armed with such evaluation algorithms, we now see how the various functionals and their gradient vectors in \eqref{eq:Helfrich_Subdivision}
can be computed.

\gap
\noindent
{\bf Formulas for $A(\mathcal{V})$ and $\nabla A(\mathcal{V})$.}
We first discuss how to compute the area $A(\mathcal{V})$ of a Loop surface
and the gradient $\nabla A(\mathcal{V})$ of $A$ with respect to $\mathcal{V}$. 
For each regular face $f$ in a control mesh,
we write $A_6(\mathbf{c}_f)$ as the area of the surface patch associated to $f$; and we view $A_6$ as a (real-valued) function
of the variables in the array $\mathbf{c}_f$. Similarly, we write $A_N(\mathbf{c}_f)$ as the area of the surface
patch associated to an irregular face $f$
with an extraordinary vertex of valence $N \neq 6$; in this case, we also write ${\rm val}(f) := N$.
So
\begin{equation}
     \begin{aligned}
       A(\mathcal{V})
                      = \underbrace{\sum_{f \,{\rm regular}} A_6(\mathbf{c}_f)} _{\text{area of regular patches}}+
                        \underbrace{\sum_{N\neq 6} \sum_{f: \,{\rm val}(f)=N} A_N(\mathbf{c}_f)}_{\text{area of irregular patches}}.
     \end{aligned}
     \label{eq:AreaA}
   \end{equation}

Recall from \eqref{eq:OrdinaryPatchC} and \eqref{eq:EigenXform},
the parametrization $\mathbf{s}_f$
can be written as
\bea \label{eq:sf}
\mathbf{s}_f = \mathbf{c}_f^T \mathbf{b}^{N}, \mbox{ where } N={\rm val}(f).
\eea
Note that $\mathbf{c}_f$ is \emph{linearly} related to $\mathcal{V}$, whereas the basis functions
$\mathbf{b}^{N}$ are \emph{independent} of $\mathcal{V}$.  These features allow for an accurate and efficient computation of $A(\mathcal{V})$
and $\nabla A(\mathcal{V})$, as we shall now see.

We define the map $P_f$ by
\bea \label{eq:Pf}
P_f \cV = \cf
\eea
which picks out the local control data around the face $f$ from the global control data $\cV$.

For convenience, drop the subscript and write $\mathbf{s}$ instead of $\mathbf{s}_f$. Also,
we write $\mathbf{s}_1$, $\mathbf{s}_2$, $\mathbf{s}_3$ for the components of $\mathbf{s}$
and $\mathbf{s}_{i,u}$,  $\mathbf{s}_{i,v}$ for their partial derivatives.

Note that
\begin{equation} \label{eq:area_singlepatch}
\begin{aligned}
A_N(\mathbf{c}_f) &= \iint_\Omega \| \mathbf{n}(v,w) \| \, dv\, dw, \mbox{ where}\\
\mathbf{n} &= \frac{\partial \mathbf{s}}{\partial v} \times \frac{\partial \mathbf{s}}{\partial w}
           = \left[\s{2}{v} \s{3}{w} - \s{2}{w} \s{3}{v}, \;\;
                   \s{3}{v} \s{1}{w} - \s{3}{w} \s{1}{v} , \;\;
                   \s{1}{v} \s{2}{w} - \s{1}{w} \s{2}{v} \right].
\end{aligned}
\end{equation}
Again, we drop the subscript $f$ and write $\mathbf{c}_{\cdot1}$, $\mathbf{c}_{\cdot2}$, $\mathbf{c}_{\cdot3}$ to refer to the
columns of $\mathbf{c}_f$.

When $f$ is a regular face, $\mathbf{s}_i = \mathbf{c}_{\cdot i}^T \mathbf{b}$, so
\begin{equation} \label{eq:siv}
\begin{aligned}
\s{i}{v} = &\mathbf{c}_{\cdot i}^T \mathbf{b}_v, \;\; \s{i}{w} = \mathbf{c}_{\cdot i}^T \mathbf{b}_w,
 \;\; i=1,2,3.
\end{aligned}
\end{equation}

The gradient of $\mathbf{n}_1$ with respect to $\mathbf{c}_f$, organized as a $12 \times 3$ array
(i.e. same dimension as $\mathbf{c}_f$), can be expressed as:\footnote{Here and below, we have to deal with a number of
scalar quantities $S$ that vary with \textbf{both} the local control data $\cf$ and parameter values $(v,w)$
(e.g.
$\|\mathbf{n}\|$, $\mathbf{n}_i$, $E$, $F$, $G$, $e$,
$f$, $g$, etc..) In order to avoid confusion, we use the notation `$\nabla_\cf S$' to
denote the gradient vector  of $S$ viewed as a function of $\cf$; the gradient `vector' is
structured as
an array of the same size as $\cf$. Likewise, the gradient `vectors'
$\nabla W(\mathcal{V})$, $\nabla {A}(\mathcal{V})$, $\nabla {V}(\mathcal{V})$, $\nabla M(\mathcal{V})$
are structured as $\#V \times 3$ arrays, i.e. the same dimensions as $\mathcal{V}$.}
\bea
\begin{aligned}
  \nabla_{\cf} \mathbf{n}_1 &= \nabla_{\cf} (\mathbf{s}_{2,v} \mathbf{s}_{3,w} - \mathbf{s}_{2,w} \mathbf{s}_{3,v})\\
  &=  \big[ \mathbf{0}, \;\; \mathbf{b}_v,\;\; \mathbf{0} \big] \mathbf{s}_{3,w} +
\mathbf{s}_{2,v}\big[ \mathbf{0},\;\; \mathbf{0},\;\; \mathbf{b}_w \big]  -
 \big[ \mathbf{0}, \;\; \mathbf{b}_w,\;\; \mathbf{0} \big] \mathbf{s}_{3,v} -
\mathbf{s}_{2,w}\big[ \mathbf{0},\;\; \mathbf{0},\;\; \mathbf{b}_v \big] \\
  & = \big[ \mathbf{0},\;\; \mathbf{s}_{3,w} \mathbf{b}_v - \mathbf{s}_{3,v} \mathbf{b}_w, \;\;
                            \mathbf{s}_{2,v} \mathbf{b}_w - \mathbf{s}_{2,w} \mathbf{b}_v \big].
\end{aligned}
\eea
Similarly,
\begin{equation}
\begin{aligned}
  \nabla_{\cf} \mathbf{n}_2 &= \big[-\mathbf{s}_{3,w} \mathbf{b}_v + \mathbf{s}_{3,v} \mathbf{b}_w,\;\; \mathbf{0},\;\;
                                    -\mathbf{s}_{1,v} \mathbf{b}_w + \mathbf{s}_{1,w} \mathbf{b}_v],\\
  \nabla_{\cf} \mathbf{n}_3 &= \big[\mathbf{s}_{2,w} \mathbf{b}_v - \mathbf{s}_{2,v} \mathbf{b}_w,\;\;
                                     \mathbf{s}_{1,v} \mathbf{b}_w - \mathbf{s}_{1,w} \mathbf{b}_v\;\; \mathbf{0}\big].
\end{aligned}
\end{equation}
Next, we have
\begin{equation} \label{eq:n_and_gradn}
\begin{aligned}
\| \mathbf{n} \| &= \sqrt{ \langle \mathbf{n}, \mathbf{n} \rangle }, \\
\nabla_{\cf} \| \mathbf{n} \| &= \nabla_{\cf} \langle \mathbf{n}, \mathbf{n} \rangle^{1/2}
= \frac{1}{2 \langle \mathbf{n}, \mathbf{n} \rangle^{1/2}} \nabla_{\cf} \langle \mathbf{n}, \mathbf{n} \rangle \\
&= \frac{1}{ \langle \mathbf{n}, \mathbf{n} \rangle^{1/2}}
\left(
\mathbf{n}_1 \nabla_{\cf} \mathbf{n}_1 + \mathbf{n}_2 \nabla_{\cf} \mathbf{n}_2 + \mathbf{n}_3 \nabla_{\cf} \mathbf{n}_3
\right).
\end{aligned}
\end{equation}
Therefore, the local area functional $A_N$ and its gradient
\begin{equation} \label{eq:grad_A6}
\begin{aligned}
\nabla A_N(\mathbf{c}_f) = \iint_\Omega \nabla_{\cf} \| \mathbf{n}(v,w) \|\, dv\, dw,
\end{aligned}
\end{equation}
can be computed based on the control data $\mathbf{c}_f = [\mathbf{c}_{\cdot 1}, \mathbf{c}_{\cdot 2}, \mathbf{c}_{\cdot 3}]$
and the basis function $\mathbf{b}$ via \eqref{eq:siv}-\eqref{eq:n_and_gradn}. Together with
\eqref{eq:Pf} and the chain rule, we can compute the total area and its gradient with respect to $\mathcal{V}$ by
\begin{equation} \label{eq:AreaA_Final}
     \begin{aligned}
       A(\mathcal{V}) = \sum_{N} \sum_{f: \, {\rm val}(f)=N} A_N(P_f \mathcal{V}), \;\;\;
     \nabla A(\mathcal{V})= \sum_{N} \sum_{f: \, {\rm val}(f)=N} P_f^T  \nabla A_N (\cf).
     \end{aligned}
\end{equation}

\noindent
{\bf Formulas for $W(\mathcal{V})$, $V(\mathcal{V})$, $M(\mathcal{V})$, and their gradients.}
Similar to $A(\mathcal{V})$, we aim to express the other three functionals in the C-H model in terms of $\cf$ and the basis functions
$\mathbf{b}$ and $\Phi$. For $W$ and $M$, we need an expression for the mean curvature.
Recall that
$$E = \langle \sfv, \sfv \rangle, \; F = \langle \sfv, \sfw \rangle, \; G = \langle \sfw, \sfw \rangle$$
represent the first fundamental form of the surface $\mathbf{s}_f$,
whereas
$$e = \langle \sfvv, \bn \rangle/\|\bn\|, \; f = \langle \sfvw, \bn \rangle/ \|\bn\|, \;
g = \langle \sfww, \bn \rangle/ \|\bn\| $$ represent the second fundamental form.
The mean curvature can be expressed as
\begin{align} \label{eq:H}
   H &= \frac{eG-2fF+gE}{2(EG-F^2)}.
     \end{align}
Therefore
\bea \label{eq:M}
M = -\iint H \, dA = -\sum_{f \in \mathcal{F}} \iint_\Omega \frac{eG-2fF+gE}{2(EG-F^2)} \| \bn \| \, dv \, dw
= -\sum_{f \in \mathcal{F}} \iint_\Omega \frac{\bar{e}G-2\bar{f}F+\bar{g}E}{2(EG-F^2)} \, dv \, dw,
\eea
where
$\bar{e} := \langle \sfvv, \bn \rangle$, 
$\bar{f} := \langle \sfvv, \bn \rangle$, 
$\bar{g} := \langle \sfvv, \bn \rangle$. 
Similarly,
\bea \label{eq:W}
W = \iint H^2 \, dA
= \sum_{f \in \mathcal{F}} \iint_\Omega \left[\frac{\bar{e}G-2\bar{f}F+\bar{g}E}{2(EG-F^2)}\right]^2
\frac{1}{\| \bn \|} \, dv \, dw.
\eea
\begin{remark} When $f$ is an irregular face, the mean curvature can potentially blows up when approaching the extraordinary vertex, however
it is proved in \cite{ReifSchroder:Curvature} that the corresponding integrals (called ${M}_N$ and ${W}_N$ below) are
always finite.
\end{remark}

By the divergence theorem, with the choice of the vector field $\vec{X}(x,y,z) = x \mathbf{i} + y \mathbf{j} + z \mathbf{k}$,
the volume enclosed by a surface can be expressed as a surface integral:
\bea \label{eq:V}
V =  \frac{1}{3} \iiint {\rm div} \vec{X} \, dx\, dy\, dz= \frac{1}{3} \iint \vec{X} \cdot \bn/\|\bn\|\, dA
= \frac{1}{3}  \sum_{f \in \mathcal{F}}\iint_\Omega \langle \mathbf{s}_f, \bn \rangle \, dv\, dw.
\eea
By \eqref{eq:sf},
we can express the rightmost integral in \eqref{eq:M}-\eqref{eq:V} in terms of $\mathbf{c}_f$ and
$\mathbf{b}^{N}$ when where $f$ ranges over all faces and $N$ ranges over all valences existing in the control mesh;
we denote the integral
by $M_N(\cf)$, $W_N(\cf)$ and $V_N(\cf)$, respectively.

We explain how to compute the gradients of $M_N(\cf)$, $W_N(\cf)$ and $V_N(\cf)$.

 The gradients of $E$, $F$, $G$, $\bar{e}$, $\bar{f}$, $\bar{g}$ can be computed as follows
\bea
\begin{aligned} \label{eq:FirstFundamental}
\mathbf{s}_i = \mathbf{c}_{\cdot i}^T \mathbf{b}, & \; \mathbf{s}_{i,v}= \mathbf{c}_{\cdot i}^T \mathbf{b}_v, \;
\mathbf{s}_{i,v}= \mathbf{c}_{\cdot i}^T \mathbf{b}_w, \;
\mathbf{s}_{i,vv}= \mathbf{c}_{\cdot i}^T \mathbf{b}_{vv}, \;
\mathbf{s}_{i,vw}= \mathbf{c}_{\cdot i}^T \mathbf{b}_{vw}, \;
\mathbf{s}_{i,ww}= \mathbf{c}_{\cdot i}^T \mathbf{b}_{ww}, \\
\nabla_\cf E &= \sum_{i=1,2,3} \nabla_\cf (\mathbf{c}_{\cdot i}^T \mathbf{b}_v)^2 =
2 \big[ \s{1}{v} \mathbf{b}_v, \;
        \s{2}{v} \mathbf{b}_v, \;
        \s{3}{v} \mathbf{b}_v \big],\\
\nabla_\cf F &= \sum_{i=1,2,3} \nabla_\cf (\mathbf{c}_{\cdot i}^T \mathbf{b}_v)(\mathbf{c}_{\cdot i}^T \mathbf{b}_w)
=\big[\s{1}{v} \mathbf{b}_w + \s{1}{w} \mathbf{b}_v, \;
\s{2}{v} \mathbf{b}_w + \s{2}{w} \mathbf{b}_v, \;
\s{3}{v} \mathbf{b}_w + \s{3}{w} \mathbf{b}_v \big],\\
\nabla_\cf G &= \sum_{i=1,2,3} \nabla_\cf (\mathbf{c}_{\cdot i}^T \mathbf{b}_w)^2 =
2 \big[\s{1}{w} \mathbf{b}_w, \; \s{2}{w} \mathbf{b}_w, \; \s{3}{w} \mathbf{b}_w \big],\\
\nabla_\cf \bar{e} &= \sum_i \nabla_\cf (\mathbf{c}_{\cdot i}^T \mathbf{b}_{vv}) \, \bn_i
= \big[\s{1}{vv}\nabla_\cf \bn_1, \;
       \s{2}{vv} \nabla_\cf \bn_2, \;
       \s{3}{vv} \nabla_\cf \bn_3 \big]
+ \big[ \bn_1 \mathbf{b}_{vv}, \; \bn_2 \mathbf{b}_{vv}, \; \bn_3 \mathbf{b}_{vv} \big], \\
\nabla_\cf \bar{f}
&= \big[\s{1}{vw}\nabla_\cf \bn_1, \;
       \s{2}{vw} \nabla_\cf \bn_2, \;
       \s{3}{vw} \nabla_\cf \bn_3 \big]
+ \big[ \bn_1 \mathbf{b}_{vw}, \; \bn_2 \mathbf{b}_{vw}, \; \bn_3 \mathbf{b}_{vw} \big], \\
\nabla_\cf \bar{g}
&= \big[\s{1}{ww}\nabla_\cf \bn_1, \;
       \s{2}{ww} \nabla_\cf \bn_2, \;
       \s{3}{ww} \nabla_\cf \bn_3 \big]
+ \big[ \bn_1 \mathbf{b}_{ww}, \; \bn_2 \mathbf{b}_{ww}, \; \bn_3 \mathbf{b}_{ww} \big]. \\
\end{aligned}
\eea
By \eqref{eq:M}, 
\bea
\begin{aligned} \label{eq:grad_M6}
\nabla M_N = \nabla_\cf \iint_\Omega \underbrace{\frac{\bar{e}G-2\bar{f}F+\bar{g}E}{2(EG-F^2)}}_{=: \bar{H}}
\, dv \, dw
= \iint_\Omega \nabla_\cf \bar{H} \, dv \, dw.
\end{aligned}
\eea
Thanks to \eqref{eq:FirstFundamental}, the integrand $\nabla_\cf \bar{H}$ can be computed by
the product and quotient rules.
After $\nabla M_N$ is computed, $\nabla W_N$ can then be computed by
\bea \nonumber
\nabla W_N &= \nabla_\cf \iint_\Omega \left[\frac{\bar{e}G-2\bar{f}F+\bar{g}E}{2(EG-F^2)}\right]^2
\frac{1}{\| \bn \|} \, dv \, dw \\
&=
\iint_\Omega \frac{2 \bar{H}}{\| \bn \|} \nabla_\cf \bar{H} -\frac{\bar{H}^2}{\| \bn \|^2} \nabla_\cf \| \bn \|
\, dv \, dw. \label{eq:W_6_2}
\eea
Note that every term in the integrand of \eqref{eq:W_6_2} was computed previously.

Finally, for $V_N$ we have
\bea \label{eq:V_6}
\begin{aligned}
\nabla V_N &= \nabla_\cf \iint_\Omega \langle \mathbf{s}_f, \bn \rangle \, dv \, dw =
\iint_\Omega \sum_i \nabla_\cf(\mathbf{c}_{\cdot i}^T \mathbf{b}) \bn_i\, dv \, dw\\
&= \iint_\Omega
\left\{ \big[\mathbf{n}_1 \mathbf{b}, \; \mathbf{n}_2 \mathbf{b}, \; \bn_3 \mathbf{b}\big] +
\mathbf{s}_1 \nabla_\cf \bn_1
+ \mathbf{s}_2 \nabla_\cf \bn_2 + \mathbf{s}_3 \nabla_\cf \bn_3 \right\}
\, dv \, dw.
\end{aligned}
\eea

With all the local functionals and their gradients with respect to the local control data
computed, the global functionals and their gradients with respect to the global control data
can be computed exactly as in \eqref{eq:AreaA_Final}:
\begin{equation} \label{eq:MWV_Final}
\begin{aligned}
M(\mathcal{V}) = \sum_{N} \sum_{f: \, {\rm val}(f)=N} M_N(P_f \mathcal{V}), \;\;\;
\nabla M(\mathcal{V}) =  \sum_{N} \sum_{f: \, {\rm val}(f)=N} P_f^T  \nabla M_N (\cf),
\end{aligned}
\end{equation}
and similarly for $W(\mathcal{V})$, $\nabla W(\mathcal{V})$, $V(\mathcal{V})$, $\nabla V(\mathcal{V})$.

\subsection{Implementation details}
\label{sec:details}
In the actual numerical computation of $A(\mathcal{V})$
and $\nabla A(\mathcal{V})$ based on \eqref{eq:AreaA_Final}, we use a symmetric 7-point Gauss quadrature rule on a triangle,
with accuracy order 5 (see for example \cite{Cowper:Quadrature}),
to approximate the integrals in $A_N(\cf)$, $V_N(\cf)$, $M_N(\cf)$, $W_N(\cf)$, and their gradients
(recall \eqref{eq:grad_A6}, \eqref{eq:grad_M6}-\eqref{eq:V_6}).
The approximation is done using a composite quadrature on $\Omega$ using a uniform grid of size $1/n$. 
Typically we use $n=8$ or $16$ at coarse subdivision levels,
and $n=1$ or $2$ at fine levels. We choose $n$ to be a dyadic integer in order to take advantage of the subdivision structure of
$\mathbf{s}_{f}$ when $N \neq 6$.

\begin{figure}[ht]
\begin{tabular}{cc}
\begin{minipage}[l]{.49\linewidth}
\includegraphics[height=2.8in]{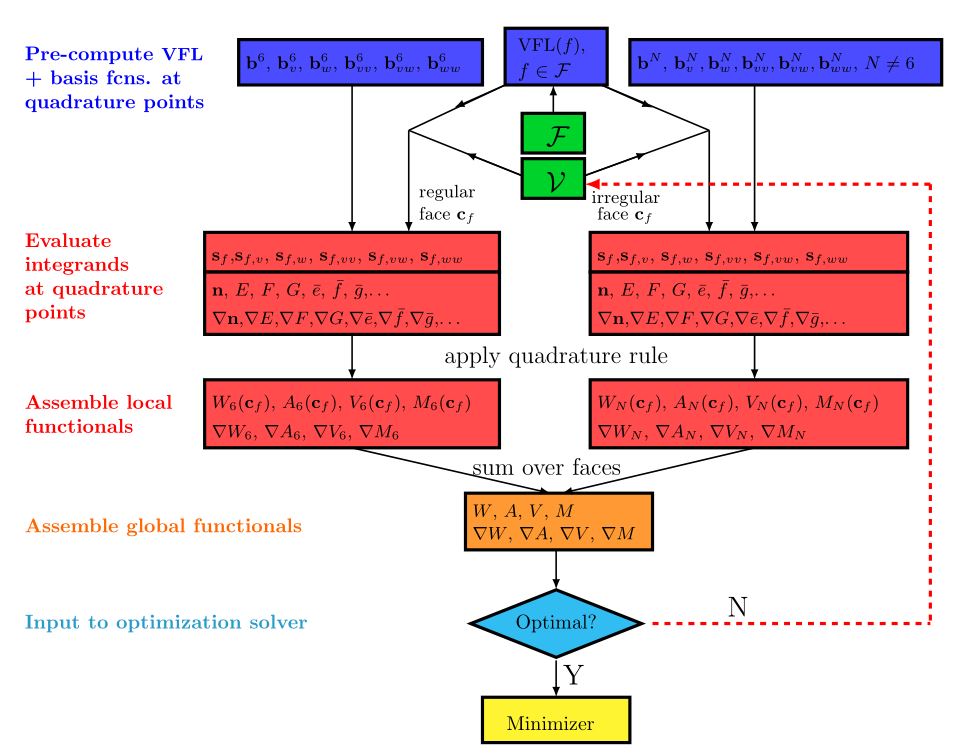}
\end{minipage}
&  \hspace{.2in}
\begin{minipage}[r]{.4\linewidth}
\includegraphics[height=1.2in]{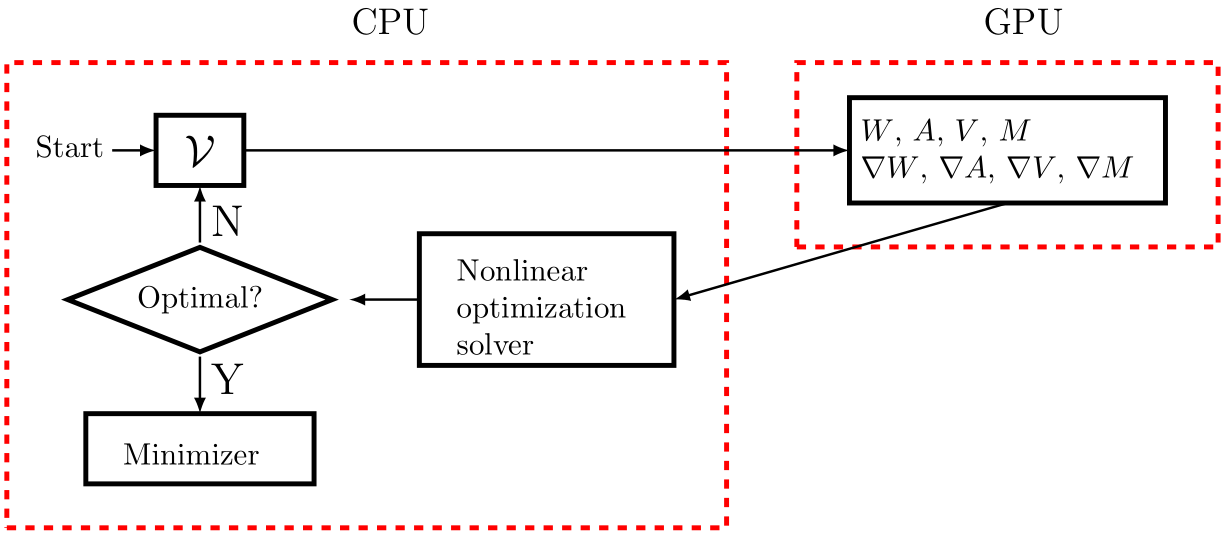}
\end{minipage}\\
\end{tabular}
\caption{(left) computation structure of the SS solver (right) architecture of the GPU enhanced solver}
 \label{fig:flowchart}
\end{figure}

A key implementation detail is that quantities independent of $\mathcal{V}$ are precomputed before entering the optimization
loop. These quantities include $\mathbf{b}^6(u,v)$ and $\mathbf{b}^N(u,v) = V_N^{-T} \Phi_N(u,v)$
(for only those extraordinary valences $N$ that show up in $\mathcal{{F}}$)
 and their derivatives
 at the quadrature points, evaluated
using Stam's algorithm (Sections~\ref{sec:SSIntro}.)
A separate pre-processing step computes the maps $P_f$ for each face $f$ in $\mathcal{F}$, as $P_f$ depends only
on the connectivity information in $\mathcal{F}$. Of course, one should not
store $P_f$ as a $(N+6) \times \#V$ (dense) matrix as suggested by \eqref{eq:Pf}.
Instead, it suffices to store the
information in $P_f$
by a list of $N+6$ integer indices in $\{1, \ldots, \#V\}$
which keeps track of the indices of $\mathbf{c}_f$ in the global vertex list $\mathcal{V}$; we denote
this list of vertex indices by
$
\mathrm{VFL}(f).
$
These preprocessing steps speed up our solver significantly already in a sequential implementation. See Figure~\ref{fig:flowchart}(left) for
the basic organization of our SS solver.
\gap

We note that both the PL and SS functionals are quite easily parallelized; we provide in {\sc Wmincon}
parallel CUDA implementations of all the functionals considered in this article.
Together with the aforementioned precomputation trick, we are able to solve a lot of instances of Willmore, Canham and Helfrich
problems not addressed in \cite{Feng2006394}.
Figure~\ref{fig:flowchart}(right) illustrates the basic architecture of our GPU enchanced solver for both the PL and SS methods.

For our implementation of the solver, we primarily use SNOPT \cite{SNOPT:SIAM} and the `SQP' option of \verb$fmincon$
in the Matlab optimization toolbox.
These solvers are
designed for smooth objectives and constraints,
which is the case for the SS method and most of the PL methods.
An interesting exception is the case of $W_{\rm Bobenko}$, which is non-smooth
whenever the PL surface has two adjacent triangles of which the four vertices lie on a common circle.
In this case, our use of the non-smooth
optimization algorithm GRANSO \cite{CMO:GRANSO} was instrumental to the
discovery of Proposition~\ref{prop:ST} below.

All three solvers use a quasi-Newtwon (BFGS) algorithm, coupled with the
sequential quadratic programming (SQP) method
\cite[Chapter 18]{NoceWrig06} for handling constraints. Being all BFGS-SQP-based algorithms, they have decisively different properties even for smooth problems, see Section~\ref{sec:Further}.

\section{Analysis} \label{sec:analysis}
In this section, we give a theoretical justification for why the SS methods succeed in solving the Willmore problem (Section~\ref{sec:positive}), but
that the PL methods presented in Section~\ref{sec:Functionals_PLSS} are bound to fail (Section~\ref{sec:negative}).
For computational experiments illustrating
the results in this section, see the companion conference paper \cite{BYY:Enumath17}.

\subsection{Convergence of SS method} \label{sec:positive}
In this section, we establish Theorem~\ref{prop:MainPositive}.
The argument is based on the existence result for Willmore surfaces  established using
 the direct method in the calculus of variations (Theorem~\ref{thm:Simon:Willmore})
\cite{Simon:Willmore,Kusner:Willmore,BK:Willmore,kusner1989}, and an
observation (Proposition~\ref{prop:MinSequenceSplines}) connecting the existence
result to our conforming subdivision methods; the connection requires a density result (Theorem~\ref{thm:Arden})
from the theory of subdivision surfaces.

A sequence $(x^k)_{k=0}^\infty$ in a space $\mathcal{X}$ is called a \emph{minimizing sequence} for a functional $\mathcal{E}: \mathcal{X} \goto \bR$
if
$$\lim_{k\goto \infty} \mathcal{E}(x^k) = \inf_{x\in \mathcal{X}} \mathcal{E}(x).$$

\begin{proposition} \label{prop:MinSequenceSplines} \normalfont
Let $\mathcal{E}: \mathcal{X} \goto \bR$ be a continuous functional on a topological space $\mathcal{X}$.
Assume that we have a nested sequence of subspaces $\{\mathcal{S}^j : j=0,1,2,\ldots \}$ of $\mathcal{X}$ such that
$\bigcup_j \mathcal{S}^j$ is dense in $\mathcal{X}$.
Assume that a minimizer
$x^\ast \in \argmin_{\mathbf{x} \in \mathcal{X}} \mathcal{E}(x)$ exists.
Then any sequence of `approximate $\mathcal{E}$-minimizers' $x^j \in \mathcal{S}^j$, i.e.
$$
\mathcal{E}(x^j) = \inf_{x \in \mathcal{S}^j} \mathcal{E}(x) + o(1), \quad j \goto \infty,
$$
is a minimizing sequence for $\mathcal{E}$.
\end{proposition}
\pf By the denseness assumption, there exists a sequence $\widetilde{x}^j \in \mathcal{S}^j$
converging to $x^\ast$ in $\mathcal{X}$. By continuity, $\mathcal{E}(\widetilde{x}^j) \goto \mathcal{E}(x^\ast)$.
But we also have
$$\inf_{x \in \mathcal{X}} \mathcal{E}(x) \leq \mathcal{E}(x^j) =  \inf_{x \in \mathcal{S}^j} \mathcal{E}(x) + o(1) \leq \mathcal{E}(\widetilde{x}^j) + o(1),$$
so $\mathcal{E}({x}^j) \goto \mathcal{E}(x^\ast)=\inf_{x\in \mathcal{X}}\mathcal{E}(x)$.
Thus, $(x^j)_{j=0}^\infty$ is a minimizing sequence for $\mathcal{E}$.
\eop

Note how this result relies on the conforming and dense nature of the spaces $\mathcal{S}^j$. Notice
also that the use of approximate minimizers
frees us from the assumption that a minimizer of $\mathcal{E}$ over each $\mathcal{S}^j$ exists.

\newcommand{\Imm}{{\rm Imm}}
\newcommand{\ImmX}{{\rm Imm}_X}

We now state the key existence result established in \cite{Simon:Willmore,BK:Willmore,Kusner:Willmore}.
Since we are analyzing a parameteric method, it is more convenient to state the result in terms of parametrizations.
For this purpose, we assume
$\Sigma$
is any (reference) genus $g$ surface with a smooth enough  differentiable structure.
A $C^{1,1}$ differentiable structure would suffice for our purpose, as it is enough to support the definition of
not only $C^{1}(\Sigma, \bR^3)$ and $C^{1,\alpha}(\Sigma, \bR^3)$, $0<\alpha\leq 1$, but also the whole range of Sobolev spaces
$W^{2,p}(\Sigma, \bR^3)$, $p \in [1,\infty]$. (See \cite{arden:2001} and \cite[Section 4.2.3]{EvansGariepy:Measure}.)
We shall work with the Banach space
$$X:= X(\Sigma):= C^{1}(\Sigma, \bR^3) \cap W^{2,2}(\Sigma, \bR^3)$$
(normed by $\|\mathbf{x}\|_X = \|\mathbf{x}\|_{C^1}+\|\mathbf{x}\|_{W^{2,2}}$)
and the nonlinear, open subspace
\begin{align} \label{eq:SpaceX}
\begin{split}
\ImmX := \ImmX(\Sigma):=\big\{ f \in C^{1}(\Sigma, \bR^3) \cap W^{2,2}(\Sigma, \bR^3) \; | \; {\rm rank} (df_x) =2, \; \forall x \in \Sigma \big\}, \quad
\end{split}
\end{align}
on which the Willmore energy
$$W: {\rm Imm}_X \goto [0,\infty)$$
is well-defined and continuous.

The following is a reformulation of the well-known existence result pioneered by L. Simon \cite{Simon:Willmore} and completed in
 \cite{BK:Willmore,Kusner:Willmore}:
\begin{theorem}[{\bf Existence of $W$-minimizer of genus $g$, via the direct method}] \label{thm:Simon:Willmore} \normalfont
Let $\Sigma$ be a closed orientable surface of genus $g$ with a $C^{1,1}$ differential structure.
For any minimizing sequence $\mathbf{x}_j \in \ImmX$ for $W$, there exists
a subsequence
$\mathbf{x}_{j'}$, and a sequence of M\"obius transformations $G_{j'}$ in $\bR^3$, such that
the sequence of immersed surfaces $G_{j'} (\mathbf{x}_{j'}(\Sigma))$ converges in Hausdorff distance
to an immersed surface $\mathbf{x}_\ast(\Sigma)$,
$\mathbf{x}_\ast \in \ImmX$.
As such, $W(\mathbf{x}_\ast) = \inf_{ \mathbf{x} \in \ImmX} W(\mathbf{x})$.
\end{theorem}

\begin{remark}
The Willmore minimizers are known to be embedded surfaces. We choose to work with general immersed surfaces because in our numerical method
we do not have any mechanism built
in to avoid self-intersections; another reason is that the solutions of the Helfrich problem for some values of $v_0$ and $m_0$
do have self-intersections.
\end{remark}

\begin{remark}
While \cite{Simon:Willmore,BK:Willmore,Kusner:Willmore} prove the existence of
a $W$-minimizer over all \emph{infinitely smooth} genus $g$ immersed surfaces, our minimization space
is taken to be the bigger \eqref{eq:SpaceX}.
This causes no problem as
$C^\infty(\Sigma, \bR^3) \cap \ImmX$ is dense in $\ImmX$
and $W: \ImmX \goto \bR$ is continuous. (If $E: A\goto \bR$ is a continuous functional on a topological space $A$,
and $B$ is a dense subspace of $A$, then any minimizer of $E: B \goto \bR$ must also be a minimizer of $E: A \goto \bR$.)
Here, the definition of $C^\infty(\Sigma, \bR^3)$, with $\Sigma$ endowed with a $C^{1,1}$ differentiable structure,
requires the fact that a maximal $C^{1}$-atlas contains a $C^\infty$-atlas. 
In fact there is a real analytic sub-atlas,\footnote{The result that every $C^1$ manifold admits compatible $C^\infty$
and analytic ($C^\omega$) structures is due to Whitney \cite{Whitney:Manifolds}.} with respect to which
the minimizer $\mathbf{x}_\ast$ in Theorem~\ref{thm:Simon:Willmore} is real-analytic \cite{Simon:Willmore}.
\end{remark}

The Loop subdivision scheme defines:
\begin{itemize}
\item a $C^{2}$-compatible
atlas, via its characteristic maps, on any base complex $K$, and
\item a linear
space of scalar-valued subdivision functions $\euS(K^j)$ at each subdivision level.
\end{itemize}
The atlas given by
the characteristic maps of Loop's scheme
turns $K$ into a $C^{2}$-manifold. For details, see \cite{arden:2001}.
 The subdivision functions are well-known to be in $C^1(K, \bR)$
\cite{Reif:Degree,Zorin:Extraordinary,Zorin:ExtraordinaryC1}. By a trivial extension of
\cite[Theorem 42]{arden:2001}, we also have
$$\euS(K^j) =: \euS(K^j, \bR) \subset W^{2,p}(K, \bR), \quad \forall p \in [1,p^\ast),$$
where $p^\ast>2$ is given by \eqref{eq:p0} below.

When we extend these subdivision functions componentwise to map into $\bR^m$, we denote the
corresponding linear space by $\euS(K^j, \bR^m)$.

Arden's thesis \cite{arden:2001} essentially establishes the following result. While he focuses on the $p=2$ case, the proof can be extended to
a range of values of $p$ that is strictly bigger than the interval $[1,2]$.
\begin{theorem} \normalfont \label{thm:Arden}
There exists a $p^\ast>2$ such that
for any $p \in [1,p^\ast)$,
the space of Loop subdivision functions at all levels $\bigcup_{j=0}^\infty \euS(K^j, \bR)$ is dense in $W^{2,p}(K, \bR)$.
\end{theorem}
\begin{remark}
The value of $p^\ast$ above depends on the set of valences
present in the underlying simplicial complex, denoted by ${\rm val}(\cF)$,
and the values of the sub- and sub-sub-dominant eigenvalues, $\lambda_n$
and $\mu_n$, of
the subdivision matrix $A$ corresponding to the different valences $n \in {\rm val}(\cF)$.
(For the definition of $A$, recall the comments after \eqref{eq:EigenXform} in Section~\ref{sec:SSIntro}.)
We have
\begin{align} \label{eq:p0}
\begin{split}
p^\ast &= \min_{n \in {\rm val}(\cF)} \frac{2 \log(\lambda_n)}{\log(\lambda_n^2/\mu_n)}, \;\;\;
\lambda_n = \frac{3}{8}+\frac{1}{4}\cos(\frac{2\pi}{n})\;\;\;
\mbox{and } \;\;\;
\mu_n= \left\{
        \begin{array}{ll}
         1/8 & \hbox{if $n=3$,} \\
          \lambda_n^2  & \hbox{if $n=4,5,6$,}\\
          \frac{3}{8}+\frac{1}{4}\cos(\frac{4\pi}{n}) & \hbox{if $n\geq 7$.}
        \end{array}
      \right.
\end{split}
\end{align}

It is easy to check that $p^\ast$ approaches 2 from above as $\max ({\rm val}(\cF)) \goto \infty$.
For details, see \cite{XieYu:ApproximationSS}.
\end{remark}
In particular, Theorem~\ref{thm:Arden} implies, by Morrey's inequality, 
that $\bigcup_{j=0}^\infty \euS(K^j, \bR^3)$ is dense in $X=C^1(K, \bR^3) \cap W^{2,2}(K, \bR^3) $.

For notational simplicity, we write
$$
\euS^j = \euS(K^j, \bR^3), \quad \Imm_{\euS^j} := \Imm_{\euS^j}(K) := \euS(K^j, \bR^3) \cap \ImmX(K).
$$

\begin{corollary}\normalfont \label{Cor:Density}
$\bigcup_j \Imm_{\euS^j}$ is dense in $\ImmX(K)$.\footnote{In fact, $\Imm_{\euS^j}$ is dense in $\euS^j$ for each $j$, so,
together with the consequence of Arden's result, we conclude that
$\bigcup_j \Imm_{\euS^j}$ is dense in $X(K)$. This stronger result requires extra arguments based on known results in the theory of
subdivision surface and an application of Sard's theorem.}
\end{corollary}
\pf Since $\bigcup_j \euS^j$ is dense in $X$, and $\ImmX$ is open in $X$,
$$\ImmX \cap \bigcup_j \euS^j = \bigcup_j \Imm_{\euS^j}$$ is dense in $\ImmX$.
\eop

Combining Corollary~\ref{Cor:Density} and Proposition~\ref{prop:MinSequenceSplines} above, we have the conclusion that
if $\mathbf{x}_j$ is an approximate $W$-minimizer over $\Imm_{\euS^j}$  in
the sense of Proposition~\ref{prop:MinSequenceSplines}, then the sequence $(\mathbf{x}_j)_j$ is a minimizing sequence in the setting of
Theorem~\ref{thm:Simon:Willmore}. So by Theorem~\ref{thm:Simon:Willmore}, we have:
\begin{theorem} \normalfont \label{prop:MainPositive}
Let $\Imm_{\euS^j}$ be the space of immersed
Loop subdivision surfaces on the $j$-times subdivided
complex $K^j$ of a genus $g$ complex $K=K^0$ (as above).
If $\mathbf{x}_j$ is an approximate $W$-minimizer over $\Imm_{\euS^j}$  in
the sense of Proposition~\ref{prop:MinSequenceSplines},
then there is a subsequence $\mathbf{x}_{j'}$ so that the surfaces $\mathbf{x}_{j'}(K)$,
with suitable M\"obius transformations applied,
converge in Hausdorff distance
to a surface $\mathbf{x}_\ast(K)$ with $\mathbf{x}_\ast \in {\rm Imm}_X$.
This $\mathbf{x}_\ast(K)$ is a genus $g$ Willmore minimizer.
\end{theorem}

This result should explain the empirical success of the SS method;
see Section~\ref{sec:Future} for a possible improvement of this result.
In contrast, we next illustrate why the PL $W$-minimizers typically have nothing to do with a continuous Willmore minimizer.

\subsection{Failure of Naive PL methods} \label{sec:negative}
We begin by establishing a negative result for Bobenko's Willmore energy, which
may sound like an unworthy result given its restrictive consistency property. 
But here we work in the setting most favorable for $W_{\rm Bobenko}$ in terms of consistency, namely the
setting for a regularly triangulated torus. We know from \cite{Bobenko05discretewillmore} that
by successively finer regular sampling and triangulation of a torodial Dupin cylcide
at its curvature lines one obtains PL tori
with $W_{\rm Bobenko}$-energy converging to the Willmore energy of the Dupin cylcide.
This applies in particular to the Clifford tori (i.e. Duplin cyclides gotten from all possible
M\"obius transformations of the torus of revolution with radii ration $1: \sqrt{2}$),
which are the only minimizers of the genus
1 Willmore problem \cite{MarquesNeves:Willmore}.
However,
we show below that we will never get an approximation of a Clifford torus by {\it minimizing} $W_{\rm Bobenko}$.

The construction used in this proof happens to be applicable, after a twist, to proving a similar negative
result for other PL $W$-energies.

\subsubsection{$W_{\rm Bobenko}$}
\begin{figure}[h]
\centerline{
\begin{tabular}{cccccc}
\includegraphics[height=1.0in]{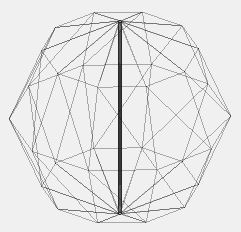}& \hspace{-.1in}
\includegraphics[height=1.0in]{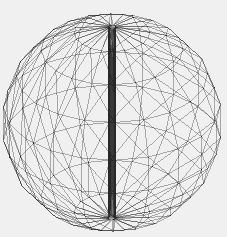}& \hspace{-.1in}
\includegraphics[height=1.0in]{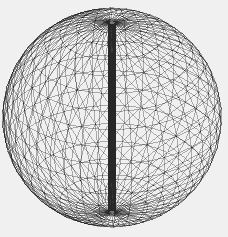}& \hspace{-.1in}
\includegraphics[height=1.0in]{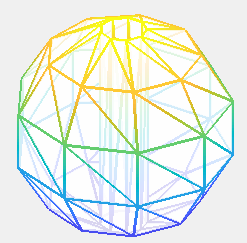}& \hspace{-.1in}
\includegraphics[height=1.0in]{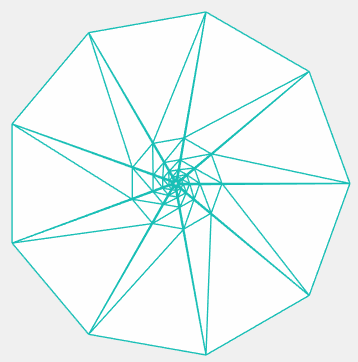}& \hspace{-.1in}
\includegraphics[height=1.0in]{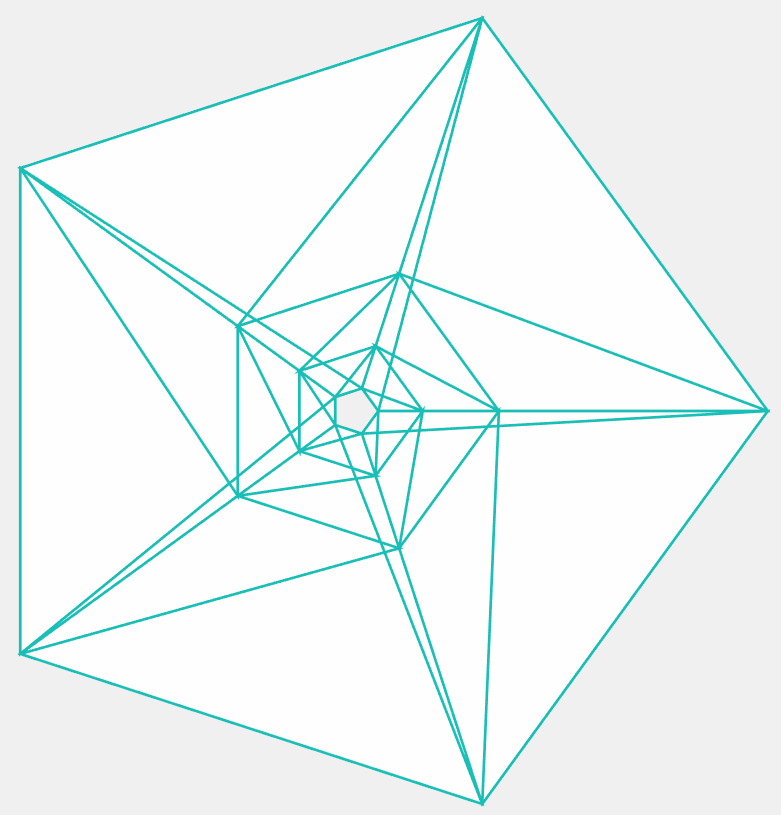} \\
(a) $T_{5,10,.01}$ & (b) $T_{10,20,.05}$ & (c) $T_{20,40,.1}$ & (d) $T_{6,9,.2}$ & (e) $T_{6,9,.2}'$ & (f) $T_{4,5,.3}'$
 \end{tabular}} 
 \caption{Various spherical tori ((a)-(d)) and flattened tori based on sphere inversions ((e)-(f))}
 \label{fig:tori}
\end{figure}

For every grid size $m$ and $n$, we define a family of
triangulated tori, denoted by $T_{m,n,\varepsilon}$,  $\varepsilon\in (0,\pi/2)$,  with the following properties:
\begin{enumerate}
\item It has $m n$ vertices, $2mn$ triangles and all vertices have valence 6.
\item All vertices lie on a sphere.
\item The diameter of the `tunnel' of the torus goes to 0 as $\varepsilon \downarrow 0$.
\item The `tunnel' of the torus is triangulated by $2n$ (long and skinny) triangles.
\end{enumerate}
See the first four panels of Figure~\ref{fig:tori}.
In details, it is defined as follows. Let $\Delta u = 2\pi/m$, $\Delta v = 2\pi/n$.
For $i=0,\dots, m-1$, $j = 0,\dots, n-1$, let
\bea
V_{i,j} =
\begin{bmatrix}
\sin(\tau(i \Delta u))\cos(j \Delta v) \\
\sin(\tau(i \Delta u))\sin(j \Delta v) \\
\cos(\tau(i \Delta u))
\end{bmatrix},  \; \mbox{ where }
\tau(u) = \displaystyle\frac{(\pi-2\varepsilon)u}{2\pi(1-1/m)}+\varepsilon.
\eea
Note that $\tau$ maps $[0,2\pi(1-{1}/{m})]$ to $[\varepsilon,\pi-\varepsilon]$.
These $m \cdot n$ points on the unit sphere are
the vertices of $T_{m,n,\varepsilon}$. For the triangulation,
each $V_{i,j}$ is connected to the six neighbors
$V_{i+1,j}$, $V_{i+1,j+1}$, $V_{i,j+1}$, $V_{i-1,j}$, $V_{i-1,j-1}$, $V_{i,j-1}$,
forming also six triangular faces incident on the vertex.
In above, the $+$/$-$ in the first and second indices are modulo $m$ and modulo $n$, respectively.

Note that the $2n$ vertices $V_{0,j}$ and $V_{m-1,j}$, $j = 0, \dots n-1$ are the vertices closest to the north and south pole, respectively.
While these two group of vertices are the furthest apart geometrically, they are connected and form a vertical tunnel of the
torus with $2n$ long skinny triangles.

Denote by $\mathcal{T}_{m,n}$ the set of all PL surfaces with connectivity of a $(m,n)$-regularly triangulated torus and, following
\cite{Bobenko05aconformal}, for any discrete Willmore energy $W_{\rm PL}$ for PL surfaces, write
$$
W_{\rm PL}(\mathcal{T}_{m,n}) = \inf_{S \in \mathcal{T}_{m,n}} W_{\rm PL}(S).
$$

\begin{proposition} \normalfont \label{prop:ST}
For any grid size $m,n \geq 3$, $W_{\rm Bobenko}(T_{m,n,\varepsilon})$ decreases monotonically to $4\pi$ as $\varepsilon \downarrow 0$.
Hence $W_{\rm Bobenko}(\mathcal{T}_{m,n}) \leq 4\pi$ regardless of $m$, $n$.
\end{proposition}

\pf We divide the proof into three steps.

\noindent
$1^\circ$ For any triangulated torus $T_{m,n,\varepsilon}$ defined above,
\bea
W_{\rm Bobenko}(T_{m,n,\varepsilon}) = \frac{1}{2} \sum_{i=0}^{m-1}\sum_{j=0}^{n-1}W_{i,j}, \;\;
\mbox{ where } W_{i,j} = \sum_{(i',j')\sim (i,j)} \beta (i',j') -2\pi
\eea
and $\beta(i',j')$ is an angle formed by the circumscribed
circles of the two triangles sharing the edge connecting $(i, j)$ to $(i',j')$
\cite[Definition 1]{Bobenko05aconformal}. We then notice that
$$
W_{i,j} = 0, \; \forall \; i \neq 0, \, m-1,
$$
since the six vertices around $(i,j)$ form a convex neighborhood lying on a common sphere \cite[Proposition 1]{Bobenko05aconformal}. By symmetry,
$W_{0,j}$ and $W_{m-1,j}$ share the same value for all $j=0,\ldots,n-1$. So we have
$$
W_{\rm Bobenko}(T_{m,n,\varepsilon}) = n\, W_{0,0}.
$$
It remains to show that $W_{0,0}$ decreases monotonicity to $4\pi/n$ as $\varepsilon \downarrow 0$.

\noindent
$2^\circ$ To simplify computation, we take advantage of the M\"obius invariance of $W_{\rm Bobenko}$ by
applying a sphere inversion that maps the unit sphere to the $z=0$ plane.
Specifically, we invert about the sphere with radius $\sqrt{2}$ and centered at the north pole $[0,0,1]^T$ of the unit sphere.
This maps  the south pole of the unit sphere to the origin, and the north pole to infinity;
and it turns our spherical torus $T_{m,n,\varepsilon}$ to a
flattened torus. See Figure~\ref{fig:tori}.

Among the 7 vertices
$V_{0,0}$, $V_{0,1}$, $V_{-1,0}$, $V_{-1,-1}$, $V_{0,-1}$, $V_{1,0}$, $V_{1,1}$
contributing to $W_{0,0}$,
$V_{0,0}$, $V_{0,1}$, $V_{0,-1}$ are close to the north pole, they are sphere inverted to points far away from the origin;
$V_{-1,0}$, $V_{-1,-1}$ are close to the south pole, they are mapped to points close to the origin. The last two neighbors
$V_{1,0}$, $V_{1,1}$ are at an approximately constant distance
from the north pole: their common polar angle is uniformly larger than, and approaches, $\pi/(m-1)$  as $\varepsilon \goto 0$.
Therefore, these 7 vertices
are mapped to $V_i' \in \bR^2$, $i=0,1,\ldots,6$ (in the same cyclic order) with the form
$$
V_0'=\rho_1(\varepsilon) [1, 0], \quad V_1'= \rho_1(\varepsilon) [c,s], \quad V_2' = \rho_2(\varepsilon) [1, 0], \quad
V_3' = \rho_2(\varepsilon) [c, -s],
$$
$$
V_4' = \rho_1(\varepsilon) [c, -s], \quad
V_5' = \rho_3(\varepsilon) [1, 0], \quad
V_6' = \rho_3(\varepsilon) [c, s], \quad c=\cos(2\pi/n), \;\;\; s = \sin(2\pi/n),
$$
where
$
\underbrace{\rho_1(\varepsilon)}_{=\omega(1)} \gg \underbrace{\rho_3(\varepsilon)}_{=\Theta(1)} \gg \underbrace{\rho_2(\varepsilon)}_{=o(1)}$, as $\varepsilon \goto 0$.
By scale invariance, $W_{0,0}$ depends on $\varepsilon$ through\vspace{-5pt}
\bea \label{eq:eps_i}
\epsilon_1 := \rho_3(\varepsilon)/\rho_1(\varepsilon), \mbox{ and }
\epsilon_2: = \rho_2(\varepsilon)/\rho_1(\varepsilon).
\eea
Clearly, $\rho_1$ (resp. $\rho_2$) increases (resp. decreases) as $\varepsilon$ decreases.
So $\epsilon_2(\varepsilon)$ is monotonic increasing for $\varepsilon \in (0,\pi/2)$.

Below, we only need the facts that
$1> \epsilon_1 > \epsilon_2 > 0$, $\epsilon_2$ is monotone in $\varepsilon$ and
$\epsilon_1, \epsilon_2 \goto 0$ as $\varepsilon\goto 0$.

\begin{figure}[h]
\centerline{
\begin{tabular}{cc}
\includegraphics[height=2.3in]{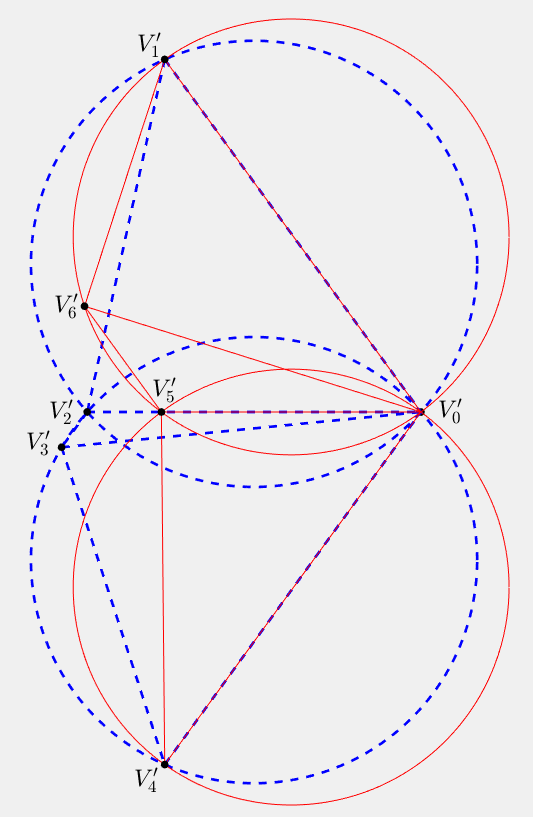} &
 \begin{minipage}[c]{.6\linewidth} \vspace{-180pt}
\caption{Computation of $W_{0,0}$ based on $V_0'= [1, 0]$, $V_1'= [c,s]$, $V_2' = \epsilon_2 [1, 0]$,
$V_3' = \epsilon_2 [c, -s]$, $V_4' = [c, -s]$, $V_5' = \epsilon_1 [1, 0]$, $V_6' = \epsilon_1 [c, s]$.
Among the six circumcircles, two pairs coincide due to the isosceles quadrilaterals  $V_0'V_2'V_3'V_4'$
and $V_0'V_5'V_6'V_1'$, thus only 4 circles are seen; they are also
divided into two groups of three (or two rather), displayed in solid and dashed line styles, which are oriented differently when viewed from the outside of the plane.
This is caused by the `folding' of the neighborhood of $V_0'$, and is also why $W_{0,0}$ does not vanish despite all the vertices are co-planar.
 In this figure, $(c,s)=(\cos(2\pi/5),\sin(2\pi/5))$, i.e. $n=5$,
$(\epsilon_1,\epsilon_2)=(.3,.1)$ and it does not
come from an $\varepsilon$ via \eqref{eq:eps_i}.
Steps $2^\circ$ and $3^\circ$ of the proof also do not rely explicitly on \eqref{eq:eps_i}.
} \label{fig:W00}
\end{minipage}
 \end{tabular}}
\end{figure}
\noindent
$3^\circ$ To calculate $W_{0,0} = \sum_{i=1}^6 \beta_i - 2\pi$ we analyze each angle $\beta_i$ between the circumcircles
of the two triangles sharing the edge $e_i=\overline{V_0'V_i'}$, now thought of as functions of $(\epsilon_1, \epsilon_2)$.

We first notice that the two triangles $V_0'V_3'V_2'$, $V_0'V_3'V_4'$ sharing $e_3$ are co-cyclic because
$V_0'V_2'V_3'V_4'$ form an isosceles quadrilateral. Therefore $\beta_3$ is either $0$ or $\pi$.
A closer inspection based on orientation (or simply applying the formula below) tells us that $\beta_3=0$.
Similarly, $\beta_6=0$.
Next, we show that
\bea \label{eq:limits}
\lim_{\epsilon_1, \epsilon_2\goto 0} \beta_1 = \pi = \lim_{\epsilon_1, \epsilon_2\goto 0} \beta_4, \quad
\lim_{\epsilon_1, \epsilon_2\goto 0} \beta_2 = \frac{2\pi}{n} = \lim_{\epsilon_1, \epsilon_2\goto 0} \beta_5.
\eea
These limits are not hard to see geometrically based on Figure~\ref{fig:W00}; but to get the finer monotonicity result
we resort to algebra and use the formula
$$
\cos(\beta_i)= \frac{\langle A,C \rangle \langle B,D \rangle - \langle A,B \rangle \langle C,D \rangle -\langle B,C \rangle \langle D,A \rangle}{\|A\| \|B\| \|C\| \|D\|},
$$
where $A =V_0'-V_{i+1}'$,  $B = V_{i-1}'-V_0'$, $C = V_{i}'-V_{i-1}'$, $D = V_{i+1}'-V_{i}'$.
(Here the $+$ and $-$ are modulo 6 addition and subtraction operated on the indices $1,\ldots,6$.)
By computation, we get
$
\cos(\beta_1)=\cos(\beta_4)=-\frac{1-\epsilon_1 c - \epsilon_2 c + \epsilon_1\epsilon_2}{\sqrt{1-2\epsilon_1 c + \epsilon_1^2} \sqrt{1-\epsilon_2 c + \epsilon_2^2}}$,
$\cos(\beta_2)= \frac{c-2\epsilon_1+c\epsilon_1^2}{1-2\epsilon_1 c + \epsilon_1^2}$, $\cos(\beta_5)=\frac{c-2\epsilon_2+c\epsilon_2^2}{1-\epsilon_2 c + \epsilon_2^2}$. The limits \eqref{eq:limits} then follow immediately
and we have proved that $\lim_{\epsilon_1, \epsilon_2\goto 0}
W_{0,0}(\epsilon_1,\epsilon_2) = 4\pi/n$. To see that the convergence is monotone in the original $\varepsilon$, write
\begin{align*}
W_{0,0}(\epsilon_1,\epsilon_2) = 
 \cos^{-1} \left(\frac{c-2\epsilon_1+c\epsilon_1^2}{1-2\epsilon_1 c + \epsilon_1^2} \right)
+ \cos^{-1} \left(\frac{c-2\epsilon_2+c\epsilon_2^2}{1-2\epsilon_2 c + \epsilon_2^2} \right) - 2
\cos^{-1} \left( \frac{1-\epsilon_1 c - \epsilon_2 c + \epsilon_1\epsilon_2}{\sqrt{1-2\epsilon_1 c + \epsilon_1^2} \sqrt{1-2\epsilon_2 c + \epsilon_2^2}}\right).
\end{align*}
By either a geometric argument or explicitly checking that
$\frac{\partial W_{0,0}}{\partial \epsilon_1} = 0$
when $1 \geq \epsilon_1 \geq \epsilon_2 \geq 0$,
$W_{0,0}(\epsilon_1,\epsilon_2)$ is independent of $\epsilon_1$ and
$W_{0,0}= 2\cos^{-1} \left(\frac{c-2\epsilon_2+c\epsilon_2^2}{1-2\epsilon_2 c + \epsilon_2^2} \right)$.
Then again by either a geometric argument or explicitly checking that
$\frac{\partial W_{0,0}}{\partial \epsilon_2} = \frac{4s}{1 - 2 \epsilon_2 c + \epsilon_2^2} > 0$,
$W_{0,0}$ is monotone in
$\epsilon_2$.
Combined with the monotonicity
of $\epsilon_2(\varepsilon)$, the proof is completed.
\eop

In the genus $g=0$ case, as long as the connectivity is generic in some sense (see \cite[Proposition 9]{Bobenko05aconformal}),
the simplicial sphere is \emph{inscribable} in a sphere and hence has a
 minimum $W_{\rm Bobenko}$-energy $4\pi$ (0 in Bobenko's definition of $W$ = our definition of $W_{\rm Bobenko} - 4\pi(1-g)$).
This is of course the same energy level as what one expects from the smooth setting.

\begin{figure}[h]
\centerline{
\begin{tabular}{cccc}
\includegraphics[height=1.2in]{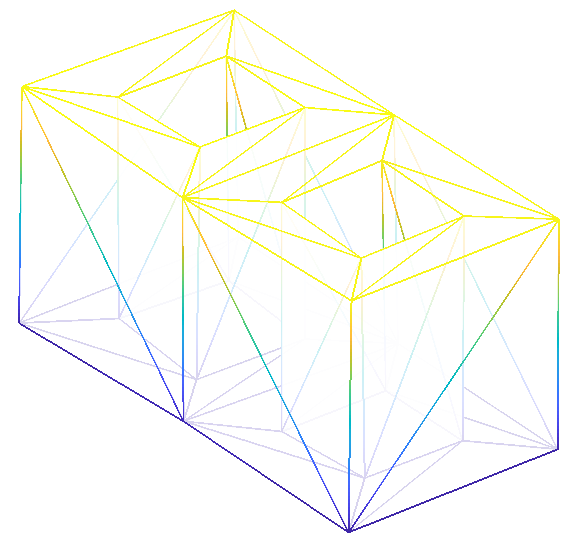}& \hspace{-.1in}
\includegraphics[height=1.2in]{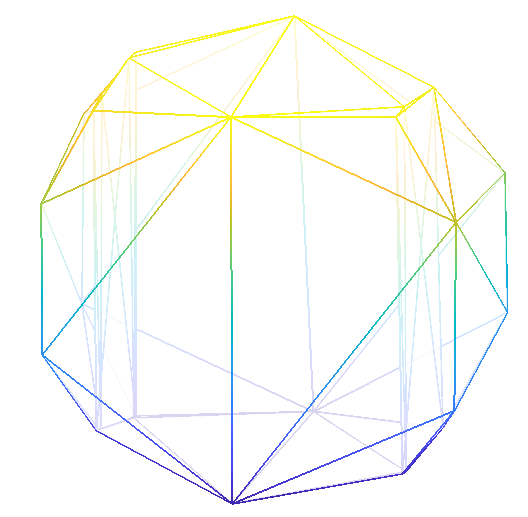}& \hspace{-.1in}
\includegraphics[height=1.2in]{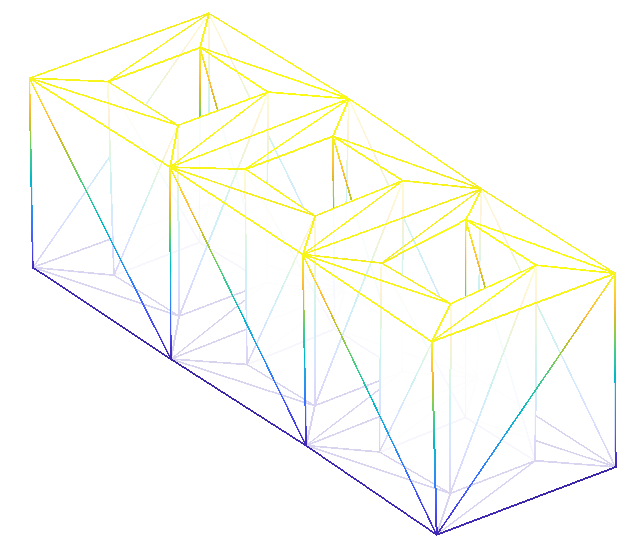}& \hspace{-.1in}
\includegraphics[height=1.2in]{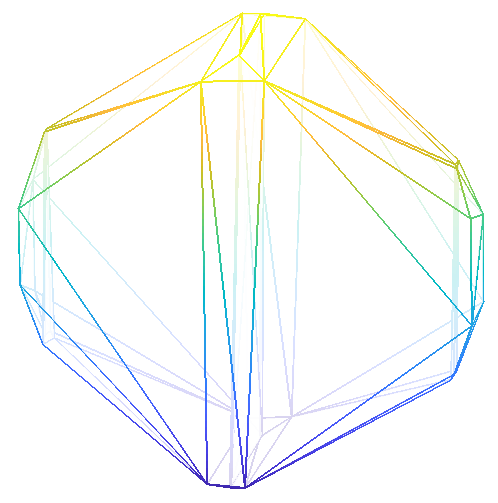}\\
(a) & (b) & (c) & (d)
 \end{tabular}} 
 \caption{
 A $W_{\rm Bobenko}$-minimization process rounds the 2- and 3-hole tori in (a) and (c) and closes up the holes,
 resulting in the empirical minimizers in (b) \& (d), both with $W_{\rm Bobenko}\approx 4\pi$.}
 \label{fig:HigherGenus}
\end{figure}

Proposition~\ref{prop:ST} shows that in the genus 1 case, the infimum is never what one expects from
the smooth setting, even with the perfectly regular connectivity and regardless of the grid size.
Moreover, Proposition~\ref{prop:ST} and computational experiments suggest the following:
\begin{conjecture} \normalfont \label{conjecture:WBobenko}
For any $m, n \geq 3$, $W_{\rm Bobenko}(\mathcal{T}_{m,n}) = 4\pi$. Moreover, the infimum is realized by degenerate spheres ($T_{m,n,\varepsilon}$,
$\varepsilon\goto 0$) and {\it not} by any embedded genus 1 PL surface.
\end{conjecture}
This conjecture suggests that a genus 1 $W_{\rm Bobenko}$-minimizer does {\it not} exist, again in
contrast to the smooth setting; cf. \cite{Simon:Willmore}.
It also seems possible to formulate and prove a generalization of Proposition~\ref{prop:ST} to genus $g>1$, as suggested by the computational
result in Figure~\ref{fig:HigherGenus}.

Back to the genus 1 setting, the proof of Proposition~\ref{prop:ST} suggests the following generalization. Consider the family of planar tori, denoted
by $$T_{m,n,\mathbf{r}},$$
with the same connectivity as before and vertices
\begin{align*}
V_{i,j} &= r_i \left[ \cos\left(\frac{2\pi j}{n} \right), \sin\left(\frac{2\pi j}{n} \right), 0 \right], \;\; i=0,\ldots,m-1, \;\; j=0,\ldots,n-1, \\
\mathbf{r} &= (r_0, r_1, \ldots, r_{m-1}), \quad r_0 > r_1 > \cdots > r_{m-1} >0.
\end{align*}
We have the following corollary of the proof.
\begin{corollary} \normalfont
For any grid size $m,n\geq 3$, and fixed $r_0> \cdots > r_{m-2}>0$,
$$\lim_{r_{m-1}\goto 0} W_{\rm Bobenko}(T_{m,n,\mathbf{r}})=4\pi$$ and the convergence is monotone.
\end{corollary}
\pf 
All vertices
in the intermediate layers ($0<i<m-1$) have zero energy. By the calculation in the previous proof,
every vertex in the outermost ($i=0$) layer has the same energy $2\cos^{-1} \left(\frac{c-2\epsilon+c\epsilon^2}{1-2\epsilon c + \epsilon^2} \right)$, $\epsilon=r_{m-1}/r_0$. The fact that this energy is independent of $r_1$, together with the M\"obius invariance of
each $W_{i,j}$, actually imply that every vertex in the innermost layer ($i=m-1$) also has the same energy $2\cos^{-1} \left(\frac{c-2\epsilon+c\epsilon^2}{1-2\epsilon c + \epsilon^2} \right)$ (independent of $r_{m-2}$); this can be seen by turning $T_{m,n,\mathbf{r}}$ inside out based on
 inverting about the sphere with radius $\sqrt{r_0 r_{m-1}}$ centered at the origin.
Therefore $W_{\rm Bobenko}(T_{m,n,\mathbf{r}}) = \frac{1}{2} \cdot 2 n \cdot  2\cos^{-1} \left(\frac{c-2\epsilon+c\epsilon^2}{1-2\epsilon c + \epsilon^2} \right)$ and the result follows by taking $\epsilon \goto 0$.
\eop

\subsubsection{$W_{\rm Centroid}$ and $W_{\rm EffAreaCur}$} \label{sec:W_PL}
For any mesh $T_{m,n,\mathbf{r}}$ with $r_0$ normalized to 1, $r_0=1>r_1>\cdots>r_{m-1}>0$,
we have the following by direct calculation: 
\begin{align} \label{eq:T_cal}
\begin{split}
&\mbox{$a_{\rm centroid}(0,0) = \frac{4- r_1 - r_1^2 - r_{m-1} - r_{m-1}^2 }{6}s$}, \quad
\mbox{$a_{\rm centroid}(m-1,0) = \frac{1 + r_{m-1} - 4 r_{m-1}^2 + r_{m-2}^2 + r_{m-1} r_{m-2}}{6}s$} \\
& \mbox{$\nabla_{(0,0)} A = \begin{bmatrix}
                    2s, & 0, & 0
                  \end{bmatrix}$},
\;\;\;
\mbox{$\nabla_{(m-1,0)} A = \begin{bmatrix}
                    -2 r_{m-1}s, & 0, & 0
                  \end{bmatrix}$},
\\
&\mbox{$\nabla_{(0,0)} V = \begin{bmatrix}
                    0, & 0, &
                    \frac{r_1 - r_{m-1} + r_1^2 - r_{m-1}^2}{6} s
                  \end{bmatrix}$},
\;\;\;
\mbox{$\nabla_{(m-1,0)} V = \begin{bmatrix}
                    0, & 0, &
                    \frac{(r_{m-2} - 1)(r_{m-1} + r_{m-2} + 1)}{6} s
                  \end{bmatrix}$},
\end{split}
\end{align}
where $s = \sin(2\pi/n)$. Since the mesh is flat, it does not have an `inside' or `outside',
we arbitrarily choose a consistent orientation in order to determine the direction of $\nabla_v V$.
(Recall Remark~\ref{remark:orientation}.)

Next, we show:
\begin{proposition} \normalfont \label{prop:W_negative}
Let the grid sizes $m,n \geq 3$ be fixed.
\begin{itemize}
\item[(i)] For any fixed $r_0>0$,
\begin{align}
\lim_{r_1 \goto 0} W_{\rm Centroid}(T_{m,n,\mathbf{r}}) = \frac{3}{2} n \sin\Big( \frac{2\pi}{n} \Big).
\end{align}
Hence
$W_{\rm Centroid}(\mathcal{T}_{m,n}) \leq 3\pi < 2\pi^2$ regardless of $m$, $n$.
\item[(ii)] For any fixed $r_0 > r_2 > \cdots > r_{m-2} >0$,
\begin{align} \label{eq:limit2}
\lim_{\substack{r_1\goto r_0 \\ r_{m-1} \goto 0}} W_{\rm EffAreaCur}(T_{m,n,\mathbf{r}})= 3 n \sin\Big( \frac{2\pi}{n} \Big).
\end{align}
Hence $W_{\rm EffAreaCur}(\mathcal{T}_{m,n}) \leq 6\pi < 2\pi^2$ regardless of $m$, $n$.
\end{itemize}
\end{proposition}

\pf
All vertices
in the intermediate layers ($0<i<m-1$) have zero energy. By symmetry, every vertex in the outermost ($i=0$) and innermost ($i=m-1$) layer has  the same energy.
So $W_{\rm Centroid}(T_{m,n,\mathbf{r}})= n (W_{0,0} + W_{m-1,0})$.
By scale invariance of $W_{\rm Centroid}$, we can assume $r_0=1$.
By \eqref{eq:T_cal},
\begin{align*}
\lim_{r_1, r_{m-1} \goto 0} W_{0,0} &= \lim_{r_1, r_{m-1}\goto 0} \left[\frac{\| \nabla_{(0,0)} A \|}{2 a_{\rm centroid}(0,0)} \right]^2 a_{\rm centroid}(0,0) = \frac{3}{2} \sin(2\pi/n), \\
\lim_{r_1, r_{m-1} \goto 0} W_{m-1,0} &= \lim_{r_1, r_{m-1}\goto 0} \left[\frac{\| \nabla_{(m-1,0)} A \|}{2 a_{\rm centroid}(m-1,0)} \right]^2 a_{\rm centroid}(m-1,0) = 0.
\end{align*}
And the proof of (i) is completed.

Similarly, $W_{\rm EffAreaCur}(T_{m,n,\mathbf{r}})= n (W_{0,0} + W_{m-1,0})$.
Assume $r_0=1$. By \eqref{eq:T_cal},
\begin{align*}
\lim_{\substack{r_1\goto 1 \\ r_{m-1} \goto 0}} W_{0,0} &= \lim_{\substack{r_1\goto 1 \\ r_{m-1} \goto 0}} \left[\frac{\| \nabla_{(0,0)} A \|}{2 \| \nabla_{(0,0)} V \|} \right]^2 \| \nabla_{(0,0)} V \| = 3 \sin(2\pi/n), \\
\lim_{\substack{r_1\goto 1 \\ r_{m-1} \goto 0}} W_{m-1,0} &= \lim_{\substack{r_1\goto 1 \\ r_{m-1} \goto 0}} \left[\frac{\| \nabla_{(m-1,0)} A \|}{2 \| \nabla_{(m-1,0)} V \|} \right]^2 \| \nabla_{(m-1,0)} V \| = 0.
\end{align*}
And the proof is completed.
\eop

In contrast to Conjecture~\ref{conjecture:WBobenko} pertaining to $W_{\rm Bobenko}$, we observe from computation and preliminary
calculations that
Proposition~\ref{prop:W_negative} is {\it not} sharp in the sense that
 the infimum values $W_{\rm Centroid}(\mathcal{T}_{m,n})$ and $W_{\rm EffAreaCur}(\mathcal{T}_{m,n})$
are strictly smaller than $3/2n\sin(2\pi/n)$ and $3 n\sin(2\pi/n)$, respectively. It is observed that
minimizers of both energies are some subtly `folded up' planar meshes.

Although we do not pursue it here, a similar negative result holds for $W_{\rm Voronoi}$.

\begin{figure}[h]
\centerline{
\begin{tabular}{cccc}
\includegraphics[height=1.2in]{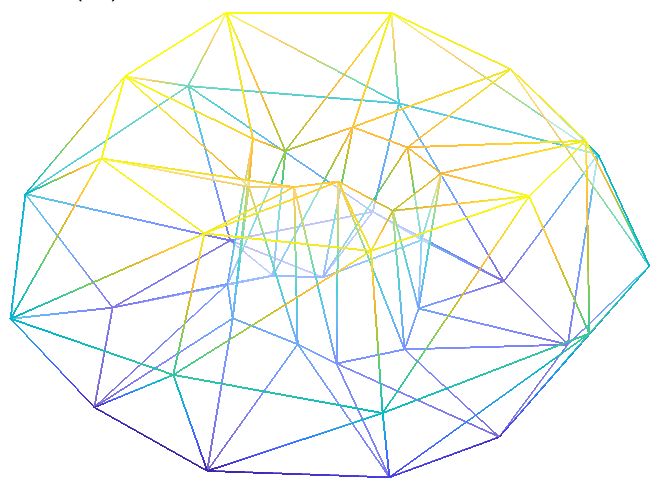}& \hspace{-.2in}
\includegraphics[height=1.0in]{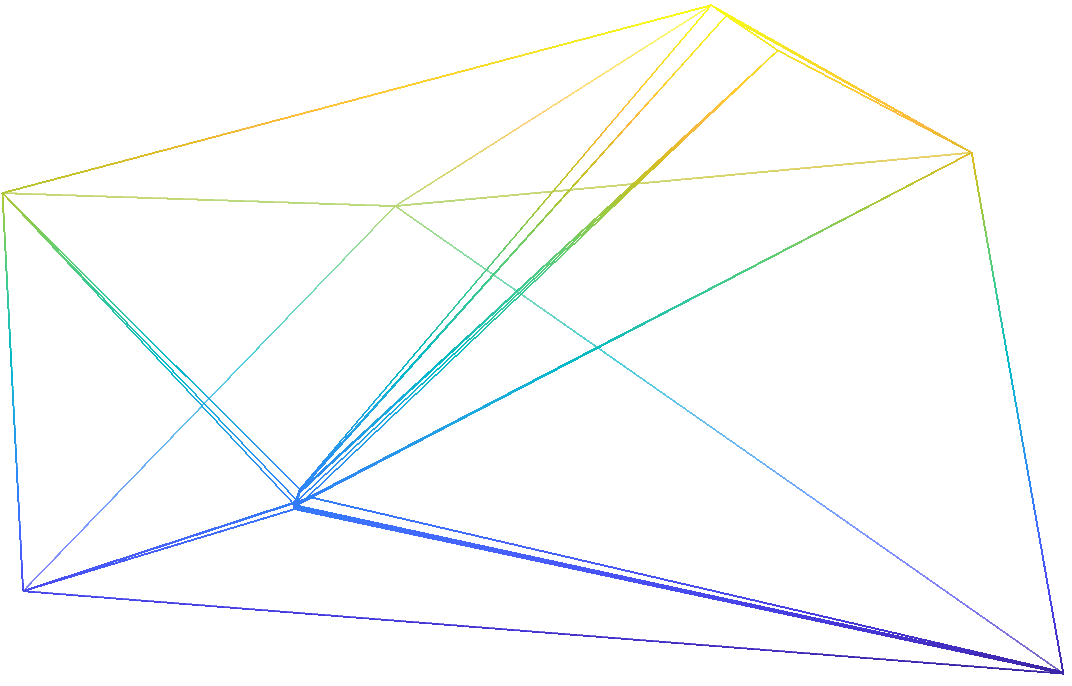}& \hspace{-.2in}
\includegraphics[height=1.2in]{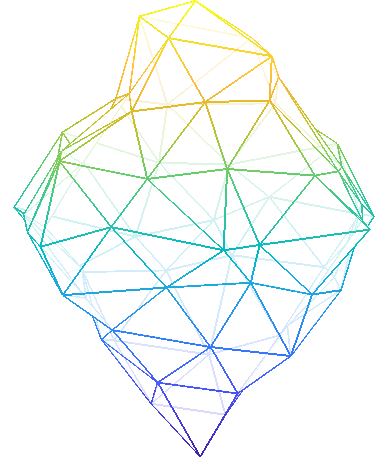}& \hspace{-.2in}
\includegraphics[height=1.2in]{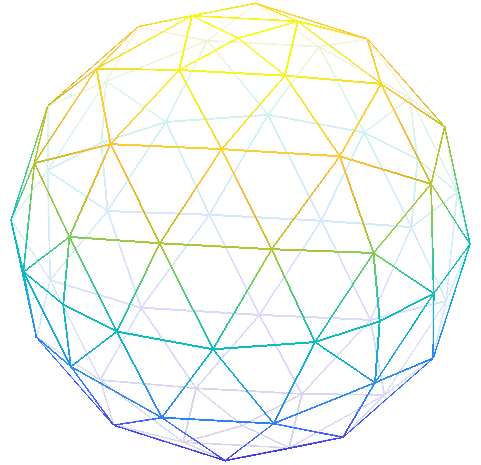}\\
(a) $W_{\rm NormalCur}=37.73$ &\hspace{-.2in} (b) $W_{\rm NormalCur}=14.85$  & \hspace{-.2in}
(c) $W_{\rm NormalCur}=45.46$ &\hspace{-.1in} (d) $W_{\rm NormalCur}=12.95$\\
& $\ll 2\pi^2$ & & $\approx 4\pi$
 \end{tabular}}
 \caption{Failure for genus 1 and success for genus 0. (a)-(b): a $W_{\rm NormalCur}$-minimization process `folds up'
the regularly triangulated torus with grid size $(m, n) = (5, 8)$ into a mesh that fails to approximate a Clifford
torus. (c)-(d): a $W_{\rm NormalCur}$-minimization process rounds out a twice-subdivided octahedron into an
approximation of a round sphere.}
 \label{fig:NormalCurMinimizer58}
\end{figure}

The situation for $W_{\rm NormalCur}$, however, is less well-understood.
Note that the outer- and inner-most vertices of $T_{m,n,\mathbf{r}}$
satisfy $$
\nabla_v V \bot \nabla_v A,  \;\;\; \nabla_v A \neq 0;
$$
recall \eqref{eq:T_cal}.
This means the corresponding local areas in $W_{\rm NormalCur}$ vanish and $W_{\rm NormalCur}=+\infty$.
Therefore the planar meshes we considered cannot make $W_{\rm NormalCur}$ small.
This is, by design, a key difference between
$W_{\rm EffAreaCur}$ and $W_{\rm NormalCur}$.
However, this does not mean $W_{\rm NormalCur}$ is rid of the type of negative results in
Propsoition~\ref{prop:ST} and \ref{prop:W_negative}. We observe from computations that minimizers
of $W_{\rm NormalCur}$ exhibit a similar behavior as those of $W_{\rm EffAreaCur}$: they are highly
`folded up', but non-planar, meshes, with energies less than $2\pi^2$; see Figure~\ref{fig:NormalCurMinimizer58}(a)-(b).

\subsubsection{Genus 0 case}
In the genus 0 case, computational experiments suggest that $W_{\rm Centroid}$, $W_{\rm Voronoi}$
and $W_{\rm EffAreaCur}$ fail in a way similar
to the genus 1 case; we believe that a genus 0 version of Proposition~\ref{prop:W_negative} can be established.

However, the situation for $W_{\rm Bobenko}$ and $W_{\rm NormalCur}$ are different.
Recall that the Willmore energy of any closed surface should always be greater than or equal to $4\pi$, with
equality attained only by the round sphere \cite{Willmore:WillmoreConjecture}. $W_{\rm Bobenko}$
is designed to satisfy this `ground-state' property:
$W_{\rm Bobenko}\geq 4\pi$, with equality holds when and only when the PL surface is a convex polyhedron inscribed in
a sphere \cite[Theorem 5]{Bobenko05aconformal}. So there cannot be a genus 0 version of the negative result
in Proposition~\ref{prop:ST}.
An empirical finding is that $W_{\rm NormalCur}$ actually works well for the genus 0 Willmore problem. See
Figure~\ref{fig:NormalCurMinimizer58}(c)-(d), which shows one of the many trials of the experiment with a randomized initial mesh; all trials of
$W_{\rm NormalCur}$-minimization
applied to a genus 0 mesh result in a near round sphere with $W_{\rm NormalCur}$ slightly larger than $4\pi$.

\section{A Regularized PL Method} \label{sec:Regularization}
All failures we observed have one thing in common: triangles with bad aspect ratios develop.
This is of course a familiar issue in FEM and mesh generation, but in the moving surface context here
the issue seems different and understudied.
In previous work by Hsu-Kusner-Sullivan  \cite{Kusner:Brakke},  procedures such as `vertex averaging', `edge
notching' and `equiangulation', implemented in the Surface Evolver \cite{Brakke:evolver,Brakke:evolverManual},
are used to fix up triangles with bad aspect ratios along the way of the optimization process.
While such `mesh smoothing'
procedures are probably well studied in the mesh generation community,
when applied to the geometric variational problems here the approach seems ad-hoc and difficult to analyze
mathematically.

Our numerical methods, based on either PL or SS, are parametric in nature. Yet,
the variational problems are geometric and their solutions are independent of parametrization.
This creates another problem instead of solving the existing one!
Our observation is that it is possible to design a method in such a way that the
 two problems cancel each other.

Since the solution is independent of parametrization, we
may request the parametrization to be (approximately)
{\it conformal},\footnote{Not to be confused with the conformality in finite-element methods.} i.e. angle preserving, with respect to an appropriate
conformal structure in the reference manifold $\Sigma$. If such a property can be built into a $W$-minimizing
process, the method would only search over PL surfaces without long skinny triangles.

The approach adopted here is inspired by the
classical uniformization theorem and results in harmonic maps \cite{EellsSampson:Review1,Eells:Review2,Eells:Review3}, as well as
the relatively recent developments in computational conformal geometry, see
\cite{GuYau:Recent,Crane:2020:DCG,IPMI03:gu,Duchamp:PL,Lui2012,Gu-2018-I,Gu-2018-II} and the references therein.

In a finite-dimensional approximation of the solution surface, we cannot expect to have exactly conformal parametrization.
To fix the idea, let us first consider the problem of finding a conformally parameterized
Willmore surface in the smooth setting.

\subsection{Ideas in the smooth setting} \label{sec:Continuous}
We use the theory of harmonic maps. From now on we assume $\Sigma$ is a Riemann surface.
For any map $\mathbf{x} \in W^{1,2}(\Sigma, \bR^3)$,  we can define its Dirichlet energy by
\bea \label{eq:DirichletEnergy}
\mathcal{E}(\mathbf{x})  =  \frac{1}{2} \int_{M} \|d\mathbf{x}\|^2 d\Sigma,
\eea
where $\|d\mathbf{x}\|$ is the Hilbert-Schmidt norm first
defined based on \textit{some} Riemannian metric
on $\Sigma$ and the standard Riemannian metric in $\bR^3$, followed by the observation
that $\mathcal{E}$ is
 invariant under any conformal change of Riemannian metric in the domain;
see \cite[Pg. 126]{EellsSampson:Review1}.
This invariance is specific to $\dim \Sigma = 2$, and it means that
the Dirichlet energy depends only on the conformal structure of $\Sigma$. Next, we have the following well-known result:


\begin{theorem}[\cite{EellsSampson:Review1}] \label{thm:ConformalEnergy}
 For any $C^1$ immersion $\mathbf{x}: \Sigma \goto \bR^n$,
  $A(\mathbf{x}) \leq \mathcal{E}(\mathbf{x})$. 
  Equality holds when and only when $\mathbf{x}$ is conformal.
\end{theorem}

%

It is already evident from the theorem that the difference $\mathcal{E}(\mathbf{x})-A(\mathbf{x})$
measures the deviation of $\bx$ being conformal. In fact, the difference can be explicitly expressed
by a conformal energy expressing the deviation of $\bx$ from satisfying the Cauchy-Riemann equations.

\begin{proposition} \normalfont \label{prop:Wmin_Regularized}
Let $\Sigma$ be a Riemann surface and $S = \{\bx \in \Imm_{C^1(\Sigma,\bR^3) \cap W^{2,2}(\Sigma,\bR^3)}: A(\bx)=A_0 \}$
for a fixed $A_0>0$.
Assume that:
\begin{itemize}
\item[(i)] a minimizer
$\bx_\lambda \in \argmin_{\bx \in S} W(\bx)+\lambda\, \mathcal{E}(\bx)$
exists for all small enough $\lambda$,
\item[(ii)] 
$\lim_{\lambda \downarrow 0} \bx_{\lambda} =: \bx_\ast$
exists. (Here convergence is in the topology of $C^1(\Sigma,\bR^3) \cap W^{2,2}(\Sigma,\bR^3)$.) 
\end{itemize}
Then
\begin{itemize}
\item[(I)] $\mathbf{x}_\ast:\Sigma \goto \bR^3$ is a Willmore minimizer, i.e. it is a
 solution of the Willmore problem
 $\min_{\bx \in S} W(\bx)$. Moreover,
  any other Willmore minimizer $\mathbf{x}_{\ast\ast}$ satisfies $A_0 \leq \mathcal{E}(\mathbf{x}_\ast)\leq \mathcal{E}(\mathbf{x}_{\ast\ast})$.
\item[(II)] If 
the conformal structure on $\Sigma$ is such that there exists a
conformal parametrization of the surface $\mathbf{x}_\ast(\Sigma)$, then $\bx_\ast$, and in fact $\bx_\lambda$ for any $\lambda>0$,
is a conformal parametrization of a Willmore minimizer,
i.e. $W(\bx_\lambda) = W(\bx_\ast)$ and
$A_0 = \mathcal{E}(\mathbf{x}_\lambda) = \mathcal{E}(\mathbf{x}_\ast)$ for any $\lambda>0$.
 \end{itemize}
\end{proposition}
\pf
We first argue that $\bx_\ast$ must be a Willmore minimizer. If not, there exists an
$\bx' \in S$ such that $W(\bx')< W(\bx_\ast)$.
But then for small enough $\lambda>0$, we
must have $W(\bx') + \lambda\, \cE(\bx') < W(\bx_\ast)+ \lambda\, \cE(\bx_\ast)$.
(If $\cE(\bx')\leq \cE(\bx_\ast)$, this holds for any $\lambda>0$, otherwise choose
$\lambda < (\cE(\bx')-\cE(\bx_\ast))/(W(\bx_\ast)-W(\bx'))$.)
By (ii) and the continuity of $W: S \goto \bR$ and $\cE: S \goto \bR$,
we have
$$
W(\bx') + \lambda\, \cE(\bx') < W(\bx_\lambda)+ \lambda\, \cE(\bx_\lambda), \quad \forall \mbox{ small enough $\lambda>0$.}
$$
This contradicts assumption (i), and we have proved the first claim in (I). We now argue the second claim again by contradiction. Let $\bx_{\ast\ast}$ be any other Willmore minimizer,
so $W(\bx_{\ast}) = W(\bx_{\ast\ast})$. Assume the contrary that
$\cE(\bx_{\ast}) > \cE(\bx_{\ast\ast})$.
Then $W(\bx_{\ast}) + \lambda\,\cE(\bx_{\ast}) > W(\bx_{\ast\ast}) + \lambda\,\cE(\bx_{\ast\ast})$ for any $\lambda$.
Again by (ii) and the continuity of $W$ and $\cE$,
$$
W(\bx_\lambda)+ \lambda\, \cE(\bx_\lambda) > W(\bx_{\ast\ast}) + \lambda\,\cE(\bx_{\ast\ast}), \quad \forall \mbox{ small enough $\lambda>0$.}
$$
This contradicts (i), so $\cE(\bx_{\ast}) \leq \cE(\bx_{\ast\ast})$. By Theorem~\ref{thm:ConformalEnergy}, we also have $A_0 \leq \cE(\bx_{\ast}) \leq \cE(\bx_{\ast\ast})$.

By Theorem~\ref{thm:ConformalEnergy} and the assumption of the conformal structure on $\Sigma$,
the minimum values of $W$ and $\cE$ over $S$, namely
$W(\bx_\ast)$ and $A_0$, are attained simultaneously by some $\bx \in S$. But any $\bx_\lambda$, being a minimizer of
$W(\cdot)+\lambda A(\cdot)$, must also be such a joint minimizer. Then, by continuity, $\bx_\ast$ is such a joint minimizer as well.
\eop

\gap
\begin{remark} \label{remark:RegularizationTheorem}
In the most ubiquitous spherical topology, there is
only one conformal structure up to conformal equivalence.
Together with the genus 0 case of the uniformization theorem, the assumption in (II) above is satisfied automatically
when $\Sigma$ is the Riemann sphere. In this case, solving the $\cE$-penalized
Willmore problem with \emph{any} penalization parameter $\lambda$ -- not necessarily small -- would deliver a conformally parameterized Willmore surface.
In the higher genus cases, there is a 1 (when $g=1$) and $3g -3$ (when $g\geq 2$) complex-dimensional space of
conformal structures \cite{jost2006riemann}, and typically
we do not know a priori the ideal conformal structure for the Willmore surface.
With the correct conformal structure, a conformally parametrized solution is impossible, but
Proposition~\ref{prop:Wmin_Regularized}
says (in (I)) that the solution of the $\cE$-penalized
Willmore problem with a \emph{small} penalization parameter $\lambda$ would still deliver a Willmore surface parameterized
as conformally as possible in the non-ideal conformal structure.
\end{remark}

We shall use this result to guide us to
develop an algorithm for fixing the naive PL method. Before we proceed, we present a version of Theorem~\ref{thm:ConformalEnergy}
for PL surfaces. 

\subsection{A PL version of Theorem~\ref{thm:ConformalEnergy}} \label{sec:PLConformalEnergy}
Theorem~\ref{thm:ConformalEnergy_PL} below is well-known to experts in computational conformal geometry;
see, for example, \cite{Crane:2020:DCG,GuYau:Recent} and references therein.
We provide a proof for it not only for the sake of self-containedness but also to help us to understand its subtle
relationship with the smooth counterpart Theorem~\ref{thm:ConformalEnergy}.

\begin{definition} \label{def:ConformalPL}
Let $K$ be a simplicial surface.
\begin{enumerate}
\item A \emph{angle assignment} on $K$ is an assignment of angles $(\theta_1^\tau, \theta_2^\tau, \theta_3^\tau)$
to every triangle $\tau$ in $K$.
\item A PL immersion $\bx: |K| \goto \bR^n$ is called \emph{angle preserving} if every triangle $\tau$ in $K$ is similar to $\bx(\tau) \in \bR^n$,
i.e. they share the same angles.
\end{enumerate}
\end{definition}

The notion of angle assignment in 1. above actually defines a conformal (or complex analytic) atlas in the traditional sense via a construction
by L. Bers. The construction
essentially exploits the fact that the map
$$
z \mapsto z^a
$$
maps the interior -- but not the boundary -- of
a cone in $\bC$ with apex angle $\theta$ bihorlomorphically to the interior of a cone with apex angle $a \theta$ (assume $0<a< 2\pi/\theta$).
For details, see \cite[Section 2]{Duchamp:PL} or \cite[Lecture 2]{Bers:Riemann}.
In these references, the conformal structure is supposed to be induced by a PL embedding
$F: |K| \goto \bR^n$. However, it is obvious that the conformal atlas
defined there
depends only on the angle assignment induced by the embedding.
We note that, however, many different angle assignments may give rise
to {\it equivalent} conformal structures on $K$. For example, in the genus 0 case, any two conformal structures are conformally equivalent.

Once such a conformal structure
is defined, and noting also
that piecewise linear maps have
square integrable derivatives,
the Dirichlet energy $\cE(\bx)$ of a PL immersion $\bx$ is well-defined according to \eqref{eq:DirichletEnergy}.
Of course, the area $A(\bx)$
is also well-defined in the standard sense.
As such, the first half of Theorem~\ref{thm:ConformalEnergy}, which requires only $\bx$ to be in $W^{1,2}$,
holds verbatim, i.e. $\cE(\bx) \leq A(\bx)$.
However, the second half of Theorem~\ref{thm:ConformalEnergy} is less obvious, as the
standard definition of conformality for smooth maps does not seem to apply directly to PL maps.
Fortunately, if one simply defines conformality for PL immersions as in Definition~\ref{def:ConformalPL}, we get the following result:
\begin{theorem}[PL version of Theorem~\ref{thm:ConformalEnergy}] \label{thm:ConformalEnergy_PL}
Let $K$ be a simplicial surface with an angle assignment.
For any PL immersion $\mathbf{x}: |K| \goto \bR^n$,
  $A(\mathbf{x}) \leq \mathcal{E}(\mathbf{x})$. 
  Equality holds when and only when $\mathbf{x}$ is angle preserving.
\end{theorem}


\pf Bers' conformal structure and the fact that $\cE$ depends solely on the conformal structure
guarantee that in the interior of every triangle $\tau$ of $K$, the
 Dirichlet energy \eqref{eq:DirichletEnergy} can be computed as
\bea
\cE(\bx)  = \frac{1}{2} \sum_{\tau \in K} \int_\tau \|d\bx\|^2 d{\rm Area}_\tau,
\eea
where
``$d{\rm Area}_\tau$'' is the area element in the coordinates of $\tau$ after we identify $\tau$ with a triangle on the plane
 with the angles assigned to $\tau$. For details, consult \cite[Section 2]{Duchamp:PL}.
For any $\tau$, $\bx|_\tau$ is a linear map and we may represent
it by a rank 2 matrix $A_\tau \in \bR^{n\times 2}$, so
\bea \label{eq:E_PL}
\cE(\bx)  = \frac{1}{2} \sum_{\tau} {\rm trace}(A_\tau^T A_\tau) {\rm Area}(\tau)
=
\sum_{\tau} \frac{\frac{1}{2} {\rm trace}(A_\tau^T A_\tau)}{\sqrt{\det(A_\tau^T A_\tau)}} {\rm Area}(\bx(\tau)).
\eea
Since
$$
\frac{\frac{1}{2} {\rm trace}(A_\tau^T A_\tau)}{\sqrt{\det(A_\tau^T A_\tau)}} =
\frac{\frac{1}{2} (\lambda^\tau_1+\lambda^\tau_2)}{\sqrt{\lambda^\tau_1 \lambda^\tau_2}} \geq 1,
$$
where $\lambda_{\tau,1}$ and $\lambda_{\tau,2}$, both positive, are the eigenvalues of $A_\tau^T A_\tau$, $\cE(\bx) \geq A(\bx)$. Moreover, equality holds above if and only if
$\lambda^\tau_1=\lambda^\tau_2$ for all triangles $\tau \in K$, which is exactly the case when $\bx$ is angle preserving. \eop

In the PL setting above, the Dirichlet energy can also be expressed by the cotangent formula
\bea \label{eq:cotangent_E_PL}
\cE(\bx) = \frac{1}{4} \sum_{\tau} \sum_{i=1,2,3} \cot(\theta^\tau_i)\, {\rm length}(\bx(e_i^\tau))^2
= \frac{1}{4} \sum_{e \in {\rm Edge}} (\cot(\alpha_e)+\cot(\beta_e))\, {\rm length}(\bx(e))^2,
\eea
where $e_i^\tau$ is the edge of $\tau$ opposite the angle $\theta^\tau_i$, and $\alpha_e$ and $\beta_e$ are the angles opposite the edge $e$
in the two incident triangles.
This follows immediately from a simple formula:
if $f:\tau \goto \bR$ is a linear function on a triangle $\tau \subset \bR^2$,
then
$$
\iint_{\tau} \left|\frac{\partial{f}}{\partial x} \right|^2 +
\left| \frac{\partial{f}}{\partial y} \right|^2 dx\, dy =
\frac{1}{2} \Big\{ \cot(\theta_1) |f(v_2)-f(v_3)|^2 +
\cot(\theta_2) |f(v_3)-f(v_1)|^2 +
\cot(\theta_3) |f(v_1)-f(v_2)|^2 \Big\},
$$
where $v_1,v_2,v_3 \in \bR^2$ are the vertices of $\tau$ and $\theta_i$ is the angle at vertex $v_i$.

\begin{remark} 
A subtle point here is that if
we simply \emph{define} the Dirichlet energy of a PL immersion by either
\eqref{eq:E_PL} or \eqref{eq:cotangent_E_PL}, then
both the formulation and the proof of Theorem~\ref{thm:ConformalEnergy_PL} completely
bypass Bers' construction. We may
then view the `Bers-free version' of Theorem~\ref{thm:ConformalEnergy_PL} as a {\it discrete analog}
 of Theorem~\ref{thm:ConformalEnergy}. This would be in line with the spirit of `discrete differential geometry' \cite{Bobenko:Book}.
We find it instructive to know that Theorem~\ref{thm:ConformalEnergy_PL} and Theorem~\ref{thm:ConformalEnergy}
 are not only analogous, but that one may view (at least the first part of) the former theorem as a {\it special case}
 of the latter classical result.
\end{remark}

Finally, we note that when the angles $(\theta_1^\tau, \theta_2^\tau, \theta_3^\tau)$ assigned to $K$ are
exactly those of the triangles $\bx(\tau)$ in $\bR^n$, then it is a tautology to say that $\bx$ is angle preserving.
 In this case, $A(\bx)$ is given by the cotangent formula \eqref{eq:cotangent_E_PL}. This, in turn,
can be used to give an alternative proof for the cotangent formula of $\nabla_v A$ in \eqref{eq:cotangent}, explaining
also why the two seemingly unrelated cotangent formulas look so alike.

\subsection{A penalized PL method}
Inspired by Proposition~\ref{prop:Wmin_Regularized},
we propose a numerical solution to the Willmore, Canham or Helfrich problem based on solving:
\bea \label{eq:PeanalizedWPLmin}
\min_{\cV} W_{\rm PL}(\cV) + \lambda\: \mathcal{E}(\cV), \mbox{ s.t. $A(\cV) = A_0$ and
relevant constraints on $V(\cV)$ and $M(\cV)$.}
\eea
Here $W_{\rm PL}$
represents any of the PL Willmore energy discussed earlier, and $\mathcal{E}$
is the Dirichlet energy of the piecewise linear immersion (determined by $\cV$)
of the domain simplicial complex $K$ (determined by the face list $\cF$)
endowed with a suitable angle assignment.

Proposition~\ref{prop:Wmin_Regularized} %
suggests that, while it is not necessary, it would be best if we use a conformal structure
for which a conformal parametrization is possible, as it would free us from the need
of making the penalization $\lambda$ arbitrarily small, thus giving us a more efficient numerical method.
In the genus 0 case, there is only one conformal structure up to conformal equivalence. So we may simply
assume every triangle in $K$ to be equilateral. Under this angle assignment, the Dirichlet energy is, by \eqref{eq:cotangent_E_PL},
\bea \label{eq:EPL}
\mathcal{E} = \frac{1}{2} \cot(60^\circ) \sum_{e \in {\rm Edge}}  {\rm length}(e)^2=
\frac{1}{2\sqrt{3}}\sum_{e \in {\rm Edge}}  {\rm length}(e)^2.
\eea
As mentioned, we can still use \eqref{eq:EPL} in the higher genus cases, as long as we use a small enough $\lambda$;
see the genus 1 results in the next section.

\subsection{Computational results}
\begin{figure}[ht]
\centerline{
\begin{tabular}{cccc}
\includegraphics[width=1.4in]{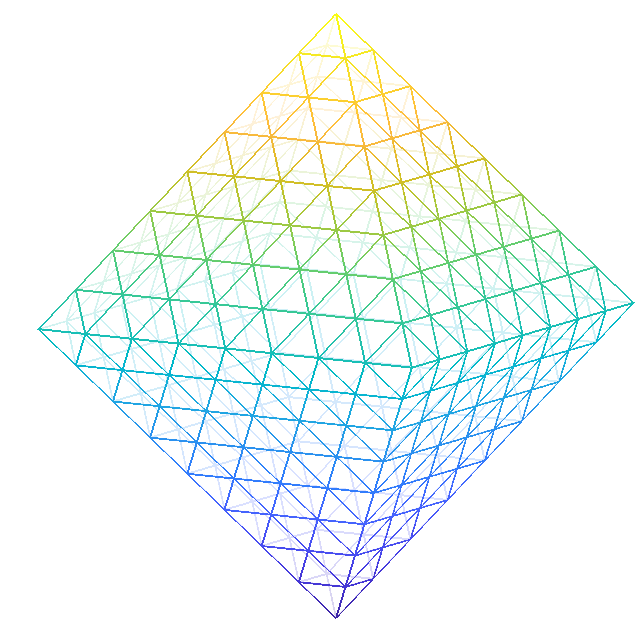} &
\includegraphics[width=1.4in]{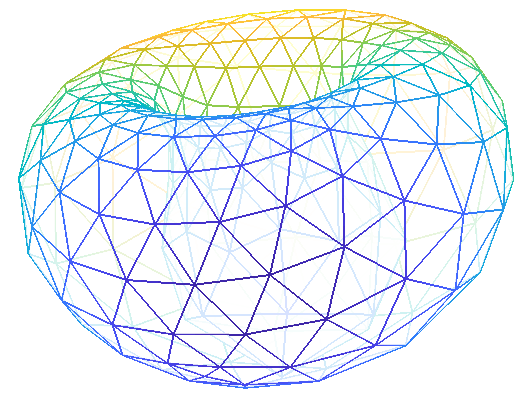} & \hspace{-.1in}
\includegraphics[width=1.4in]{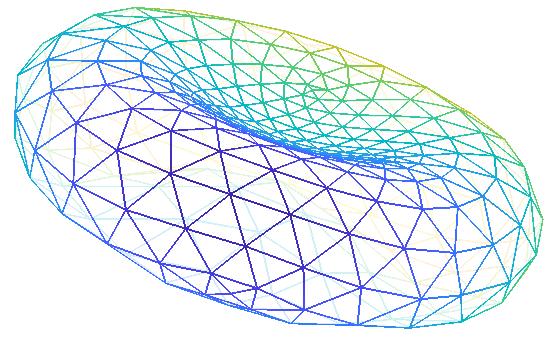}& \hspace{-.08in}
\includegraphics[width=1.4in]{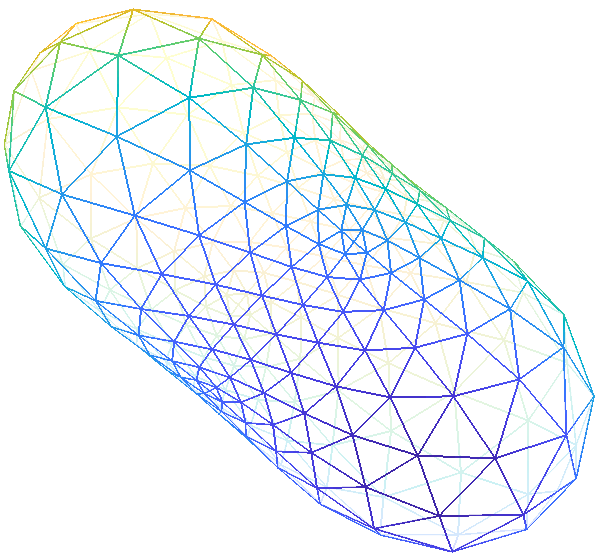} \\
(a) initial mesh & (b) $v_0= .50$, stomatocyte & (c) $v_0=.62$, oblate & (d) $v_0=.85$, prolate \\
$A=\cE=6.93$ & $A=6.93,\; \cE = 7.12$ & $A=6.93,\; \cE = 7.04$ & $A=6.93,\; \cE = 7.00$\\
& $W_{\rm NormalCur}=26.83$ & $W_{\rm NormalCur}=24.62$ & $W_{\rm NormalCur}=16.39$
 \end{tabular}}
\caption{Genus 0 $v_0$-constrained $W_{\rm NormalCur}$-minimizers}
\label{fig:Failure4}
\end{figure}

Figure~\ref{fig:Failure4} shows how
the $\cE$-penalized method \eqref{eq:PeanalizedWPLmin}
performs on the genus 0 Canham problem
with reduced volume/isoperimetric ratio
$v_0 := \frac{3V_0}{4\pi} (\frac{A_0}{4\pi})^{-\frac{3}{2}}$
in three different intervals known to give the shapes of a stomatocyte, discocyte (red-blood cell) and prolate.
The results are consistent with those reported in the biophysics literature, and those produced by
SS methods. See Section~\ref{sec:Further} for a more thorough comparison of different methods applied to these problems.

In
each case, we found that the numerical solution changes little with the choice of $\lambda$,
which is consistent to the conclusion of Proposition~\ref{prop:Wmin_Regularized}.
Also from the fact that $\cE$ is pretty close to the area $A$, we conclude from
Theorem~\ref{thm:ConformalEnergy_PL} that the PL embeddings are close to being angle preserving.

In the higher genus case not all conformal structures
are conformally equivalent. In the genus 1 case, the uniformization theorem implies that our solution surface can be mapped
conformally to a flat torus. However, there is a one-complex dimensional family of non-equivalent
conformal structures on a flat torus.

\begin{figure}[ht]
\centerline{
\begin{tabular}{cccc}
\includegraphics[width=1.5in]{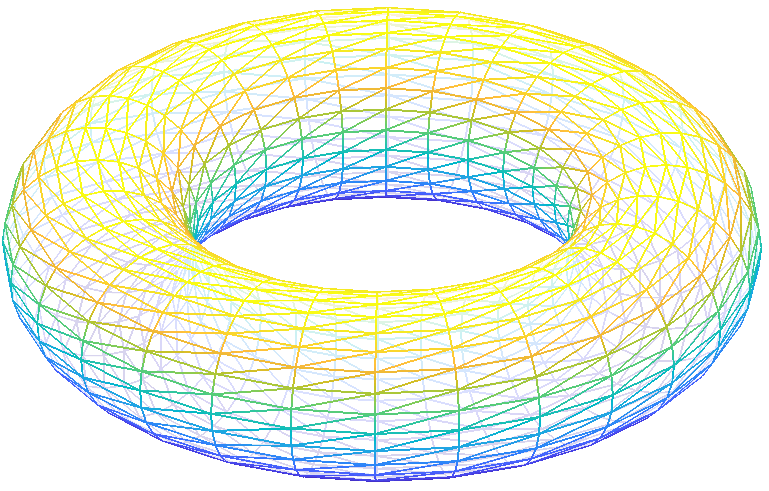} & \hspace{-.2in}
\includegraphics[width=1.5in]{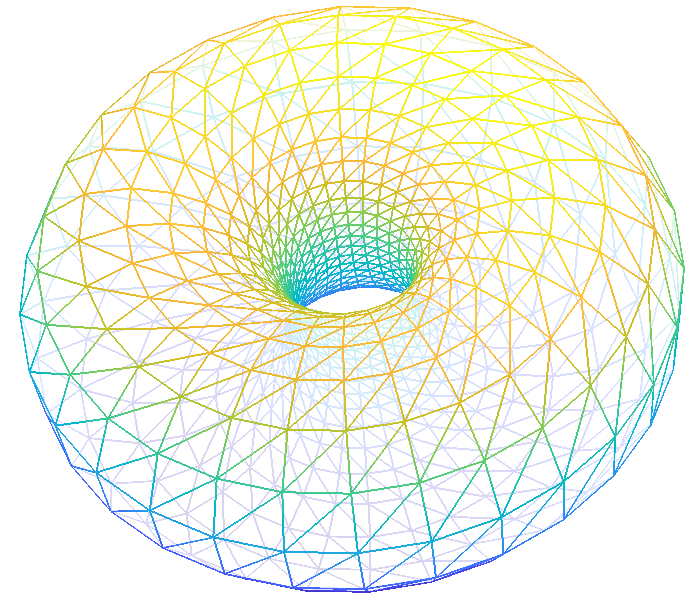} & \hspace{-.1in}
\includegraphics[width=1.5in]{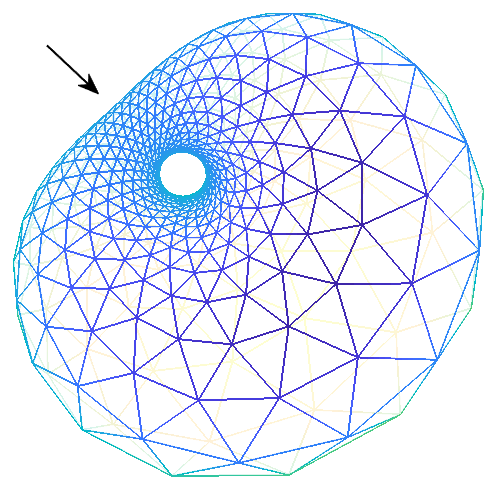} & \hspace{-.1in}
\includegraphics[width=1.5in]{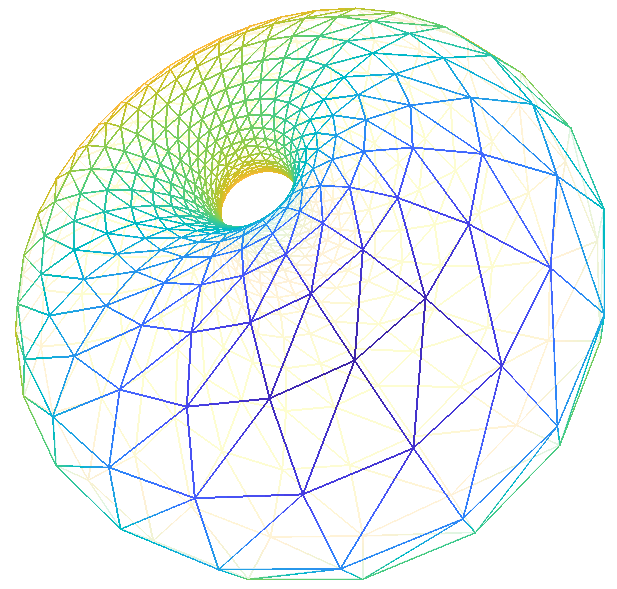}  \\
(a)  initial mesh, $A_0=117.57$  &
\hspace{-.1in} (b) $W_{\rm EffAreaCur}=20.41$  &
\hspace{-.1in} (c) $W_{\rm EffAreaCur}=20.71$  &
\hspace{-.1in} (d) $W_{\rm EffAreaCur}=20.62$
\\
$28\times28\times2$ triangles &
\hspace{-.1in} $A=A_0$, $\cE_{e^{i\pi/2}}=118.29$ &
\hspace{-.1in} $A=A_0$, $\cE_{e^{i\pi/3}}=121.50$ &
\hspace{-.1in} $A=A_0$, $\cE_{e^{i\pi/2}}=119.95$ \\
 \end{tabular}}
\caption{(b) \& (c): Unconstrained Genus 1 $(W_{\rm EffArealCur}+\frac{1}{2}\, \cE_\omega)$-minimizers
with $\cE_\omega$ computed based on the square
conformal structure and the $60^\circ$ conformal structure. (d): $v_0$-constrained genus 1
$(W_{\rm EffArealCur}+\frac{1}{2}\, \cE)$-minimizer with $\cE_\omega$ computed based on the square
conformal structure, $v_0=0.8$.}
\label{fig:CliffordTorus_CS}
\end{figure}

Here, we use $W_{\rm EffAreaCur}$ to solve the genus 1 Willmore problem problem. And we choose
the domain Riemann surface to be a regularly triangulated parallelogram with sides $1$ and $\omega$ on the complex plane.
We denote by $\cE_\omega$ -- easily computable according to \eqref{eq:cotangent_E_PL} --
the corresponding Dirichlet energy.
In this
case we know that the solution is the Clifford torus, with Willmore energy $2\pi^2$, and the ideal conformal structure is that of a square,
i.e. $\omega = e^{i\pi/2}$.
We compare the result with what we get by using the $60^\circ$ conformal structure, i.e. $\omega =  e^{i\pi/3}$.
See the results in Figure~\ref{fig:CliffordTorus_CS} (b) \& (c).
We see that by using the ideal conformal structure
 the PL Willmore energy is closer to the expected $2\pi^2 \approx 19.73$,
and $\cE - A$ is closer to zero. (With the incorrect conformal structure, we cannot expect
 to have a conformal parametrization, and hence we cannot expect $\cE - A \goto 0$ as the grid sizes grow.)
Also with the incorrect conformal structure the axis-symmetry is broken, the
surface looks
 like a M\"obius transformation of the surface of revolution Clifford torus (a Dupin cyclide),
but with an indentation, indicated by the arrow in Figure~\ref{fig:CliffordTorus_CS} (c).
The indentation suggests that the method is not capturing the shape accurately. However,
by reducing $\lambda$ from $1/2$ to $1/10$,
the indentation disappears and the resulted surface becomes a more accurate approximation to a Clifford torus. This is completely in line
with what Proposition~\ref{prop:Wmin_Regularized} says.

In Figure~\ref{fig:CliffordTorus_CS}(d), we show a numerical solution of
the genus 1 Canham problem with reduced volume $v_0 = 0.8$, it is believed that
the solution is a \textbf{unique} M\"obius transformation of the surface of revolution Clifford torus; see \cite{YuChen:Uniqueness}.
And we also know that the ideal
conformal structure to use is the square one. The solution is accurate even with a relatively big penalization factor. Again, it is
in agreement with Proposition~\ref{prop:Wmin_Regularized}.


%

\section{Further Computational Results and Symmetry Preservation} \label{sec:Further}
A group of biophysicists obtained a plethora of computational results for the Canham and Helfrich problems for
low genus $g$ and various constraint parameters $(v_0, m_0)$
\cite{MichaletBensimon:Genus2,PhysRevLett.72.168,LIPOWSKY:Nature,LIPOWSKY:Genus2,Seifert:Starfish,PhysRevE.52.6623}.
Many interesting observations are made
 about the uniqueness, non-uniqueness, symmetry and phase transition properties for these problems.
So far few of these observations have been justified mathematically.
In fact, the mere existence of a solution of the Canham problem, i.e. existence of
 a Willmore minimizer with a prescribed isoperimetric ratio and genus, is an unsolved
 problem in geometric analysis when the genus is larger than 0 \cite{KMR:Willmore}; for the genus 0 case,
 see  \cite{Schygulla:Willmore}.
The existence of solution of the Helfrich problem has not been addressed so far.
Moreover, in the biophysics literature, the numerical behavior of the optimization method is never addressed; the only
information we have is that Brakke's Surface Evolver was used, as in the experiments done in \cite{Kusner:Brakke}.

In this final technical section, we give a few comparisons of the different implementations of
the PL and SS methods presented earlier. We
must begin with a confession:

The algorithms developed in Section 2 and 4 only address how to discretize the variational problem
into a standard finite-dimensional constrained optimization problem of the form
$$\min_{x \in \bR^{N}} f(x) \quad \mbox{s.t.} \quad g_i(x) = 0,$$
and the theory in Section 3 sheds some light on what it means to the variational problem \emph{if}
we manage to solve the optimization problem. Needless to say, solving the latter problem itself
--- a high-dimensional, nonlinear, nonconvex, constrained optimization problem --- is a major
challenge and there are many algorithms and solvers available. In our experiments,
we use the following three solvers:
(i) \verb$fmincon$ in the Matlab optimization toolbox, (ii) SNOPT \cite{SNOPT:SIAM} and
(iii) GRANSO \cite{CMO:GRANSO}.

Given the complexity of these solvers, there are countless issues to explore and compare in conjunction with our
problems. 
We shall focus mainly on the approximation and symmetry properties of the optimization {\it problems} arising from the PL and SS methods, and will briefly touch upon
the existence issue of these optimization problems.
With regrets, we will not address the important question of the efficiency of solving these problems by various optimization alogirithms/solvers.

\subsection{Comparison I: PL vs SS} The original expectation is that a PL method is less accurate than the
higher order SS counterparts, only that we found in Section~\ref{sec:negative} and \ref{sec:Regularization}
that the non-conforming nature of PL methods require us to introduce regularization.
We are finally in a position to compare the accuracy of our PL and SS methods.

It is known from Schygulla \cite{Schygulla:Willmore} that  $W$-minimizer of genus 0 with any prescribed
isoperimetric ratio $v_0 \in (0,1)$ exists, and that the minimum Willmore energy is strictly less than $8\pi$.
In Figure~\ref{fig:Failure4}, we solve for these minimizers using a regularized PL method for
$v_0=0.5$, $0.62$ and $0.85$. The $v_0=0.5$ case is arguably the most difficult one as an `invagination'
develops in the vesicle; in this case the PL Willmore energy computed is bigger than $8\pi$, see
the caption of Figure~\ref{fig:Failure4}(b).
We now solve the same three problems but using our Loop SS method.
\begin{figure}[ht]
\centerline{
\begin{tabular}{cccc}
\includegraphics[width=1.6in]{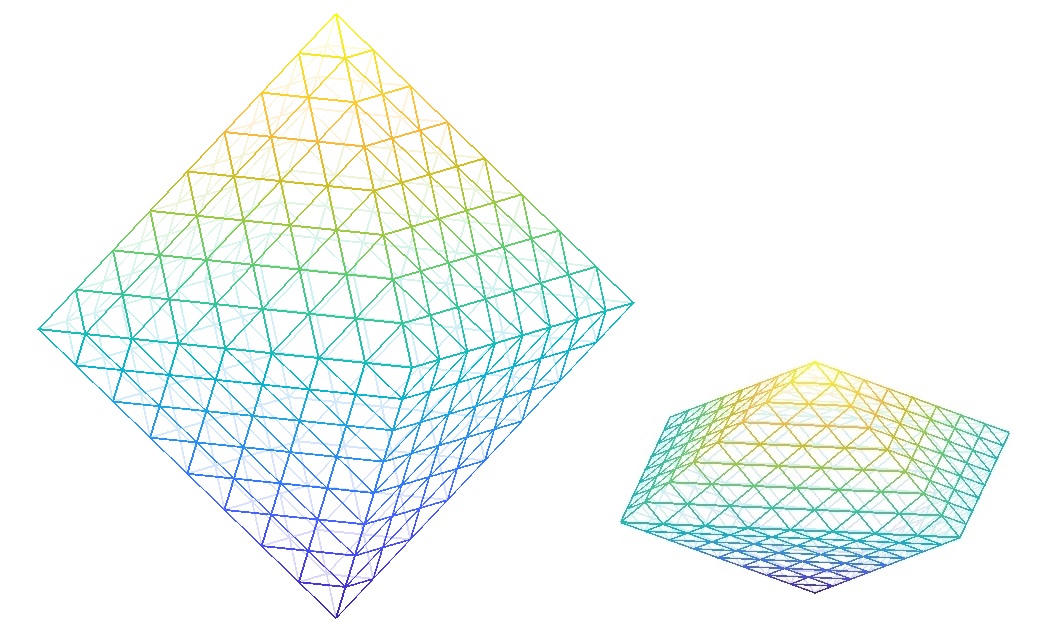} &
\includegraphics[width=1.2in]{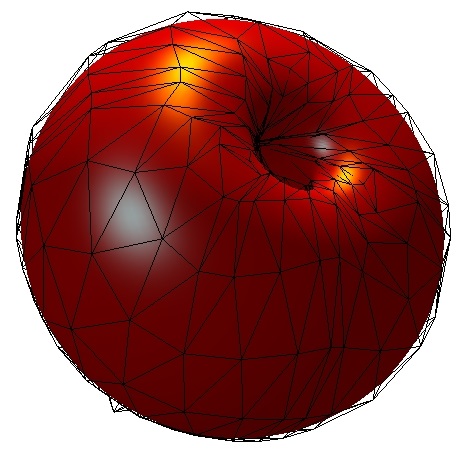} & \hspace{-.1in}
\includegraphics[width=1.6in]{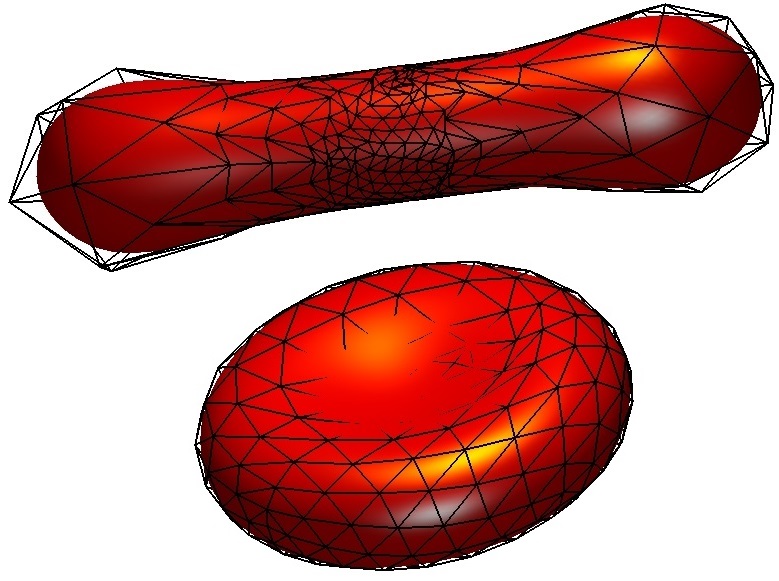}& \hspace{-.08in}
\includegraphics[width=1.4in]{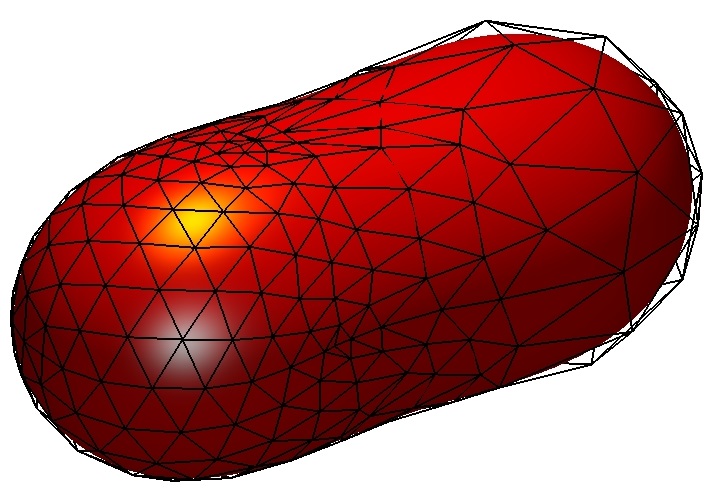} \\
(a) initial meshes  & (b) $v_0= .50$ & (c) $v_0=.62$  & (d) $v_0=.85$\\
                   &  $W=25.05$                            &   $W=24.66$ (top), &  $W=16.16$ \\
& & $W=23.92$ (bottom) &
 \end{tabular}}
\caption{Genus 0 $v_0$-constrained $W$-minimizers computed based on our Loop SS method. (a) two initial meshes: a
3-times subdivided octahedron and its flattened counterpart. Both initial meshes, and both the SNOPT and fmincon optimization solvers
give the stomatocyte in (b) and the prolate in (d).
For $v_0 = 0.62$,
the initial mesh with full octahedral symmetry happens to give the non-global local minimizer on the top of (c) (a prolate), whereas
the flattened initial mesh (with $\mathcal{D}_{4h}$ symmetry, in Schoenflies notation)
 gives the typical bi-concave red-blood cell shape (an oblate) at the bottom.
In (b)-(d), the black lines depict the control mesh of the SS surface, the red surface is the SS surface itself.}
\label{fig:SchygullaSpheres}
\end{figure}
Without an explicit representation of the solution, we cannot directly compare the accuracy of the PL and SS methods.
However, we see from our computation
that in each case, the (true) Willmore energy of the SS approximation is lower than the (PL) Willmore energy
of the PL approximation. In the $v_0=0.5$ case, the Willmore energy of the SS approximation is less than $8\pi$,
as it should according to Schygulla's result. These observations are consistent with the higher accuracy order of SS
compared to PL approximations.

\subsection{Comparison II: Symmetry preservation vs symmetry breaking}
\begin{figure}[ht]
\centerline{
\begin{tabular}{cccccc}
\includegraphics[width=1.0in]{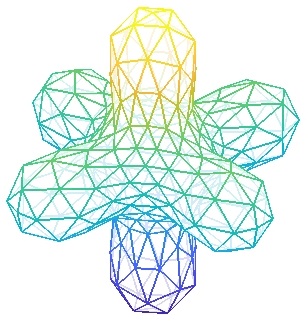} & \hspace{-.1in}
\includegraphics[width=1.0in]{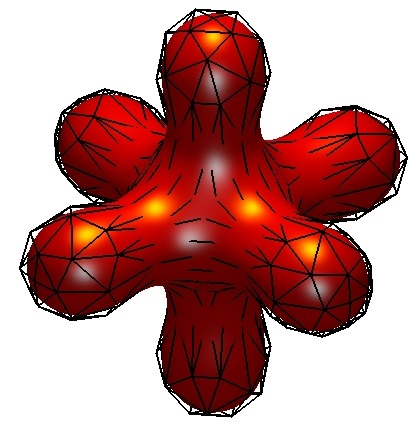} & \hspace{-.1in}
\includegraphics[width=1.0in]{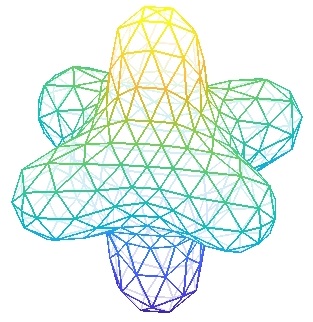}& \hspace{-.1in}
\includegraphics[width=1.0in]{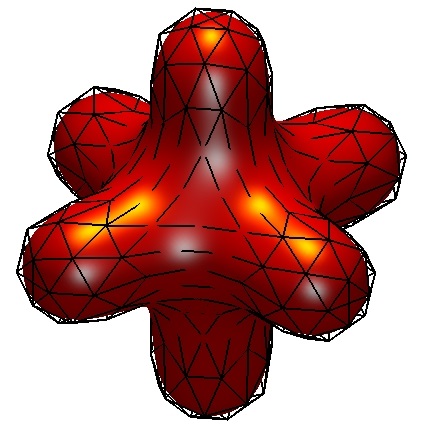} & \hspace{-.1in}
\includegraphics[width=1.0in]{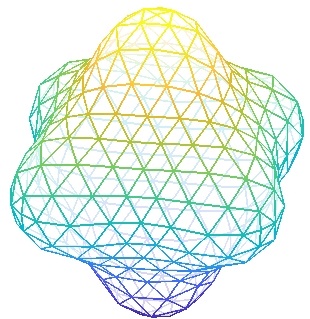}& \hspace{-.1in}
\includegraphics[width=1.0in]{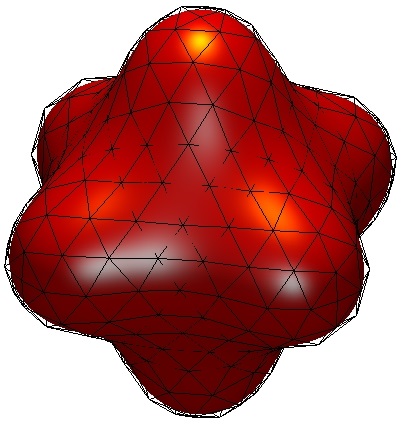} \\
(a) $W_{\rm NormalCur}$ & (b) $W=50.73$ & (c) $W_{\rm NormalCur}$ & (d) $W=41.40$ & (e) $W_{\rm NormalCur}$ & (f) $W=24.36$\\
$=51.97$ & & $=41.55$ & & $=24.04$ &
 \end{tabular}}
\caption{Genus 0 $v_0$-constrained local $W$-minimizers computed with PL and SS methods, the GRANSO solver and
the octahedral initial mesh in Figure~\ref{fig:Failure4}(a):
(a)-(b): $v_0=0.50$,
(c)-(d): $v_0=0.62$,
(e)-(f): $v_0=0.85$; (a),(c),(e) are based on
$W_{\rm NormalCur}$ and Dirichlet energy penalization with $\lambda=2$,
(b),(d),(f) are based on Loop SS. It is evident that {\it all} local minimizers are not global minimizers -- they
inherit the (incorrect) octahedral symmetry from the initial mesh.}
\label{fig:SchygullaSpheres_Fail}
\end{figure}

The results in Figure~\ref{fig:Failure4} and \ref{fig:SchygullaSpheres} are based on the \verb$fmincon$ and SNOPT solvers.
With our third solver, GRANSO, and the same octahedral initial mesh in Figure~\ref{fig:Failure4}(a),
we got totally different, non-global local minimizers which inherit the octahedral symmetry from the initial mesh;
see Figure~\ref{fig:SchygullaSpheres_Fail}.

Apparently, the underlying constrained optimization algorithm employed in the GRANO solver has a {\it symmetry preserving}
property. To the best of the authors' knowledge, this kind of
 symmetry preservation and breaking
properties is not well-addressed in the optimization literature. Here, we
 observe a fundamental difference among the three solvers: GRANSO is capable of preserving
symmetry in the sense below,\footnote{The authors of GRANSO designed the solver with the intention that it can handle non-smooth problems, but {\it without} the intention for symmetry preservation.
Our original use of it is due to the non-smoothness of $W_{\rm Bobenko}$; recall Section~\ref{sec:negative}.
The problems in this section are smooth and we are only exploiting its symmetry preserving property.}
while \verb$fmincon$ and SNOPT are capable of breaking symmetry.

\gap
\noindent
{\bf Symmetry preservation of gradient and Newton descent.}
To define `symmetry preservation' precisely, we begin with
the following fact about gradient flow which is well-known to geometers:
Let $M$ be a Riemannian manifold,
and $G$ be a group of isometries acting on $M$. If we have a smooth $G$-invariant functional
$F:M \goto \bR$, i.e. $F(g \cdot x)=F(x)$ for all
$x\in M$ and $g\in G$,
then the gradient flow map $\Phi(x,t)$ is also $G$-invariant:
$\Phi(g\cdot x,t)=g \cdot \Phi(x,t)$. In particular, if
$x_0$ is a symmetric point (i.e. $g\cdot x_0=x_0$ for all $g \in G$), then so is $\Phi(x_0,t)$ for any $t$.

Here we prove a  version of this fact tailored for our setting; the proof easily
extends to methods beyond gradient descent.

Any one of our geometric functionals $F =$ $W$, $A$, $V$ or $M$ is a $O(3)$-invariant
functional,
i.e.
\bea \label{eq:Sym1}
\forall g \in O(3), \quad
F(g (\cV) ) = F(\cV), \quad g (\mathcal{V}) := (g v_1,\dots,g v_{N}) \mbox{ for }
\cV = (v_1,\cdots,v_N).
\eea
There is yet another invariance, namely invariance under simplicial isomorphisms.
For a face list $\cF$ in a triangle mesh that specifies the simplicial complex structure of the mesh, there is a
subgroup, denoted by $S(\cF)$, of the permutation group of $1,\ldots, N$ that corresponds to the group of simplicial isomorphisms.
In other words, re-labelling the vertex indices according to the permutations in $S(\cF)$ gives the same triangulation (i.e. the geometric realization of
the simplicial complex.) Our geometric functionals must satisfy
\bea \label{eq:Sym2}
\forall \pi \in S(\cF), \quad
F(\pi (\cV) ) = F(\cV), \quad \pi( \cV) := (v_{\pi(1)},\ldots,v_{\pi(N)}).
\eea
\begin{definition} \label{def:Symmetry}
Let $G$ be a finite subgroup of $O(3)$.\footnote{The Schoenflies notation can be used to refer to any one of the possibilities of such a $G \triangleleft O(3)$.}
A simplicial complex, specified by a face list $\cF$ of a mesh, is said to \textbf{support (the symmetry group) $G$} if $G$ is isomorphic
 to some subgroup of $S(\cF)$.
In this case, we denote the correspondence by
$$
G \ni g \longleftrightarrow \pi_g \in S(\cF).
$$
A mesh $(\cV,\cF)$ in $\bR^3$ is called $G$-symmetric
if for every $g \in G$, $g(\cV) = \pi_g(\cV)$. Equivalently, if we define the group action of $G$ on the manifold $M= \bR^{N\times 3}$ by
$$g \cdot \cV := \pi_g^{-1} ( g ( \cV )), \quad\; g \in G,$$ then (with $\cF$ fixed) $\cV$ is $G$-symmetric if and only
$g \cdot \cV = \cV$ for all  $g\in G$.
\end{definition}
In this setting, the set of symmetric points is a linear submanifold of $\bR^{N\times 3}$.
Note that whenever the underlying simplicial complex supports the symmetry group $G$, then by
\eqref{eq:Sym1} and \eqref{eq:Sym2} $F$ is $G$-invariant in the sense that
\bea \label{eq:Ginvar}
F( g \cdot \cV) = F(\cV), \quad \forall g \in G.
\eea

A gradient flow
$\mathcal{V}(t)$ of $F$ satisfies $\dot{\cV}(t) = -\nabla F (\cV(t))$ and the corresponding
gradient descent algorithm satisfies
$$
\cV_{k+1} = \cV_k - \alpha_k \, \nabla F (\cV_k).
$$

By the chain rule applied to \eqref{eq:Ginvar}, $g^T \cdot \nabla F(g \cdot \cV) = \nabla F(\cV)$. Then, by orthogonality,
$\nabla F(g \cdot \cV) = g\cdot \nabla F(\cV)$.
We can see that a gradient descent, which we formally denote by ${\rm GD}(\cV, \alpha):= \cV - \alpha \, \nabla F (\cV)$,
satisfies the transformation property:
$$
{\rm GD}( g \cdot \cV, \alpha ) = g \cdot \cV - \alpha \overbrace{\nabla F( g\cdot \cV)}^{g\cdot \nabla F( \cV)} =
g \cdot {\rm GD}( \cV, \alpha ).
$$
In particular, gradient descent preserves the symmetry of the previous iterate.
It is also easy to check that the quasi-Newton BFGS method is symmetry preserving.

In our application, $F$ can be chosen to be the penalty function of a Helfrich problem:
$$
F(\cV; \mu) :=  W(\cV) + \frac{\mu}{2} \big\{ (A(\cV)-A_0)^2 +(V(\cV)-V_0)^2 + (M(\cV)-M_0)^2 \big\}, \quad \mu>0.
$$
Since $F$ inherits the $O(3)$-invariance from its constituent objective and constraint functions, any one of the
gradient descent, Newton or BGFS methods  applied to $F$ would furnish a symmetry preserving algorithm
for approximately solving the constrained optimization problem. Note that the symmetry preservation property holds regardless of the
penalization parameter $\mu$ or the line search parameters $\alpha$, which vary from iteration to iteration.
However, the penalization parameter $\mu$ has be big enough in order for the constraints to be approximately satisfied.

Sophisticated BFGS-SQP based
solvers such as SNOPT, GRANSO and \verb$fmincon$ are not simply based on applying the BFGS method
to the penalty function,
and those who engineer these solvers do not have the symmetry preserving property
in mind.
It is therefore interesting to see from the experiments that the BFGS-SQP method used in GRANSO preserves symmetry almost perfectly, even in the presence
of roundoff errors.
We shall present a formal justification of this observation in a separate report.

In principle, we can optimize in a symmetry preserving
way by optimizing only over the degrees of freedom that determine the control mesh up to the desired symmetry, as is done in
Brakke's Surface Evolver or \cite{ZHAO2017164}. This has the added advantage of reducing the
dimensionality of the problem, but requires an extra  effort in coding and algorithmic development.
GRANSO, or any solver with the same symmetry preserving property, frees us from the latter.


\subsection{Comparison III: Symmetric vs `Best' Minimizers} \label{sec:Pear}
The example in the previous section says the obvious:
applying a symmetry preserving optimization algorithm to an initial guess with the \emph{wrong} symmetry
is not going to solve the problem. What if we have the correct symmetry?
Our next experiment, based on comparing different solvers, shall reveal a rather subtle feature \emph{not} of the optimization algorithms,
but of the optimization problems themselves.

\begin{figure}[ht]
\centerline{
\begin{tabular}{cccc}
\begin{tabular}{c}
\includegraphics[width=.6in]{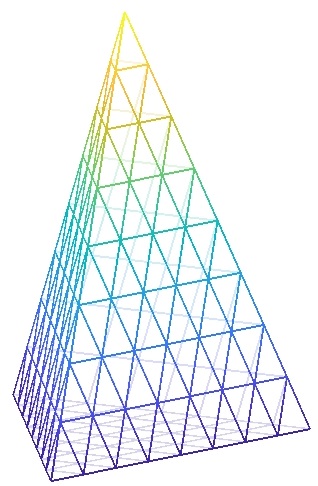} \\
initial mesh
\end{tabular}
& \hspace{-.3in}
\begin{tabular}{ccc}
C2g0 & \hspace{-.2in}Loop &\hspace{-.2in} PL \\
\includegraphics[width=.6in]{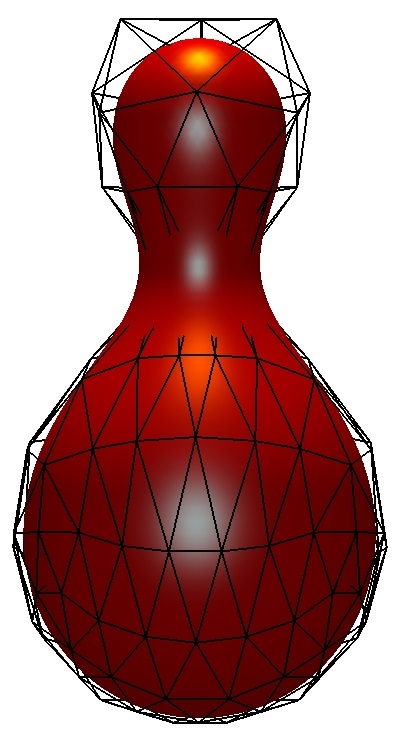} & \hspace{-.2in}
\includegraphics[width=.6in]{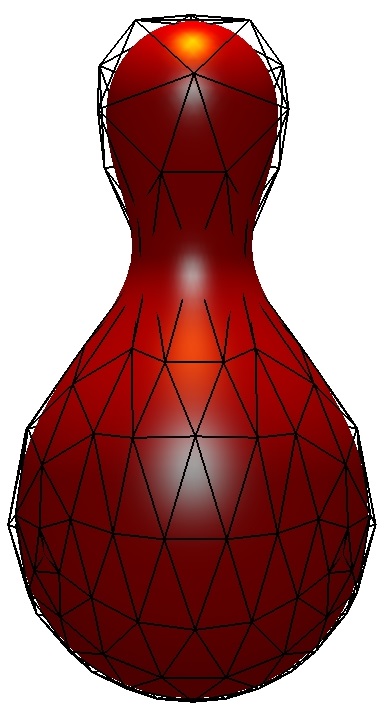}& \hspace{-.2in}
\includegraphics[width=.6in]{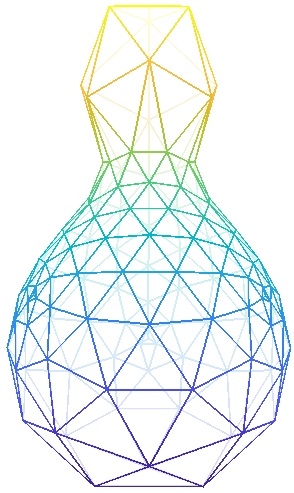} \\
 20.432 & \hspace{-.2in}  20.48 &\hspace{-.2in}  21.89 \\
\end{tabular}
& \hspace{-.3in}
\begin{tabular}{ccc}
C2g0 &\hspace{-.2in} Loop &\hspace{-.2in} PL \\
\includegraphics[width=.6in]{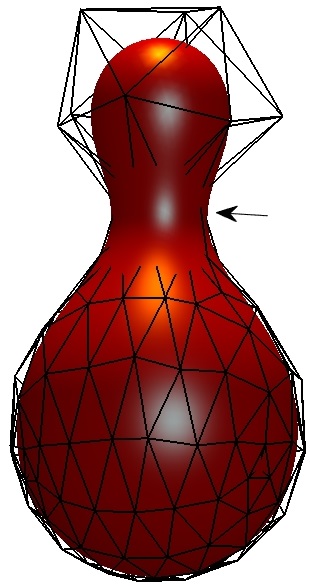}& \hspace{-.2in}
\includegraphics[width=.6in]{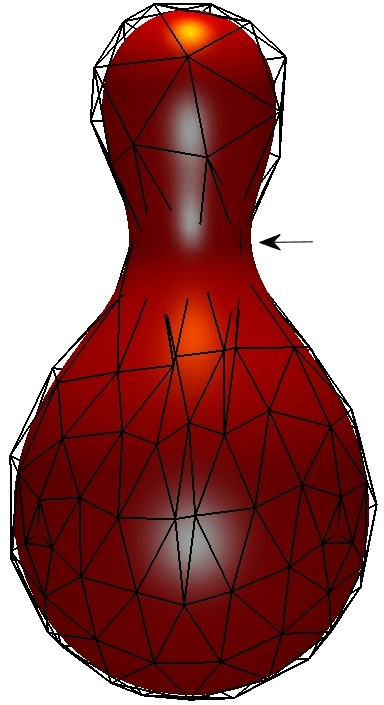} & \hspace{-.2in}
\includegraphics[width=.6in]{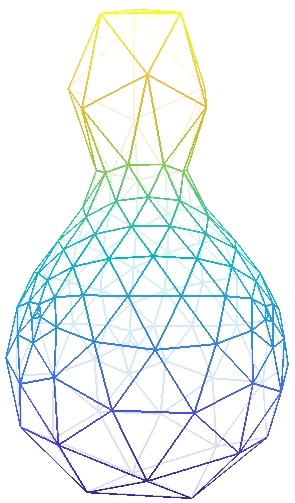}\\
20.427 & \hspace{-.2in}  20.47 &\hspace{-.2in} 21.87 \\
\end{tabular}
& \hspace{-.3in}
\begin{tabular}{ccc}
C2g0 &\hspace{-.2in} Loop &\hspace{-.2in} PL \\
\includegraphics[width=.6in]{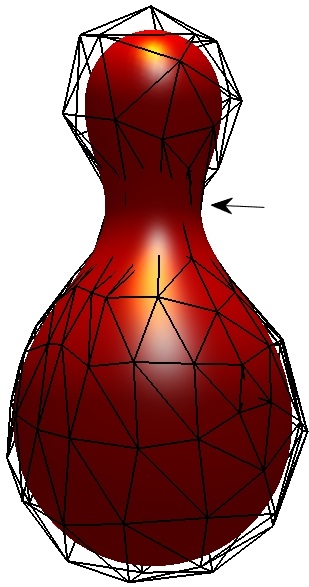} & \hspace{-.2in}
\includegraphics[width=.6in]{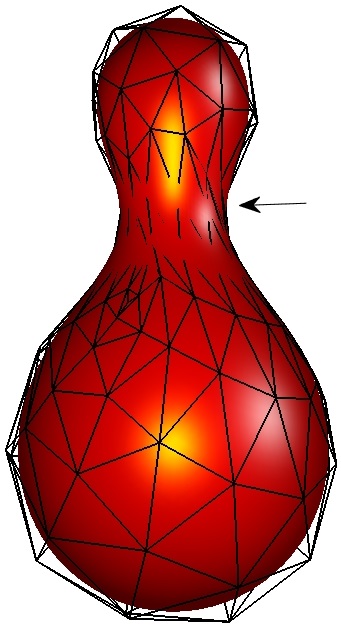}& \hspace{-.2in}
\includegraphics[width=.6in]{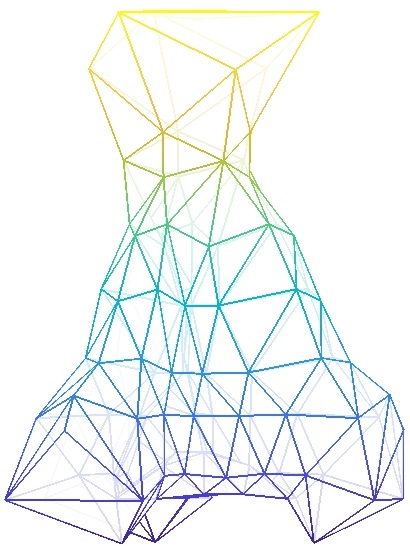}\\
20.426 & \hspace{-.2in} 20.44 &\hspace{-.2in} infeasible \\
\end{tabular}
\\
 &\hspace{-.3in} (i) GRANSO &\hspace{-.3in} (ii) SNOPT &\hspace{-.3in} (iii) fmincon
 \end{tabular}}
\caption{Numerical solutions of the Helfrich problem with $(v_0,m_0)=(0.8, 1.2)$. For the PL method, we use
$W_{\rm EffAreaCur}$ and $M_{\rm Steiner}$  and $\lambda=2$, our use of fmincon fails to give a feasible point, while
GRANSO and SNOPT work well.
GRANSO preserves the $\mathcal{D}_3$ symmetry of the initial mesh, while SNOPT and fmincon break the symmetry.
The numerical values are the corresponding (true or PL) Willmore energies. The symmetry-breaking cases
 give slightly lower Willmore energies than the symmetry-preserving counterparts.}
\label{fig:Pear}
\end{figure}

We consider the genus 0 Helfrich problem with $(v_0,m_0)=(0.8, 1.2)$. It is believed that the minimizer
is a surface of revolution with a pear shape. The two initial meshes
in Figure~\ref{fig:SchygullaSpheres}(a) have
extra symmetries (octahedral and $\mathcal{D}_{4h}$ symmetry) not possessed by a general surface of revolution.
Applying the GRANSO solver, which presumably preserves symmetry, would fail to yield the correct solution; it is indeed
what we observe from computation. If we use instead GRANSO with
a $\mathcal{D}_{n}$-symmetric initial mesh with no extra symmetry,
we expect to approach the correct minimizer with a $\mathcal{D}_{n}$-symmetric approximation,
and it is again what we observe from computation; see Figure~\ref{fig:Pear}(i).
In the case of the C2g0 and Loop SS methods, SNOPT and fmincon break the symmetry, as shown in
Figure~\ref{fig:Pear}(ii) and (iii).
If one zooms into the control meshes, one sees that the $\mathcal{D}_{3}$ symmetry is clearly broken around
the neck area (indicated by
the arrows in the relevant panels). Note that the neck area is also where the absolute Gauss curvature of the pear
 surface is the highest.

The asymmetric approximations to the (presumably symmetric) solution, produced by SNOPT and fmincon, have slightly
smaller Willmore energies than those of the symmetric approximations produced by GRANSO.
We 
believe that this is not caused by roundoff errors or truncation errors from
 numerical integration.

In (i), GRANSO terminates gracefully
with a stationarity condition satisfied up to a tolerance.
This suggests that
GRANSO produces a $\mathcal{D}_3$-symmetric critical point of the problem that is
a local minimizer of the  $(v_0,m_0)=(0.8,1.2)$ Helfrich problem with the added $\mathcal{D}_3$-symmetry constraint. 
%
 By the principle of symmetric criticality \cite{palais1979}, this critical point must also be
a \emph{critical point} of the full Helfrich problem
 (without any symmetry constraint.) Since (ii) and (iii) strongly indicate that this critical point is neither a local nor global minimizer,
 we are led to believe that it
is a \emph{saddle point} of the Helfrich problem.

 We may contrast this experiment with that in Figure~\ref{fig:NormalCurMinimizer58}(c)-(d).
We believe that the PL surface with an octahedral symmetry, as shown in
Figure~\ref{fig:NormalCurMinimizer58}(d),
is an absolute minimizer there.

\subsection{$W$-minimizers of genus $g$ and Lawson's $\xi_{g,1}$ surfaces} \label{sec:Lawson}
\begin{figure}[ht]
\centerline{
\begin{tabular}{cccccc}
\includegraphics[height=1.2in]{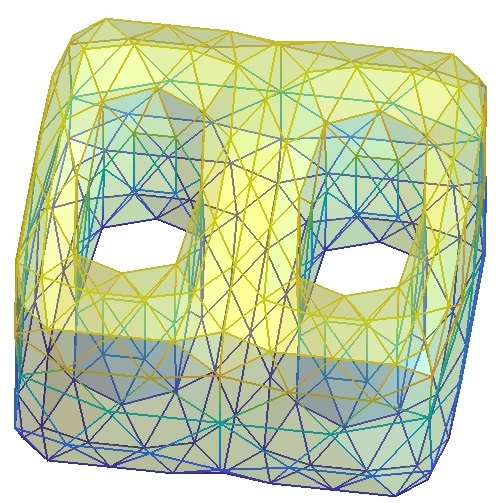} & \hspace{-.2in}
\includegraphics[height=1.2in]{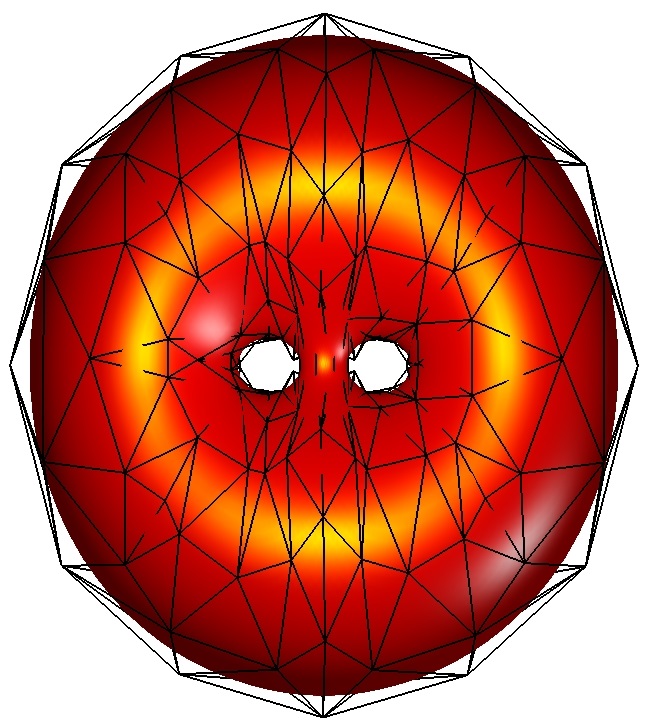} & \hspace{-.2in}
\includegraphics[height=1.2in]{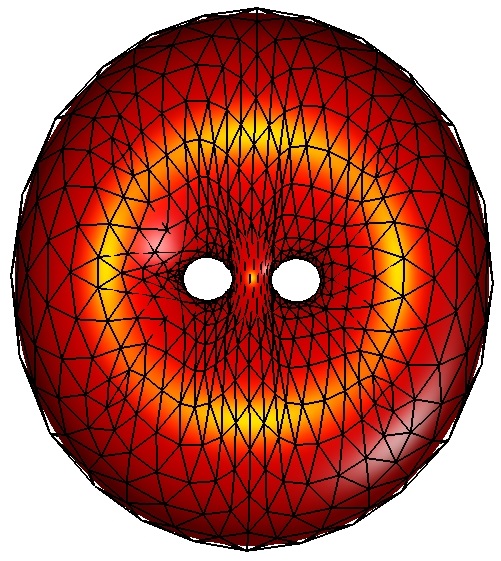}& \hspace{-.2in}
\includegraphics[height=1.2in]{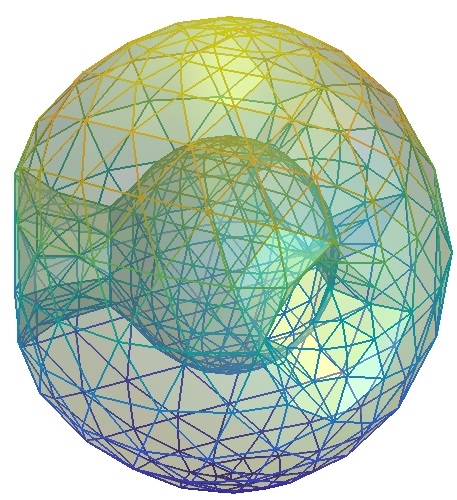} & \hspace{-.2in}
\includegraphics[height=1.2in]{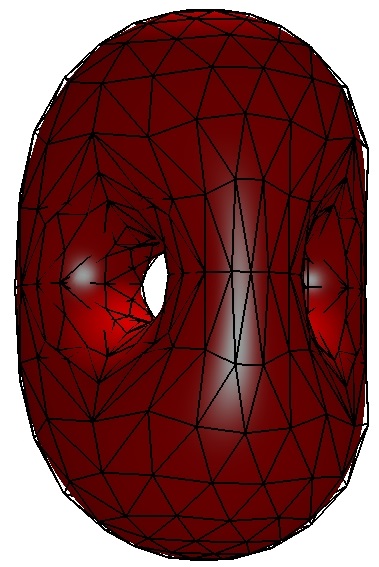}& \hspace{-.2in}
\includegraphics[height=1.2in]{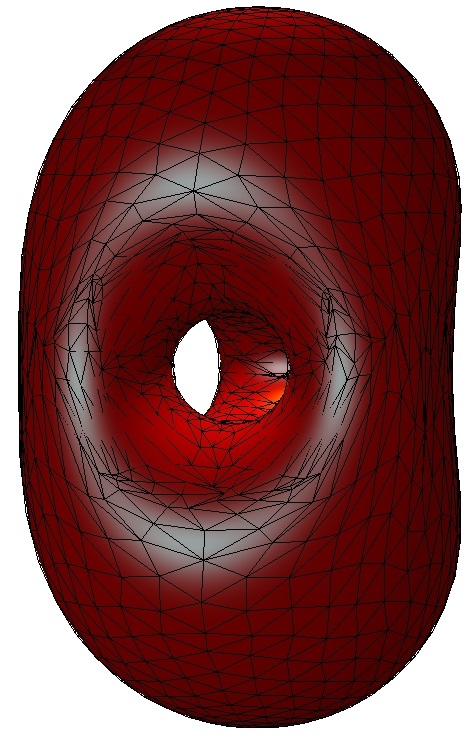} \\
(a) $\mathcal{D}_{2h}$ initial mesh &\hspace{-.2in} (b) $W=22.48$ & \hspace{-.2in} (c) $W=21.98$
&\hspace{-.2in} (d) $\mathcal{D}_{3h}$ initial mesh &\hspace{-.2in} (e) $W=22.17$ &\hspace{-.2in} (f) $W=21.97$ \\
$230$ vertices &\hspace{-.1in}$230\times 3$ d.o.f.  & $926\times 3$ d.o.f. & $430$ vertices & \hspace{-.1in}$430\times 3$ d.o.f. & $1726\times 3$ d.o.f.\\
 \end{tabular}}
\caption{Genus 2 Willmore minimizers with $\mathcal{D}_{2h}$ and $\mathcal{D}_{3h}$ symmetry. In each case, after an optimization is completed at
the original resolution (as shown in (a)-(b) and (d)-(e)), the resulting
control mesh is Loop subdivided once and used as the initial mesh for a second round of optimization.
The more accurate results at the finer resolutions are shown in (c) and (f), respectively. The GRANSO solver is used
so that the desired symmetries are preserved.}
\label{fig:Lawson}
\end{figure}

It is conjectured that the stereographic images of Lawson's
minimial surface $\xi_{g,1}$ in $\mathbb{S}^3$ \cite{Lawson:Minimal}
are the only $W$-minimizer of genus $g$ in $\bR^3$.
The term `Willmore conjecture' -- now the celebrated Marques-Neves theorem \cite{MarquesNeves:Willmore} --
refers to the $g=1$ case of this more general conjecture.
Here, we use our SS method to illustrate this conjecture for $g= 2$.
There are two stereographic projections of $\xi_{g,1}$ from $\mathbb{S}^3$ to $\bR^3$ that give
the resulting surfaces a $\mathcal{D}_{2h}$ or $\mathcal{D}_{3h}$ symmetry. We create initial control meshes
with these two symmetries and solve the Willmore problem using our Loop SS method and the symmetry-preserving
GRANSO solver. We also apply an iterative refinement, exploiting the underlying subdivision structure, to improve the accuracy
of the solutions. The resulting surfaces are visually the same as the corresponding
stereographic projections of Lawson's $\xi_{2,1}$, and have the same Willmore energy
of approximately 21.9.
(The two surfaces are meant to be M\"obius transformations of each other.)
See Figure~\ref{fig:Lawson}.

The same code for generating Figure~\ref{fig:Lawson}, available in the Wmincon package, works for any genus $g\geq 1$.
We have used it to empirically verify the generalized Willmore conjecture for up to genus $g=6$.

\section{Conclusion and Future Work}\label{sec:Future}
Admittedly,
our analyses in Section~\ref{sec:analysis}
only address one fundamental aspect of what we may expect from the PL and SS methods.
Even so, both the positive and negative results there have much room for improvements.

The positive result Theorem~\ref{prop:MainPositive} asserts that some
subsequence of a sequence of discrete minimizers converges to a
minimizer of the continuous problem; the result does not say anything about the rest of the sequence.
Does it mean that some elements of our hard-earned sequence $\mathbf{x}_j$
may have nothing whatsoever to do with
\emph{any} minimizer of
the continuous Willmore problem, even
for arbitrarily large $j$? We do not believe so,
as we may use Theorem~\ref{prop:MainPositive} repeatedly in the following way:
\begin{itemize}
\item Apply the theorem to extract a convergent subsequence (convergent in the sense stated in the theorem) from $\mathbf{x}_j$
\item Remove this subsequence from the original sequence, the remaining sequence is still a minimizing subsequence.
\item Apply the theorem again to extract a convergent subsequence from the remaining sequence.
\item ...
\end{itemize}
Loosely speaking, this process should exhaust the whole sequence, and hence it seems plausible we can partition the original sequence $\mathbf{x}_j$
into subsequences each of which converges to \textit{some} Willmore minimizer. We therefore expect the following to hold true:
\begin{conjecture} \label{prop:MainPositive2} \normalfont
Let $\bx_j$ be any minimizing sequence in the setting of Theorem~\ref{thm:Simon:Willmore}. (In particular, $\bx_j$ can be
an approximate $W$-minimizer over $\Imm_{\euS^j}$, as in Theorem~\ref{prop:MainPositive}.)
\begin{itemize}
\item[{(I)}]
There is a sequence of M\"obius transformations $G_j$ in $\bR^3$ such that the sequence $G_j \circ \mathbf{x}_j(K)$ of surfaces
can be partitioned into subsequences each of which
 converges in Hausdorff distance to \emph{some} M\"obius representative of a genus $g$ Willmore minimizer.

 \item[{(II)}] If we further assume that
the Willmore minimizer of genus $g$ is unique up to M\"obius transformations, as would be implied by the generalized
Willmore conjecture (see Section~\ref{sec:Lawson}), then the whole sequence of surfaces $G_j \circ \mathbf{x}_j(K)$ in (I)
converges in Hausdorff distance to \emph{some} M\"obius representative of a genus $g$ Willmore minimizer.
\end{itemize}
\end{conjecture}


The situation for the Canham and Helfrich problems is more challenging. On the geometric analysis side, even the existence problem for the Canham problem
 is not resolved
for every positive genus and all isoperimetric ratios; see \cite{Schygulla:Willmore} and \cite{KMR:Willmore}. (In the genus zero case,
Schygulla \cite{Schygulla:Willmore} solves the existence
problem for all isoperimetric ratios.)
On the numerical analysis side, a key difficulty is to prove the corresponding density
result Corollary~\ref{Cor:Density} with a fixed isoperimetric ratio constraint, i.e.
\begin{conjecture}\normalfont
For any $v \in (0,1),$
$\bigcup_j \Imm_{\euS^j}^v$ is dense in $\ImmX^v(K)$,
where
$\Imm_{\euS^j}^v$ is the space of all immersed Loop subdivision surfaces over a closed oriented $K$ with isoperimetric ratio
$v$ and $\Imm_{W^{2,2}\cap C^1}^v(K)$ is the space of elements in $\Imm_{W^{2,2}\cap C^1}(K)$ with isoperimetric ratio $v$.
\end{conjecture}
 If this difficulty can be overcome, we expect that
a result similar to Theorem~\ref{prop:MainPositive} can be obtained for the Canham problem based on the existence results for
the Canham problem established in \cite{Schygulla:Willmore} and \cite{KMR:Willmore}. The Helfrich problem is out of reach for now.

The negative results Proposition~\ref{prop:ST} and ~\ref{prop:W_negative} are established on a case by case basis.
Perhaps there is a universal negative
result asserting that {\it any} consistent PL Willmore energy would fail in a similar manner.
Besides the very special genus 0 cases for $W_{\rm NormalCur}$ and $W_{\rm Bobenko}$
(recall Figure~\ref{fig:NormalCurMinimizer58}),
it appears that the space of PL surfaces (of a fixed combinatorial type) is `too big' that minimizing any PL $W$-energy
over it would always take us to some PL surface inconsistent with any smooth surface.
It is an open question to formulate and prove this speculation.

A deeper mathematical study of these issues and the proposed regularization in Section~\ref{sec:Regularization}
 may lead to understandings of other non-conforming methods, such as the ones proposed
in \cite{Bonito:Biomembranes,Schumacher:WillmoreFlow}. However, as mentioned in the introduction, the methods in
\cite{Bonito:Biomembranes,Schumacher:WillmoreFlow}
are in the spirit of `discretizing a minimization,' whereas the methods studied in this article are in the
spirit of `minimizing a discretization.' It is an open question to see if
 these different methods can be analyzed in a coherent way.

It will be interesting to address rate of convergence issues in the future.
Arden's result \cite[Theorem 2]{arden:2001} establishes a rate of convergence of Loop subdivision functions in $W^{2,2}$,
but we believe that the rates established
there are suboptimal. Also, the approximation rates of other schemes (e.g. \cite{ChenGrundelYu:Sphere}) are yet to be explored.
Regardless, further techniques are needed for transferring any approximation rate result in the pure approximation setting to
an approximation rate result
for the geometric variational problems. To this end, we expect knowledge on the second variation of the Willmore energy,
such as the results by Weiner \cite{Weiner:Willmore}, will be helpful.

\end{document}